\documentclass[a4paper, 9pt]{amsart}     % onecolumn (second format)
\usepackage{graphicx}
\usepackage{amsmath, amssymb, amsfonts}
\usepackage{enumitem}
\usepackage{hyperref}
\usepackage[shortalphabetic]{amsrefs}
\usepackage{mathptmx}    
\newtheorem{conditions}{Conditions}
\newtheorem{theorem}{Theorem}
\newtheorem{lemma}{Lemma}
\newtheorem{corollary}{Corollary}
\newtheorem{proposition}{Proposition}
\newtheorem{conjecture}{Conjecture}
\theoremstyle{definition}

\newtheorem{remark}{Remark}
\newtheorem{example}{Example}
\newcommand{\RR}{{\mathbb R}}
\newcommand{\NN}{{\mathbb N}}
\newcommand{\Sn}{{\mathbb S}}
\newcommand{\rr}{{\frak r}}
\newcommand{\ep}{\varepsilon}
\newcommand{\E}{\mathrm{e}}
\newcommand{\Haus}{d_{\mathcal{H}}}
\newcommand{\dHaus}{{\hat d}_{\mathcal{H}}}
\newcommand{\smartqed}{\qed}
\begin{document}

\title{Contracting convex hypersurfaces by curvature }

\author{Ben Andrews}
\address{
              Mathematical Sciences Institute, 
              Australian National University, ACT 0200 Australia, and\\
              Mathematical Sciences Center, Tsinghua University, Beijing 100084 China, and \\
              Morningside Center for Mathematics, China Academy of Sciences, Beijing 100190 China}
               \email{Ben.Andrews@anu.edu.au}
               \thanks{This research was supported by Discovery Grants DP0556211 and DP0985802 of the Australian Research Council.  The second author's research was also assisted by a University of Wollongong URC Small Grant entitled \emph{Curvature contraction of nonsmooth, weakly convex hypersurfaces into spheres}.}
\author{James McCoy}
\address{Institute for Mathematics and its Applications, 
School of Mathematics and Applied Statistics, 
University of Wollongong,
Wollongong, NSW 2522, 
Australia}\email{jamesm@uow.edu.au}
\author{
Yu Zheng}
\address{
Department of Mathematics,
East China Normal University,
Shanghai 200062,
China}
\email{zhyu@math.ecnu.edu.cn}


\begin{abstract}
We consider compact convex hypersurfaces contracting by functions of their curvature.  Under the mean curvature flow, uniformly convex smooth initial hypersurfaces evolve to remain smooth and uniformly convex, and contract to points after finite time.  The same holds if the initial data is only weakly convex or non-smooth, and the limiting shape at the final time is spherical.  We provide a surprisingly large family of flows for which such results fail, by a variety of mechanisms:   Uniformly convex hypersurfaces may become non-convex, and smooth ones may develop curvature singularities; even where this does not occur, non-uniformly convex regions and singular parts in the initial hypersurface may persist, including flat sides, ridges of infinite curvature, or `cylindrical' regions where some of the principal curvatures vanish;  such cylindrical regions may persist even if the speed is positive, and in such cases the hypersurface may even collapse to a line segment or higher-dimensional disc rather than to a point.  We provide sufficient conditions for these various disasters to occur, and by avoiding these arrive at a class of flows for which arbitrary weakly convex initial hypersurfaces immediately become smooth and uniformly convex and contract to points.
\end{abstract}

\maketitle

%% INTRO
\tableofcontents

%%% To do:   pictures:  proofs of displacement upper and lower bounds, persistence of flat sides
%%%			      construction of shrinking cylinders
%%%                                construction of shrinking ridges/discs
%%%                anisotropy?
%%%                expansion?
%%%                convex spacelike in Minkowski?

\section{Introduction}
\label{S:intro}
In this paper we study hypersurfaces evolving accoding to parabolic flows in which the speed of motion is governed by the curvature.  The best known example of such a flow is the mean curvature flow, which was studied in an influential paper of Gerhard Huisken \cite{Hu1}.  Given a smooth, uniformly convex initial hypersurface, Huisken proved existence of a unique solution of the mean curvature flow, which contracts to a point in a finite time.  He also proved that the solution becomes asymptotically spherical, so that rescaling about the final point by a suitable time-dependent factor gives smooth convergence to the unit sphere.  This behaviour remains true even when the initial hypersurface is non-smooth or only weakly convex, for example when it is given as the boundary of an arbitrary open bounded convex region.  A variety of other related flows have been considered:
Kaiseng Chou \cite{T} considered the flow by Gauss curvature, and proved that smooth, uniformly convex initial hypersurfaces contract to points in finite time.  
Ben Chow extended Chou's method to prove that convex hypersurfaces contract to points under speeds given by positive powers of the Gauss curvature, and adapted Huisken's arguments to prove that the
limiting shape is spherical under motion by the $n$th root of the Gauss curvature \cite{Ch1} and later by the square root of the scalar curvature \cite{Ch2}.  In the latter case he required an additional curvature condition on the initial hypersurface.  The first author generalised these results to a natural class of curvature flows \cite{A1}, with speeds homogeneous of degree one in the principal curvatures of the evolving hypersurface and satisfying some other natural conditions.  The class for which such results are known was widened further in \cites{A3,A4} which, in particular, removed the curvature condition on the initial hypersurface in \cite{Ch2}.   The
weaker statement, that smooth and uniformly convex initial hypersurfaces give rise to solutions which contract to points in finite time, is known for a wider class of flows \cite{Han}.  

In this paper we consider a wide class of such flows, and investigate the factors which determine whether the general picture of behaviour observed in the mean curvature flow remains valid.  

After the introductory sections \ref{S:result} and \ref{sec:soln}
which introduce our notations and definitions and discuss elementary properties of the flows and the kinds of solution we consider, we begin with two sections which produce counterexamples to the smooth contraction of uniformly convex hypersurfaces to points:

In section \ref{sec:lossreg} we find flows for which uniformly convex, smooth initial hypersurfaces can evolve to develop curvature singularities while the inradius remains positive.  For this we use a parametrization of convex hypersurfaces using the inverse of the Gauss map, in which the hypersurface is determined by its support function.  A hypersurface can be non-smooth but uniformly convex while its support function remains smooth, and the solutions we construct display exactly this behaviour.  The class of flows where smoothness is lost in this manner is surprisingly large, and contains many flows which seem both simple and natural.  We also see that the flows in which the speed is a homogeneous degree one function of the principal curvatures are easier to control, so there are relatively few flows of this kind which lose smoothness.

In section \ref{sec:lossconv}, we use a similar construction to
provide a family of flows for which there are weakly convex smooth hypersurfaces which evolve to become non-convex.  Under our assumptions this also guarantees that there are uniformly
convex smooth hypersurfaces which evolve to become non-convex.  For flows which are uniformly parabolic at cylindrical points, our condition is both necessary and sufficient for this behaviour to occur.

Having developed some understanding of how smoothness or convexity may be lost, we find a large class of flows for which initial smoothness and strict convexity is maintained until the hypersurface shrinks to a point (in some cases we can only prove the weaker statement that the inradius shrinks to zero).  We carry this out in section \ref{sec:topoints}.  Our main result (Theorem \ref{thm:contract.point}) is sharp in the class of flows for which the function $f_*$ which defines the speed of motion as a function of the principal radii of curvature is concave and extends beyond the boundary of the positive cone except where it is zero.  We also give some more restrictive results for concave speeds and for the case $n=2$ without concavity assumptions.

For flows which do not either lose smoothness or lose convexity (that is, those covered by the results of section \ref{sec:topoints}), we prove an existence and uniqueness statement for barrier solutions (a kind of generalized solution which is natural for the flows we consider), with initial data given by the boundary of any open, bounded convex region.  In the remaining sections of the paper we consider the question of which flows have the property that all such generalized solutions immediately become smooth and uniformly convex.  

Some examples where this fails are already known, beginning with the example of Hamilton \cite{Ha3} who showed that a flat side of a convex surface (that is, an open subset of the surface which is contained in a plane) persists for some time under the Gauss curvature flow.  Together with Daskalopoulos he investigated this phenomenon in some detail \cite{DH}.  The results of \cite{A8} imply that solutions eventually become smooth and uniformly convex.  Caputo and Daskalopoulos \cite{CD} considered flow by the harmonic mean of principal curvatures for surfaces which are initially uniformly convex with the exception of one flat side.  They showed that under this flow, the boundary of the flat side evolves by the curve shortening flow, so again flat sides persist for some time.   With Natasa Sesum they extended this to hypersurfaces moving with speed given by a quotient of successive elementary symmetric functions \cite{CDS}.

Previous examples where there is an affirmative answer to our question have arisen in \cite{A4}, where such a result was proved for
flows by powers $\alpha\leq \frac{1}{n}$ of the Gauss curvature, and in \cite{A7}, which includes strict convexity estimates for a class of `mixed discriminant' curvature flows.  A related result was proved by Dieter \cite{D}, who considered curvature flows with smooth, weakly convex initial data, where the speeds were quotients of successive elementary symmetric functions of the principal curvatures.  She required that for the initial hypersurface, the denominator of the quotient is positive, and proved that the hypersurface becomes immediately uniformly convex under the evolution.  

In section \ref{S:moving} we give necessary and sufficient conditions determining whether flat sides persist (the result is a dichotomy between flows where flat sides persist and those for which there is a lower bound on the displacement for positive times, depending only on the initial inradius and circumradius).  In particular we prove lower displacement bounds for flows with speed homogeneous of degree $\alpha<1$ in the principal curvatures, and prove that flat sides always persist for speed homogeneous of degree $\alpha>1$.  In the critical case $\alpha=1$ the two types of behaviour are distinguished by a single easily checked condition.  In sections \ref{S:boundspeed} and \ref{S:posspeed} we translate these displacement bounds into upper and lower bounds on the speed for positive times.  In particular the upper bounds on displacement and speed hold for essentially arbitrary flows.  We apply this observation to prove the main existence and unqueness result for barrier solutions in section \ref{sec:barrier}.

Sections \ref{S:cylinders} and \ref{sec:nonsmooth} provide examples where unique barrier solutions exist, but are sometimes not smooth or not strictly convex.

The examples discussed in section \ref{S:cylinders} display a new kind of behaviour for such flows:  Cylindrical regions can persist even in flows where flat sides do not (and so where the speed necessarily becomes strictly positive).  
We construct flows in which cylindrical regions in the evolving hypersurfaces have strictly positive speed, but do not become uniformly convex.  Instead these persist as cylindrical regions while shrinking in radius.  This has not been observed previously and shows the necessity of further conditions in order to deduce strict convexity and smoothness of solutions.  Remarkably, this construction also produces examples of weakly convex hypersurfaces which do not even contract to points --- for example under flow with speed equal to $|A|$, the norm of the second fundamental form, there are smooth weakly convex $n$-dimensional initial hypersurfaces which collapse to a $k$-dimensional disk of positive radius, for any $k\in \{1,\dots,n-1\}$.
The examples in section \ref{sec:nonsmooth} are similar in spirit, but in this case involve regions in the Gauss image where the minimum radius of curvature is zero, which persist while shrinking under the flow.  These correspond to ridges of infinite curvature which persist in the evolving hypersurface as it shrinks.

In the final section, we conclude by producing classes of flows where barrier solutions do indeed become smooth and strictly convex for positive times.  Informed by the conditions in the previous sections which allow the construction of counterexamples, we produce as wide a class as possible, though our results in this section are probably not exhaustive.   We state a conjecture which would amount to an essentially sharp result, but which seems to be beyond the reach of our methods.

The second and third authors are grateful to the Centre for Mathematics and its Applications at the Australian National University for its hospitality during their visits, when some of this work was completed. 

\section{Notation and preliminary results} \label{S:result}

In this paper we consider the behaviour of solutions $X : M^{n} \times \left[ 0, T \right)\to\RR^{n+1}$ of curvature flow equations of the form
\begin{equation} \label{E:theflow}
  \frac{\partial X}{\partial t}\left( x, t\right) = - S(A(x,t),g(x,t)) \nu\left( x, t\right) \mbox{,}
\end{equation}
where $X$ is a smooth family of immersions (that is, $X$ is smooth, and for each $t$ the map $X_t=X(.,t)$ is an immersion), $\nu(x,t)$ is the (outer) unit normal to the image hypersurface $M_t=X_t(M,t)$ at $X(x,t)$, $A(x,t)$ is the second fundamental form of $M_t$, and $g(x,t)$ is the induced metric.  Thus $A(x,t)$ and $g(x,t)$ are symmetric bilinear forms acting on the tangent space $T_xM$.  We require the flow to be geometrically invariant and smooth, which means that the function $S$ is a smooth function of the components of $A$ and $g$ which is invariant under change of basis, i.e. for any $L\in GL(n)$ we have $S(L^T A L,L^T g L)=S(A,g)$.  We also insist that $S$ be homogeneous of some positive degree $\alpha$ in the first argument, so that $S(\lambda A,g)=\lambda^\alpha S(A,g)$ for any $\lambda>0$ (this amounts to the requirement that the flow \eqref{E:theflow} has an invariance under spatial scalings).  

The invariance of $S$ implies that $S(A,g) = S(g^{-1/2}Ag^{-1/2},I)$ where $I$ is the identity matrix.  Thus $S$ is reduced to an $O(n)$-invariant function of one argument.  Henceforward when we write $S(A)$ we will mean $S(A,I)$.  For convenience we write $S=F^\alpha$. 
It follows from a theorem of Schwarz \cite{Schw} that $F(A,g) = f(\kappa(A,g))$,
where $f$ is a smooth symmetric function of $n$ variables, and $\kappa=(\kappa_1,\dots,\kappa_n)$ is the map which takes $A$ and $g$ to the vector of eigenvalues of $A$ with respect to $g$ (the \emph{principal curvatures}, defined only up to order).  Conversely, any such smooth symmetric function $f$ gives rise to a smooth $GL(n)$-invariant function $F$ (by a theorem of Glaeser \cite{Glaeser}).
There are several other conditions we will always impose on the function $f$:

\begin{conditions} \label{T:Fconds} 
\mbox{}
\begin{enumerate}[label={(\roman*)}, ref={(\roman*)}]
  \item\label{condi} $f$ is a smooth symmetric function defined on the positive cone 
  $$\Gamma_+ = \left\{ \kappa= \left( \kappa_{1}, \ldots, \kappa_{n} \right) \in \mathbb{R}^{n}: \kappa_{i} >0 \mbox{ for all } i = 1, 2, \ldots, n\right\} \mbox{.}$$
  \item\label{condii} $f$ is strictly increasing in each argument: $\dot f^i=\frac{\partial f}{\partial \kappa_{i}} > 0$ for each $i=1, \ldots, n$ at every point of $\Gamma_+$.
  \item\label{condiii} $f$ is homogeneous of degree $1$: $f\left( \lambda \kappa \right) = \lambda f\left( \kappa\right)$ for any $\lambda>0$ and $\kappa\in\Gamma_+$.
  \item\label{condiv} $f$ is strictly positive on $\Gamma_+$ and normalised:  $f\left( 1, \ldots, 1\right) = 1$.
  \end{enumerate}
\end{conditions}

Note that condition \ref{condiv} is merely a normalization, since the positivity follows from \ref{condi}--\ref{condiii}.
Several additional conditions will be used as need dictates:

\begin{conditions}\label{cond:extra}
\mbox{}
\begin{enumerate}[label={(\roman*)}, ref={(\roman*)}, start=5]
\item\label{cond:conv} $f$ is concave on $\Gamma_+$.
\item\label{condv} $f$ is \emph{inverse concave}.  That is, the function 
      $$f_*\left( x_{1}, \ldots, x_{n} \right) = f\left( x_{1}^{-1}, \ldots, x_{n}^{-1} \right)^{-1}$$
       is concave on $\Gamma_+$.
%\item\label{cond:vanishbdry} $f$ vanishes on the boundary of $\Gamma_+$.
\item\label{cond:invvanishbdry} $f_*$ vanishes on the boundary of $\Gamma_+$.
\item\label{cond:invconvbdry} The restriction of $f$ to the boundary of the positive cone is inverse-concave.  That is, if $f(1,\dots,1,0)>0$ then the function $\tilde f(x_1,\dots,x_{n-1})=f(x_1,\dots,x_{n-1},0)$ is inverse-concave on the $(n-1)$-dimensional positive cone $\Gamma_+^{(n-1)}$.
\item\label{cond:invconvinvbdry} The restriction of $f_*$ to the boundary of the positive cone is inverse-concave.
\item\label{cond:extendf} If $\kappa\in\partial\Gamma_+$ with $f(\kappa)>0$ then $\dot f^i(\kappa)>0$ for each $i$ and $f$ is smooth near $\kappa$.
\item\label{cond:extendf*} If $\rr\in\partial\Gamma_+$ with $f_*(\kappa)>0$ then $\dot f_*^i(\rr)>0$ for each $i$ and $f_*$ is smooth near $\rr$.
\item\label{cond:admit.reg} The flow \eqref{E:theflow} admits second derivative H\"older estimates in $\Gamma_+$ (see \ref{rmk:sec.hold.est} below).
\item\label{cond:admit.reg.extend} The flow \eqref{E:theflow} admits second derivative H\"older estimates in $\bar\Gamma_+\cap\{f>0\}$.
\end{enumerate}
\end{conditions}

\begin{remark}\label{rmk.conds}
\mbox{}
\begin{enumerate}[label={(\arabic*).}, ref={(\arabic*)}]
\item The definition of inverse-concavity here is slightly different from that in \cite{A3}, but is equivalent.
Previous results on hypersurface flows have always assumed one of the concavity conditions \ref{cond:conv} or \ref{condv}, with the exception of flows of surfaces \cite{A4} (see also \cite{AM} where no concavity is required but a small pinching ratio is assumed).  There are several reasons for this:  First, such conditions have been used in arguments to control curvature,  in particular in curvature pinching estimates \cites{A1, A3, Ch1, Ch2, M1, M2}, or in preserving convexity or bounding curvature \cite{Han, U}.  Secondly, some such convexity or concavity condition is necessary in order to deduce H\"older continuity of the second fundamental form using the estimates of Krylov \cite{Krylov} and Evans \cite{Evans}.  For this either condition \ref{cond:conv} or \ref{condv} will suffice.  In the case of flows of surfaces, 
no such concavity assumption is necessary \cite{Andrews2D}.  In any case
the problem of proving H\"older continuity of second derivatives for non-concave parabolic operators is quite distinct from the rest of our considerations here, so it suffices for our purposes to deal with this issue by assuming conditions \ref{cond:admit.reg} or \ref{cond:admit.reg.extend}, which we now discuss in more detail.
\item\label{rmk:sec.hold.est} Conditions \ref{cond:admit.reg} and \ref{cond:admit.reg.extend} are intended as follows:
Let $K$ be any compact subset of $\Gamma$ (or $\bar\Gamma_+\cap\{f>0\}$ in case \ref{cond:admit.reg.extend}), and $t_0>0$.  Then there exists $\beta\in(0,1]$ and $C>0$ depending only on $K$ and $t_0$ such that for any solution $X:\ M^n\times [0,t_0]\to\RR^{n+1}$ of equation \eqref{E:theflow} with principal curvatures in $K$ at every point, $\|A\|_{C^{0,\beta}(M\times[t_0/2,t_0])}\leq C$.  The following assumptions each suffice for \ref{cond:admit.reg}:
\begin{itemize}[label=\textbullet]
\item $f$ is concave or convex;
\item $f_*$ is concave or convex;
\item $n=2$ (for details see \cite{Andrews2D});
\item $f$ is nearly isotropic on $\Gamma_+$, in the sense that $\frac{\partial f}{\partial \kappa_i}\leq C(n)\frac{\partial f}{\partial\kappa_j}$ for $1\leq i,j\leq n$ at each point of $\Gamma_+$, where $C(n)$ is an explicit constant (see \cite{AM}*{Theorem 7.3});
\item $f_*$ is nearly isotropic, in the same sense.
\end{itemize}
If condition \ref{cond:extendf} holds and $f$ is concave or convex, $n=2$, or $f$ is nearly isotropic, then condition \ref{cond:admit.reg.extend} also holds.
\item\label{rmk:smooth.extend}  Condition \ref{cond:extendf} is equivalent to the assumption that $f$ can be extended to satisfy conditions \ref{T:Fconds} on an open cone $\Gamma$ containing $\Gamma_+$, such that $\kappa\in\partial\Gamma\cap\partial\Gamma_+$ implies $f(\kappa)=0$.  Condition \ref{cond:extendf*} amounts to a similar condition for $f_*$.
\item\label{rmk:extend.to.bdry} In order to make sense of conditions \ref{cond:invvanishbdry}, \ref{cond:invconvbdry} and \ref{cond:invconvinvbdry} we must make the observation that $f$ and $f_*$ each have continuous extensions to the closure of $\Gamma_+$.  We give the argument for $f$, but the same applies equally well for $f_*$:  

\begin{lemma}\label{lem:cts.ext}
A function $f$ satisfying Conditions \ref{T:Fconds}\ref{condi}--\ref{condiii} extends continuously to the closure $\bar\Gamma_+$.
\end{lemma}

\begin{proof}
Since $f$ is monotone in each argument by \ref{condii}, one can define an extension to the boundary of the positive cone by setting $f(A) = \lim_{s\to 0_+}f(A+sI)$ (where $I=(1,\dots,1)$).  This always exists since $f(A+sI)$ is nondecreasing in $s$ and is bounded below by zero (by the Euler identity using \ref{condii} and \ref{condiii}). This extension is continuous:  Continuity at the origin is
simple, since $0\leq f(b_1,\dots,b_n)\leq \max_{1\leq i\leq n}b_if(I)=\|b\|_\infty f(I)$.  It follows that $|f(B)-f(0)|\leq \|B\|_\infty f(I)$.  Next we show continuity at a non-zero boundary point $A=(0,\dots,0,a_{k+1},\dots,a_n)$:  We write
$a=\min\{a_{k+1},\dots,a_n\}>0$.  Given $\varepsilon>0$, choose $\delta>0$ so that $f(A+\delta I)<f(A)+\varepsilon$, and so that $\delta<\frac{a\varepsilon}{f(A)}$ if $f(A)>0$.  Then if $\|B-A\|_\infty<\delta$, we have $b_i<a_i+\delta$ for all $i$, so $f(B)\leq f(A+\delta I)<f(A)+\varepsilon$.  If $f(A)=0$ then we also trivially have $f(B)\geq f(A)>f(A)-\varepsilon$, while if $f(A)\neq 0$ then we have $b_i\geq \max\{a_i-\delta,0\}\geq (1-\delta/a)a_i$, so by homogeneity and monotonicity we have $f(B)\geq f((1-\delta/a)A)=(1-\delta/a)f(A)=f(A)-\frac{\delta}{a}f(A)>f(A)-\varepsilon$.  Thus in either case $|f(B)-f(A)|<\ep$ for $\|B-A\|_\infty<\delta$.
\qed\end{proof}
\item Concavity of $f$ on $\Gamma_+$ is equivalent to concavity of $F$ as a function of the components of $A$ on the space of positive definite symmetric matrices (see \cite{CNS}*{Section 3}, or apply Lemmas \ref{lem:fdotineq} and \ref{lem:d2Fvsd2f}).  Similarly, inverse-concavity of $f$ is equivalent to concavity of $f_*$, which is equivalent to concavity of the function
$F_*$ defined by $F_*(A) = F(A^{-1})^{-1}$.
\item The boundary inverse-concavity condition \ref{cond:invconvbdry} also implies inverse-concavity of the restriction of $f$ to any face of the positive cone on which it is non-zero:    For any $k>0$, we can consider the restriction of $f$ to a $k$-dimensional face, which is a function $\tilde f_k$ on the $k$-dimensional positive cone $\Gamma_+^{(k)}$, given by $\tilde f_k(x_1,\dots,x_k)=f(x_1,\dots,x_k,0,\dots,0)$.   If $\tilde f_k$ is not identically zero, then it is positive on $\Gamma_+^{(k)}$ (by monotonicity and homogeneity). To prove that $\tilde f_k$ is inverse-concave, we note that the inverse-concavity condition for $\tilde f$ can be written using Jensen's inequality in the form $(\tilde f)_*(\lambda x+(1-\lambda)y)\geq \lambda (\tilde f)_*(x)+(1-\lambda)(\tilde f)_*(y)$
for any $x,y\in\Gamma_+^{(n-1)}$ and any $\lambda\in(0,1)$.  But now for any $x\in\Gamma_+^{(k)}$, let $x_\ep = (x_1,\dots,x_k,\ep^{-1},\dots,\ep^{-1})\in\Gamma_+^{(n-1)}$.  Then we have
\begin{align*}
\lim_{\ep\to 0}(\tilde f)_*(x_\ep) &= \lim_{\ep\to 0}\frac{1}{f(x_1^{-1},\dots,x_k^{-1},\ep,\dots,\ep,0)}\\
&=\frac{1}{f(x_1^{-1},\dots,x_k^{-1},0,\dots,0)}\\
&=(\tilde f_k)_*(x),
\end{align*}
where we used Lemma \ref{lem:cts.ext} and the fact that $f>0$ in taking the limit.
Furthermore since we have $(\lambda x+(1-\lambda)y)_\ep = \lambda x_\ep+(1-\lambda)y_\ep$, we have
\begin{align*}
(\tilde f_k)_*(\lambda x+&(1-\lambda) y)- \lambda (\tilde f_k)_*(x)-(1-\lambda)(\tilde f_k)_*(y)\\
&= \lim_{\ep\to 0}\left(\tilde f)_*((\lambda x+(1-\lambda)y)_\ep)-\lambda(\tilde f)_*(x_\ep)
-(1-\lambda)(\tilde f)_*(y_\ep)\right)\\
&=\lim_{\ep\to 0}\left((\tilde f)_*(\lambda x_\ep+(1-\lambda)y_\ep)-\lambda(\tilde f)_*(x_\ep)
-(1-\lambda)(\tilde f)_*(y_\ep)\right)\\
&\geq 0,
\end{align*}
so $\tilde f_k$ is inverse-concave.
\item We will show in the course of the paper that the concavity conditions \ref{cond:conv} and \ref{condv} are stronger than needed to keep curvature bounded or to prevent loss of convexity.  In particular, in the case $\alpha=1$ we prove that Condition \ref{cond:invconvinvbdry} is essentially necessary (Theorem \ref{thm:alpha1losesmooth}) and sufficient (Theorem \ref{thm:contract.point}\ref{item:ii2}\ref{subitem4}) for curvature to remain bounded, while for $\alpha\neq 1$ the stronger condition \ref{cond:invvanishbdry} is essentially necessary (Theorem \ref{thm:curv_blowup}) and sufficient (Theorem \ref{thm:contract.point}\ref{item:i2}).
In all cases the condition \ref{cond:invconvbdry} is necessary (Theorem \ref{thm:loss.convex}) and sufficient (Theorem \ref{thm:fconc} in the case $\alpha=1$) to avoid loss of convexity of solutions (however see the caveat in remark \ref{rmk:caveat}).  Condition \ref{cond:invvanishbdry} also arises as a necessary and sufficient condition for flat sides to immediately begin to move (Theorem \ref{thm:alpha1flatsides}) in the case $\alpha=1$.
\end{enumerate}
\end{remark}

Conditions \ref{T:Fconds} and \ref{cond:extra} correspond to properties of the speed $S$, as discussed in \cite{A3}.

In this paper we will use notation similar to that in \cite{A3}, \cite{M2}, \cite{Hu1} and \cite{U}.  In particular, $g = \left\{ g_{ij}\right\}$ and $A = \left\{ h_{ij}\right\}$ denote respectively the metric and second fundamental form of $M_{t}$.  The mean curvature of $M$ is
$$H= g^{ij}h_{ij}$$
where $g^{ij}$ is the $\left( i, j\right)$-entry of the inverse of the matrix $\left( g_{ij}\right)$.  Throughout this paper we sum over repeated indices from $1$ to $n$ unless otherwise indicated.  In computations on the hypersurface $M_{t}$, raised indices indicate contraction with the metric.

We will denote by $\left( \dot{S}^{kl} \right)$ the matrix of first partial derivatives of $S$ with respect to the components of its first argument:
$$\left. \frac{\partial}{\partial s} S\left( A+sB,g\right) \right|_{s=0} = \dot S^{kl} \left( A,g\right) B_{kl} \mbox{.}$$
Similarly for the second partial derivatives of $S$ we write
$$\left. \frac{\partial^{2}}{\partial s^{2}} S\left( A+sB,g\right) \right|_{s=0} = \ddot S^{kl, rs} \left( A,g\right) B_{kl} B_{rs} \mbox{.}$$
We will normally suppress the argument $(A,g)$, so that $\dot S^{ij}$ will always be understood to be evaluated at $(A,g)$ and partial derivatives of $f$ at $\kappa\left(A,g\right)$.  We will use similar notation and conventions for other functions of matrices and eigenvalues when we differentiate them.

It will sometimes be convenient to work with the support function, which is the normal component of the position vector $s\left( x, t\right)= \left< X\left( x, t\right), \nu\left( x, t\right) \right>$.  We will also use map $\Psi$ which describes the speed as a function of the principal radii of curvature $r_i=\kappa_i^{-1}$, given by
\begin{equation} \label{E:defnPsi}
  \Psi \left( A^{-1}\right) = -S\left(A\right)=-F(A)^{\alpha} = -F_*(A^{-1})^{-\alpha}.
\end{equation}
Note that $\Psi(R) =\psi({\mathfrak r})= -f_*(r_1,\dots,r_n)^{-\alpha}$, where ${\mathfrak r}=(r_1,\dots,r_n)$ are the eigenvalues of the symmetric matrix $R$.  Note that concavity of $f$ is equivalent to inverse-concavity of $f_*$, and similarly inverse-concavity of $f$ is equivalent to concavity of $f_*$ (or of $\psi$).

We now collect some useful results concerning functions $f$ satisfying Conditions \ref{T:Fconds} and various parts of Conditions \ref{cond:extra}:

\begin{lemma}[\cite{EckerHuisken}*{Lemma 2}, \cite{A1}*{Lemma 2.20}]\label{lem:fdotineq}
If $f$ is concave on $\Gamma_+$, then for any $i\neq j$, $(\dot f^i-\dot f^j)(\kappa_i-\kappa_j)\leq 0$.
\end{lemma}

\begin{lemma}[\cite{A3}*{Theorem 2.1}] \label{T:ic}
A function $f$ satisfying Conditions \ref{T:Fconds} is inverse concave if and only if the following matrix inequality holds:
  $$Q\left[ f\right] := \left( \frac{\partial^{2} f}{\partial x_{i} \partial x_{j}} + \frac{2}{x_{i}} \frac{\partial f}{\partial x_{i}} \delta_{ij} \right) \geq 0 \mbox{.}$$
\end{lemma} 

\begin{proof}
Since $f$ (and hence also $f_*$) is smooth,
inverse-concavity is equivalent to the requirement that $\frac{d^2}{ds^2}f_*(z(s))\big|_{s=0}\leq 0$, where $z_i(s)=x_i^{-1}+sx_i^{-2}b_i$, for all $x\in\Gamma_+$ and $b\in\RR^n$.  
A direct computation gives
$$
\frac{d^2}{ds^2}f_*(z(s))\Big|_{s=0} = \frac{1}{f(x)^2}\left[-\ddot f^{ij}-\frac{2\dot f^i\delta_{ij}}{x_i}
+\frac{2}{f}\dot f^i\dot f^j\right]b_ib_j=-\frac{1}{f^2}Z(b,b).
$$
Any such $b$ can be written in the form $b=ax+w$ where $a=\frac{\dot f^iv_i}{f}$ and $\dot f^iw_i=0$.  The homogeneity of $f$ implies $\ddot f^{ij}x_i=0$, so $Z(b,b)=Q[f](w,w)$.  Thus if $Q[f]\geq 0$ then $Z\geq 0$ and $f$ is inverse-concave.  Conversely, $Q[f](v,v) = 2a^2f+Q[f](w,w)=2a^2f+Z(b,b)$, so if $f$ is inverse-concave then $Z\geq 0$ and $Q[f]\geq 0$.
\qed\end{proof}

\begin{lemma} \label{T:lowerdf}
  If $f$ satisfies Conditions \ref{T:Fconds} and is concave, then
  $\sum_{i=1}^n\dot f^i\geq 1$ on $\Gamma_+$.
\end{lemma}

\begin{proof}
By concavity we have $0\geq \frac{d^2}{ds^2}f(A+sI) = \frac{d}{ds}\left(\sum_{i=1}^n\dot f^i(A+sI)\right)$, and hence
$$
\sum_{i=1}^n\dot f^i(A)\geq \lim_{s\to\infty}\sum_{i=1}^n\dot f^i(A+sI) = \lim_{s\to\infty}\sum_{i=1}^n\dot f^i(I+A/s)=\sum_{i=1}^n\dot f^i(I),
$$
where we used the fact that $\dot f^i$ is homogeneous of degree zero.
The normalization condition \ref{condiv} and the Euler identity imply $\sum_{i=1}\dot f^i(I)= f(I)=1$.
\qed\end{proof}

\begin{lemma} \label{lem:lowerdf2}
If $f$ satisfies Conditions \ref{T:Fconds} and is inverse-concave, then
$\sum_{i=1}^n\dot f^i\kappa_i^2\geq f^2$.
\end{lemma}

\begin{proof} 
Since $f_*$ is concave, we have by Lemma \ref{T:lowerdf}
\begin{align*}
1&\leq \frac{d}{ds}\left(\dot f^*(\kappa_1^{-1}+s,\dots,\kappa_n^{-1}+s\right)\Big|_{s=0}\\
&=\frac{d}{ds}\frac{1}{f\left((\kappa_1^{-1}+s)^{-1},\dots,(\kappa_n^{-1}+s)^{-1}\right)}\Big|_{s=0}\\
&=\frac{1}{f^2}\sum_{i=1}^n\dot f^i\kappa_i^2.
\end{align*}
\qed
\end{proof}

\begin{lemma}[\cite{A1}*{Equation 2.23}, \cite{Gerhardt}*{Lemma 1.1}, \cite{A3}*{Theorem 5.1}]\label{lem:d2Fvsd2f}

If $F(A)=f(\kappa(A))$, where $f$ is a smooth symmetric function, and $A$ is diagonal with $\kappa_i(A)\neq\kappa_j(A)$ for $i\neq j$, then
$$
\ddot F(B,B) = \sum_{i,j=1}^n\ddot f^{ij}B_{ii}B_{jj}+\sum_{i\neq j}\frac{\dot f^i-\dot f^j}{\kappa_i-\kappa_j}B_{ij}^2.
$$
\end{lemma}

\begin{lemma}\label{lem:fconvvsFconv}
If $f$ is a function satisfying Conditions \ref{T:Fconds}, and $F(A)=f(\kappa(A))$, then $f$ is concave on $\Gamma_+$ if 
and only if $F$ is concave on the cone of positive definite symmetric matrices.
\end{lemma}

\begin{lemma}\label{lem:fstardotineq}
For any inverse-concave function $f(x_1,\dots,x_n)$ on $\Gamma_+$, we have
for any $i\neq j$ and at each point $\vec{x}=(x_1,\dots,x_n)$ of $\Gamma_+$, 
$\left(\dot f^i x_i^2-\dot f^j x_j^2\right)\left(x_i-x_j\right)\geq 0$.
\end{lemma}

\begin{proof}
By definition, $f_*$ is concave at $(x_1^{-1},\dots,x_n^{-1})$.  Lemma \ref{lem:fdotineq} applies to $f_*$ to give 
$\left(\dot f_*^i-\dot f_*^j\right)\left(x_i^{-1}-x_j^{-1}\right)\leq 0$ for $i\neq j$.  The lemma follows since
$\dot f_*^i = f^{-2}\dot f^ix_i^2$.
\qed\end{proof}

We conclude this section by mentioning some examples of functions which satisfy Conditions \ref{T:Fconds} and \ref{cond:extra}:

\begin{itemize}[label=\textbullet]
  \item $f=E_{k}^{1/k}$ for any $k=1, \ldots n$, where $E_{k}$ is the $k$th elementary symmetric function of the principal curvatures of $M_{t}$,
\begin{equation*}
  E_{k}=  \frac{1}{\binom{n}{k}}\sum_{1\leq i_{1} < \ldots < i_{k} \leq n} \kappa_{i_{1}} 
  \kappa_{i_{2}}\cdots \kappa_{i_{k}}  \mbox{.}
\end{equation*}  
  These include $E_{1}=H/n$ (the mean curvature), $E_{2} =R/(n(n-1))$ (the scalar curvature) and $E_{n}=K$ (the Gauss curvature).
  \item The power means $f=H_{r} = \left( \frac1n\sum_{i=1}^{n} \kappa_{i}^{r} \right)^\frac{1}{r}$ for any real $r$.  
   The case $r=0$, corresponding to $f=K^{\frac{1}{n}}$ as explained in \cite{A3},  is covered above.
  \item $f= \frac{E_{k+1}}{E_{k}}$, for any $k=1, \ldots, n-1$ ($k=0$ covered above).   
  \item  $f= \left( \frac{E_{k}}{E_{l}} \right)^{\frac{1}{k-l}}$, $n \geq k > l \geq 0$. 
  \item Convex combinations of the above examples. 
  \item Weighted geometric means of the above examples, i.e. $f=\prod_{i=1}^k f_i^{\alpha_i}$, where $\alpha_i>0$ and $\sum_{i=1}^k\alpha_i=1$, and $f_1,\dots,f_k$ come from the examples above.
   \end{itemize}
All of the explicit examples above except the power means with $r>1$  are concave, and all except for the power means with $r<-1$ are inverse-concave.  Taking convex combinations or geometric means of concave examples produce more concave examples, and convex combinations and geometric means of inverse-concave examples produce inverse-concave examples \cite{A3}*{Section 2}.

Other examples are constructed throughout the paper (see examples \ref{ex:loss.conv}, \ref{ex:cyl3}, and \ref{ex:cyl8}).

\section{Properties of solutions}\label{sec:soln}

Two fundamental properties of solutions of \eqref{E:theflow} are the following.  Firstly, because the speed of motion is homogeneous in the principal curvatures, solutions have a scaling property.  Indeed, if $X: M^{n} \times \left[ 0, T\right)\to\RR^{n+1}$  satisfies \eqref{E:theflow}, then for each $\lambda>0$, another solution $X_{\lambda} :M^{n} \times \left[ 0, \lambda^{1+\alpha} T \right)\to\RR^{n+1}$ is given by
\begin{equation} \label{E:scaling}
  X_{\lambda}\left( x, t\right) = \lambda \, X\left( x, \lambda^{-\left( 1 + \alpha \right)} \, t\right) \mbox{.}
\end{equation}
Secondly, solutions of \eqref{E:theflow} satisfy the comparison principle: If $\left\{ M_{t}^{\left( 1\right)} \right\}$ and $\left\{ M_{t}^{\left( 2\right)} \right\}$ are two families of smooth, uniformly convex hypersurfaces, each moving under \eqref{E:theflow}, and $M_{0}^{\left( 1\right)} \cap M_{0}^{\left( 2\right)} = \emptyset$, then $M_{t}^{\left( 1\right)} \cap M_{t}^{\left( 2\right)} = \emptyset$ for all $t>0$ in the common time  interval of existence.  A local version also holds: if  $\left\{ M_{t}^{\left( 1\right)} \right\}$ and $\left\{ M_{t}^{\left( 2\right)} \right\}$ are two families of smooth, uniformly convex hypersurfaces with boundary, each moving under \eqref{E:theflow}, and $M_{0}^{\left( 1\right)} \cap M_{0}^{\left( 2\right)} = \emptyset$, $M_{t}^{\left( 1\right)} \cap \partial M_{t}^{\left( 2\right)} = M_{t}^{\left( 2\right)} \cap \partial M_{t}^{\left( 1\right)} = \emptyset$, for $t\in \left[ 0, T\right)$, then $M_{t}^{\left( 1\right)} \cap M_{t}^{\left( 2\right)} = \emptyset$ for $t\in \left[ 0, T\right)$. 

The following evolution equations on $M_{t}$ are derived in \cite{Hu1}, \cite{A1} and \cite{M2}.  

\begin{lemma} \label{T:evlneqns}
Under the flow \eqref{E:theflow},
\begin{enumerate}[label={(\roman*)}]
\item\label{eq:dtg}
$\frac{\partial}{\partial t}g_{ij}=-2Sh_{ij}$.
   \item\label{eq:dtF}
   $\frac{\partial}{\partial t}S = \mathcal{L}S + \dot{S}^{kl} h_{km} h^{m}_{\ \> l} S$, 
  \item\label{eq:dts} $ \frac{\partial}{\partial t} s = \mathcal{L} s + \dot{S}^{kl} h_{km} h^{m}{}_{l} s  - \left( 1 + \alpha \right) S$. 
  \item\label{eq:dth}
  $\frac{\partial}{\partial t}h_{ij} = \mathcal{L}h_{ij}+\ddot S^{kl,mn}\nabla_ih_{kl}\nabla_jh_{pq}
 +h_{ij}\dot S^{kl}h_{kp}h^p{}_l -(1+\alpha)Sh_{ip}h^p{}_j$, or equivalently
 \begin{align*}
 \frac{\partial}{\partial t}h_{ij} &=\mathcal{L}h_{ij} + \alpha F^{\alpha-1}\ddot F^{kl,mn}\nabla_ih_{kl}\nabla_jh_{mn}+\alpha(\alpha-1)F^{\alpha-2}\dot F^{kl}\dot F^{mn}\nabla_ih_{kl}\nabla_jh_{mn}\\
 &\quad\null+h_{ij}\dot S^{kl}h_{kp}h^p{}_l -(1+\alpha)Sh_{ip}h^p{}_j.
 \end{align*}
  \end{enumerate}
  \end{lemma}

Here we denote by $\mathcal{L}$ the operator given by $\mathcal{L} =  \dot{S}^{kl} \nabla_{k} \nabla_{l}=\alpha F^{\alpha-1}\dot F^{kl}\nabla_k\nabla_l$, where $\nabla$ is the covariant derivative on $M_{t}$.  Condition \ref{T:Fconds}\ref{condii} ensures that $\mathcal{L}$ is an elliptic operator.

For convex hypersurfaces we may alternatively adopt the Gauss map parametrisation of the flow.  This parametrisation of convex hypersurfaces was extensively used in \cite{A6}.   The support function $s:\ \mathbb{S}^n\to\RR$ of a smooth, uniformly convex hypersurface $M_0\subset\RR^{n+1}$ (given by the boundary of a convex region $\Omega$) is defined by $s(z)=\sup\{\langle x,z\rangle:\ x\in\Omega\}$.  The hypersurface $M_0$ is then given by the embedding
\begin{equation}\label{eq:X}
X \left( z\right) = s\left( z\right)z + \overline{\nabla} s\left( z\right) \mbox{,}
\end{equation}
which takes $z$ to the unique point in $M_0$ with normal direction $z$ (here $\overline{\nabla}$ is the gradient with respect to the standard metric and connection on $\mathbb{S}^{n}$).  The derivative of this map is given by
$$
\partial_iX = \rr_{ik}\bar g^{kl}\partial_lz,
$$
where $\rr_{ij}$ is given as follows in terms of the support function and its derivatives:
\begin{equation} \label{E:r}
\rr_{ij} = \overline{\nabla}_{i} \overline{\nabla}_{j} s + \overline{g}_{ij} s \mbox{,}
\end{equation}
%%%% remark on derivative of $X$
where $\overline{g}_{ij}$ is the standard metric on $\mathbb{S}^{n}$.     In particular the eigenvalues of $\rr$ with respect to the metric $\bar g$ are the inverses of the principal curvatures, or the \emph{principal radii of curvature}.
The solution of the evolution equation \eqref{E:theflow} is then given, up to a time-dependent diffeomorphism, by solving the scalar parabolic equation
 \begin{equation}\label{eq:evolsuppfn}
\frac{\partial s}{\partial t} = \Psi\left(\overline{\nabla}_{i} \overline{\nabla}_{j}s+\overline{g}_{ij} s\right)= -\left(F_*\left(\overline{\nabla}_{i} \overline{\nabla}_{j}s+\overline{g}_{ij} s\right)\right)^{-\alpha}\mbox{,}
\end{equation}
for the support function $s$, and then defining $X \left( x, t \right)$ according to Equation \eqref{eq:X} and reparametrizing to make $\frac{\partial X}{\partial t}$ normal to the hypersurface.

\begin{lemma} \label{T:evlnr}
Under equation \eqref{eq:evolsuppfn}, the following evolution equations hold:  The speed $\Psi$ evolves according to
\begin{equation*}
\frac{\partial\Psi}{\partial t}=\overline{\mathcal{L}}\Psi + \Psi\dot\Psi^{ij}\bar g_{ij},
\end{equation*}
while the
inverse second fundamental form $\rr_{ij}$ evolves by
\begin{equation*}
  \frac{\partial}{\partial t} \rr_{ij} = \overline{\mathcal L} \rr_{ij} + \ddot \Psi^{kl, pq} \overline{\nabla}_{i}\rr_{pq} \overline{\nabla}_{j} \rr_{kl} - \dot\Psi^{kl}\bar g_{kl}\rr_{ij}+(1-\alpha)\Psi \bar g_{ij}
  \end{equation*}
 where $\overline{\mathcal L}=\dot \Psi^{kl}  \overline{\nabla}_{k} \overline{\nabla}_{l}
 =\alpha F_*^{-(1+\alpha)}\dot F_*^{kl}\overline{\nabla}_k\overline{\nabla}_l$, and equivalently
 \begin{align*}
  \frac{\partial}{\partial t} \rr_{ij} &=\overline{\mathcal L} \rr_{ij} + \alpha F_*^{- \left( \alpha +1\right)} \ddot F_*^{kl, pq} \overline{\nabla}_{i}\rr_{pq} \overline{\nabla}_{j} \rr_{kl}
  - \alpha \left( \alpha +1\right) F_*^{- \left( \alpha+2\right)} \overline{\nabla}_{i} F_* \overline{\nabla}_{j} F_*\\
&\quad\null  - \alpha F_*^{- \left( \alpha+1 \right)}  \dot F_*^{kl} \overline{g}_{kl} \rr_{ij} - \left( 1-\alpha \right) F_*^{-\alpha} \overline{g}_{ij} \mbox{.} 
  \end{align*}
\end{lemma}

\begin{proof} The evolution equation for $\Psi=\frac{\partial s}{\partial t}$ follows directly by differentiating equation \eqref{eq:evolsuppfn} with respect to time.  For the evolution of $\rr_{ij}$ we differentiate equation \eqref{E:r} with respect to time, yielding
$$\frac{\partial}{\partial t} \rr_{ij} = \overline{\nabla}_{i} \overline{\nabla}_{j}  \frac{\partial s}{\partial t} + \overline{g}_{ij} \frac{\partial s}{\partial t} \mbox{.}$$
A straightforward calculation involving interchanging covariant derivatives on the first term and using symmetry yields the result.  For details of Codazzi and Simons type identities applying in this setting we refer the reader to \cite{A2}*{Equation 12}. \qed\end{proof}

Note that condition \ref{T:Fconds}\ref{condii} ensures $\overline{\mathcal L}$ is an elliptic operator.

%%% remark on existence of solutions for str convex initial data
Generalized solutions can be defined for equation \eqref{E:theflow} by a variety of constructions:  One can write the evolving hypersurfaces as level sets of a function evolving by a degenerate parabolic equation, and find viscosity solutions of this equation by elliptic regularization --- this is the so-called level set method, which has been carried out for equations of the type we consider in \cite{Ishii-Souganidis} and \cite{Goto}, building on earlier work for a more restricted class of flows in \cite{CGG} (and \cite{ES} in the case of mean curvature flow).  Alternatively, one can use the method of \emph{minimal barriers} formulated by de Giorgi, as discussed in \cite{Bellettini-Novaga}.  In the setting of this paper where we consider convex hypersurfaces, a further possible approach would be to find viscosity solutions for the equation which describes the evolving hypersurfaces as spherical graphs over some enclosed origin, or for the equation \eqref{eq:evolsuppfn} for the evolution of the support function.   The relations between these
various notions are not fully established, though in many important cases they agree.  In particular, 
the relation between level-set viscosity solutions and the minimal barriers is understood \cite{Bellettini-Novaga-Comparison}, although these do not agree in general with the generalized solutions constructed using Equation \eqref{eq:X} from the viscosity solutions of \eqref{eq:evolsuppfn} --- see section \ref{sec:lossreg} for examples where this is the case.   Many of the results we prove depend principally on comparison arguments, and so apply to level-set viscosity solutions and to barrier solutions (defined below), and we note that uniformly convex classical solutions of the flow are instances of all of the other notions of solution described above.

In most of the results in this paper we will use the following notion of solution (similar in spirit to the set-theoretic subsolutions of mean curvature flow used by Ilmanen in \cite{Ilmanen}, and including in particular the concept of minimal barriers):     Let ${\mathcal F}^-_J$ be the space of smooth admissible sub-solutions of Equation \eqref{E:theflow} on a time interval $J$, by which we mean subsets $\Omega\subset\RR^{n+1}\times J$  such that each time-slice $\Omega(t)=\Omega\cap(\RR^{n+1}\times\{t\})$ is a closed set with smooth boundary depending smoothly on $t\in J$,
such that for each $z\in\partial\Omega(t)$ the principal curvatures $\kappa=(\kappa_1,\dots,\kappa_n)$ to $\partial\Omega(t)$ at $z$ lie in $\Gamma$, and the inward normal velocity of $\partial\Omega(t)$ at $z$ is no less than $f(\kappa)$.  An \emph{outer barrier} on a time interval $I$ is a subset $A\subset \RR^{n+1}\times I$, such that whenever 
$\Omega\in{\mathcal F}^-_J$ with $J=[a,b]\subset I$, and $\Omega(a)\subseteq A(a)$, then $\Omega\subseteq A$.  The space of 
such outer barriers is denoted by ${\mathcal B}_+(I)$.  Similarly, let ${\mathcal F}^+_J$ be the space of smooth super-solutions, which are subsets $\Omega\subset\RR^{n+1}\times J$  such that each time-slice $\Omega(t)=\Omega\cap(\RR^{n+1}\times\{t\})$ is a closed set with smooth boundary depending smoothly on $t\in J$,
such that for each $z\in\partial\Omega(t)$ the principal curvatures $\kappa=(\kappa_1,\dots,\kappa_n)$ to $\partial\Omega(t)$ at $z$ lie in $\Gamma$, and the inward normal velocity of $\partial\Omega(t)$ at $z$ is no greater than $f(\kappa)$.  An \emph{inner barrier} on a time interval $I$ is a subset $A\subset \RR^{n+1}\times I$, such that whenever $\Omega\in{\mathcal F}^+_J$ with $J=[a,b]\subset I$, and $A(a)\subset\Omega(a)$, then $A\subset\Omega$.  The space of inner barriers we denote by ${\mathcal B}_-(I)$.  A \emph{barrier solution} is a set which is both an outer barrier and an inner barrier.  In particular, for given $t_0\in\RR$ and $E\subset\RR^{n+1}$ we define the \emph{minimal barrier} ${\mathcal M}^-_{t_0}(E)$ and the \emph{maximal barrier} ${\mathcal M}^+_{t_0}(E)$ by
\begin{align*}
{\mathcal M}^-_{t_0}(E) &= \bigcap\left\{A:\ A\in{\mathcal B}_+([t_0,\infty)),\ A(t_0)\supseteq E\right\}\\
{\mathcal M}^+_{t_0}(E) &=\bigcup\left\{A:\ A\in{\mathcal B}_-([t_0,\infty)),\ A(t_0)\subseteq E\right\}.
\end{align*}
For convenience we also write ${\mathcal M}^{\pm}_{t_0,t}(E)={\mathcal M}^{\pm}_{t_0}(E)\cap(\RR^{n+1}\times\{t\})$ for the minimal and maximal barriers at time $t$.
Note that intersections of outer barriers are outer barriers, and unions of inner barriers are inner barriers.  In particular the minimal barrier is itself an outer barrier, and also by construction an inner barrier, and so is a barrier solution.  Similarly, the maximal barrier is a barrier solution, and barrier solutions are unique if and only if these two coincide.
The comparison principle implies that classical supersolutions (subsolutions) are outer (inner) barriers and that classical solutions are barrier solutions.  Also note that the construction produces a unique minimal barrier and a unique maximal barrier for an arbitrary initial subset $E\subset\RR^{n+1}$, and these obey the semigroup property
${\mathcal M}^{\pm}_{t_1,t}({\mathcal M}^{\pm}{t_0,t_1}(E))=
{\mathcal M}^{\pm}_{t_0,t}(E)$ for $t>t_1>t_0$.
Further properties of barrier solutions will be developed in section \ref{sec:barrier}.

\section{Examples:  Hypersurfaces which lose regularity}\label{sec:lossreg}

In this section we discuss a class of examples of singular behaviour in flows of convex hypersurfaces, in which hypersurfaces which are
initially smooth and uniformly convex develop a curvature singularity before they shrink to points.  The class of flows for which this happens is surprisingly large, and in particular many apparently very natural
flows admit such behaviour (see remark \ref{rmk:examples} below).

We begin by considering the cases $\alpha\neq 1$:

\begin{theorem}\label{thm:curv_blowup}
Suppose $\alpha\neq 1$, and 
condition \ref{cond:extendf*} holds but condition \ref{cond:invvanishbdry} fails.  Then there exists a uniformly convex non-smooth hypersurface $M_0$ with $C^\infty$ support function $s_0$ such that the solution of equation \eqref{eq:evolsuppfn} has $\min\{\rr(v,v):\ v\in TS^n,\ \|v\|=1\}<0$ for all small $t>0$.
\end{theorem}

\begin{corollary}\label{cor:losesmooth}
Under the assumptions of Theorem \ref{thm:curv_blowup} there exist smooth, uniformly convex initial hypersurfaces for which the solution of the flow \eqref{E:theflow} has $\max_{M_t}|A|^2\to\infty$ as $t\to T$, but the inradius of $M_t$ does not approach zero as $t\to T$.
\end{corollary}

The idea of the corollary is that one can find smooth, uniformly convex hypersurfaces with support function as close as desired to the support function of the non-smooth hypersurface constructed in Theorem \ref{thm:curv_blowup}.  By continuous dependence on initial data for equation \eqref{eq:evolsuppfn} the support functions of these will evolve to have a zero radius of curvature after some time, so that the principal curvatures become unbounded.

%\begin{corollary}\label{cor:weaknonsmooth}
%Under the assumptions of Theorem \ref{thm:curv_blowup}
%there exist uniformly convex non-smooth initial hypersurfaces for there is no barrier solution of \eqref{E:theflow} which is smooth for positive times.
%\end{corollary}

\begin{remark}\label{rmk:examples} Before embarking on the proof, we remark that many apparently very natural flows
satisfy the conditions of the theorem:  For example the flows defined by 
$f=\left(\frac{E_m}{E_\ell}\right)^{\frac{\alpha}{m-\ell}}$ with $0< \ell<m\leq n$ and $\alpha\neq 1$ are all
covered by the theorem.   In particular flow by positive powers $\alpha\neq 1$ of the harmonic mean curvature do not contract smooth uniformly convex hypersurfaces to points without encountering curvature singularities along the way. This is remarkable since in the case $\alpha=1$ all of these flows contract
arbitrary smooth, uniformly convex initial hypersurfaces to points with spherical asymptotic shape, as proved in \cite{A3}. 

The theorem does not apply to the flows with $f=E_m^{\frac{\alpha}{m}}$, and for these we conjecture that any weakly convex initial hypersurface does indeed evolve to immediately become smooth and uniformly convex, and contracts to a point (this applies only for $0<\alpha<1$ since we prove in Theorem \ref{thm:alphalarge} that flat sides always persist in flows with $\alpha>1$).  We prove the corresponding result for $\alpha=1$ in Theorem \ref{thm:n>2intermediate}.

Although the construction is rather similar in the two cases $0<\alpha<1$ and $\alpha>1$, the reasons for the failure in the two cases is quite different:  In the case 
$0<\alpha<1$ the loss of smoothness is driven by the `reaction'  terms in the evolution of curvature, while in case $\alpha>1$ the bad terms involve gradients of curvature, and reflect the fact that the speed $S$ is not concave (or more precisely, that $\Psi$ is not inverse-concave on the boundary of the positive cone).
\end{remark}

\begin{proof}
We will consider rotationally symmetric support functions, that is support functions which depend only on the angle from the equatorial $(n-1)$-sphere in $S^n$.  We parametrize $S^n$ using the map $
\varphi: S^{n-1}\times(-\pi/2,\pi/2)\to S^n\subset \RR^n\times\RR$ defined by
$$
\varphi(z,\theta) = (\cos\theta z,\sin\theta).
$$
We denote by $\bar g$ the standard metric on $S^{n-1}$.  Then the metric $g$ induced on $S^{n-1}\times(-\pi/2,\pi/2)$ by $\varphi$ is given by
\begin{align*}
g_{ij}&=(\cos\theta)^2\bar g_{ij};\\
g_{i\theta}&=0;\\
g_{\theta\theta}&=1.
\end{align*}
If $s$ is independent of $z$ we can compute the Hessian with respect to $g$:
\begin{align*}
\bar\nabla_i\bar\nabla_j s &= -\sin\theta\cos\theta\bar g_{ij}s_\theta;\\
\bar\nabla_i\bar\nabla_\theta s &=0;\\
\bar\nabla_\theta\bar\nabla_\theta s &=s_{\theta\theta},
\end{align*}
where $s_\theta$ denotes the derivative of $s$ with respect to $\theta$.  From this we deduce an expression for the matrix of radii of curvature from Equation \eqref{E:r}:
\begin{align*}
\rr_{ij}&=\left(s(\cos\theta)^2-s_\theta\sin\theta\cos\theta\right)\bar g_{ij};\\
\rr_{i\theta}&=0;\\
\rr_{\theta\theta}&=s_{\theta\theta}+s.
\end{align*}
The eigenvalues of this matrix with respect to the matrix $g$ are the principal radii of curvature, and these are given by
\begin{align*}
r_1 &= s_{\theta\theta}+s;\\
r_2=\dots=r_n &= s-\tan\theta s_\theta.
\end{align*} 
Thus in this symmetric situation the evolution equation \eqref{eq:evolsuppfn} becomes the following
scalar parabolic equation in one spatial variable $\theta$:
\begin{equation}\label{eq:symm_evol_s}
\frac{\partial s}{\partial t} = \psi\left(s_{\theta\theta}+s,s-s_\theta\tan\theta,\dots,s-s_\theta\tan\theta\right).
\end{equation}

Noting that $\partial_\theta r_2 = \tan\theta(r_2-r_1)$, and using the homogeneity of $\psi$ via the Euler relation, we differentiate equation \eqref{eq:symm_evol_s} twice to deduce the following evolution equation for $r=r_1=s_{\theta\theta}+s$:
\begin{align*}
\frac{\partial r}{\partial t} &=\dot\psi^1r_{\theta\theta}
+\left(\ddot\psi^{11}r_\theta+2(n-1)\tan\theta\ddot\psi^{12}(r_2-r)-(n-1)\tan\theta\dot\psi^2\right)r_\theta\\
&\quad\null +
\left((n-1)^2\tan^2\theta(r-2r_2)\ddot\psi^{22}-(n-1)(1+2\tan^2\theta)\dot\psi^2-(n-1)\tan^2\theta\ddot\psi^{12}r_2\right.\\
&\quad\null\left.\qquad+(1+\alpha)\tan^2\theta\dot\psi^1-(1+2\tan^2\theta)\dot\psi^1\right)r\\
&\quad\null +(1-\alpha)(1-\alpha\tan^2\theta)\psi.
\end{align*}

The result now follows from a simple computation for a particular choice of initial data:  In the case $\alpha\in(0,1)$ we choose $\theta_0<\pi/2$ such that $1-\alpha\tan^2\theta_0>0$, and let $r_0$ be a $C^\infty$ function which is even, $\pi$-periodic, zero on $[-\theta_0,\theta_0]$, and positive on $(\theta_0,\pi/2]$.  

In the case $\alpha>1$ we choose $0<\theta_1<\theta_2<\pi/2$ such that $1-\alpha\tan^2\theta_1<0$, and let $r_0$ be $C^\infty$, even, $\pi$-periodic, zero on $[\theta_1,\theta_2]$, and positive on $[0,\theta_1)$ and $(\theta_2,\pi/2)$.

In both cases we let $s_0$ be the unique even $\pi$-periodic solution of $(s_0)_{\theta\theta}+s_0=r_0$, given explicitly by
$$
s_0(\theta) = \sin\theta\int_0^\theta r_0(\tau)\cos\tau d\tau + \cos\theta \int_\theta^{\pi/2}r_0(\tau)\sin\tau d\tau.
$$
Then $s_0$ is $C^\infty$, and equation \eqref{eq:evolsuppfn} is strictly parabolic at $s_0$.  Therefore there exists a unique smooth solution $s(\theta,t)$ on some time interval $[0,\delta)$.  If $\theta\in [0,\theta_0]$ (in the case $0<\alpha<1$) or $\theta\in [\theta_1,\theta_2]$ (in the case $\alpha>1$) then 
$r(\theta,0)=r_\theta(\theta,0)=r_{\theta\theta}(\theta,0)=0$, but $\psi(\theta,0)>0$, and therefore
$\frac{\partial r}{\partial t}(\theta,0)<0$ and $r(\theta,t)<0$ for $t>0$ sufficiently small.

\begin{figure}\hskip 1 cm
\includegraphics[scale=0.3]{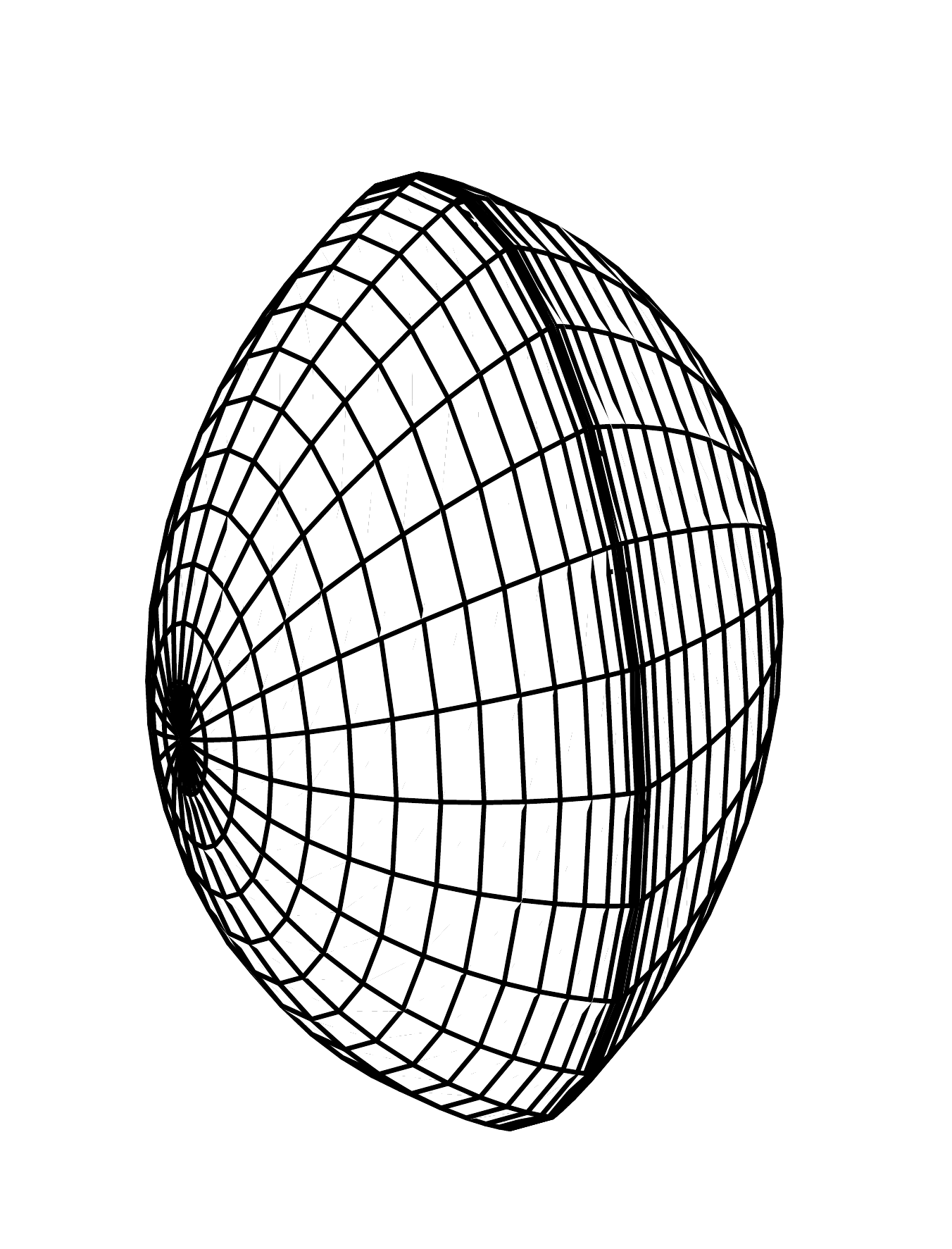}\hskip 1.5 cm
\includegraphics[scale=0.3]{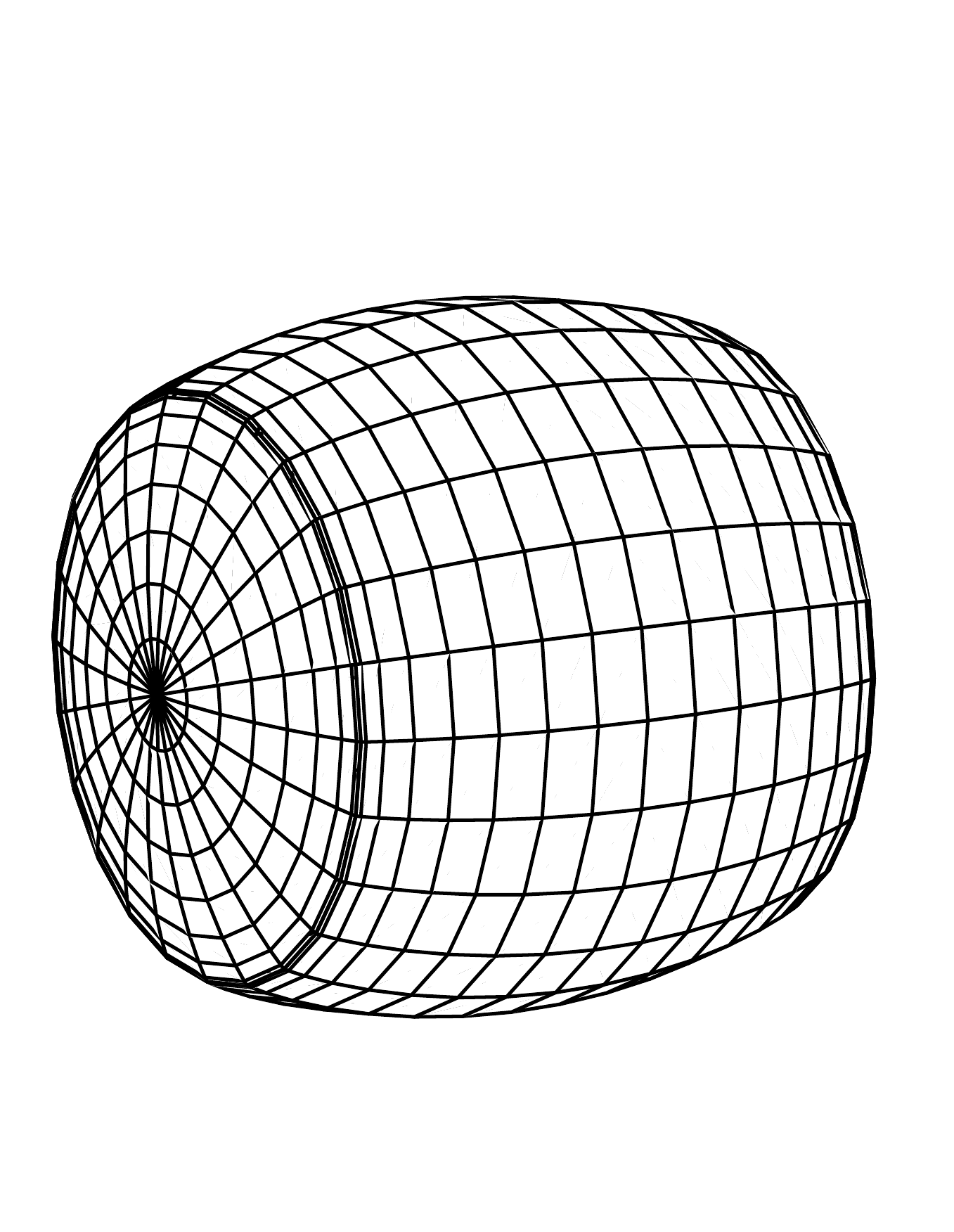}
\caption{Examples of surfaces which develop negative radius of curvature in cases $0<\alpha<1$ and $\alpha>1$.}\label{fig:loss_smooth}
\end{figure}

\begin{figure}\hskip 1 cm
\includegraphics[scale=0.3]{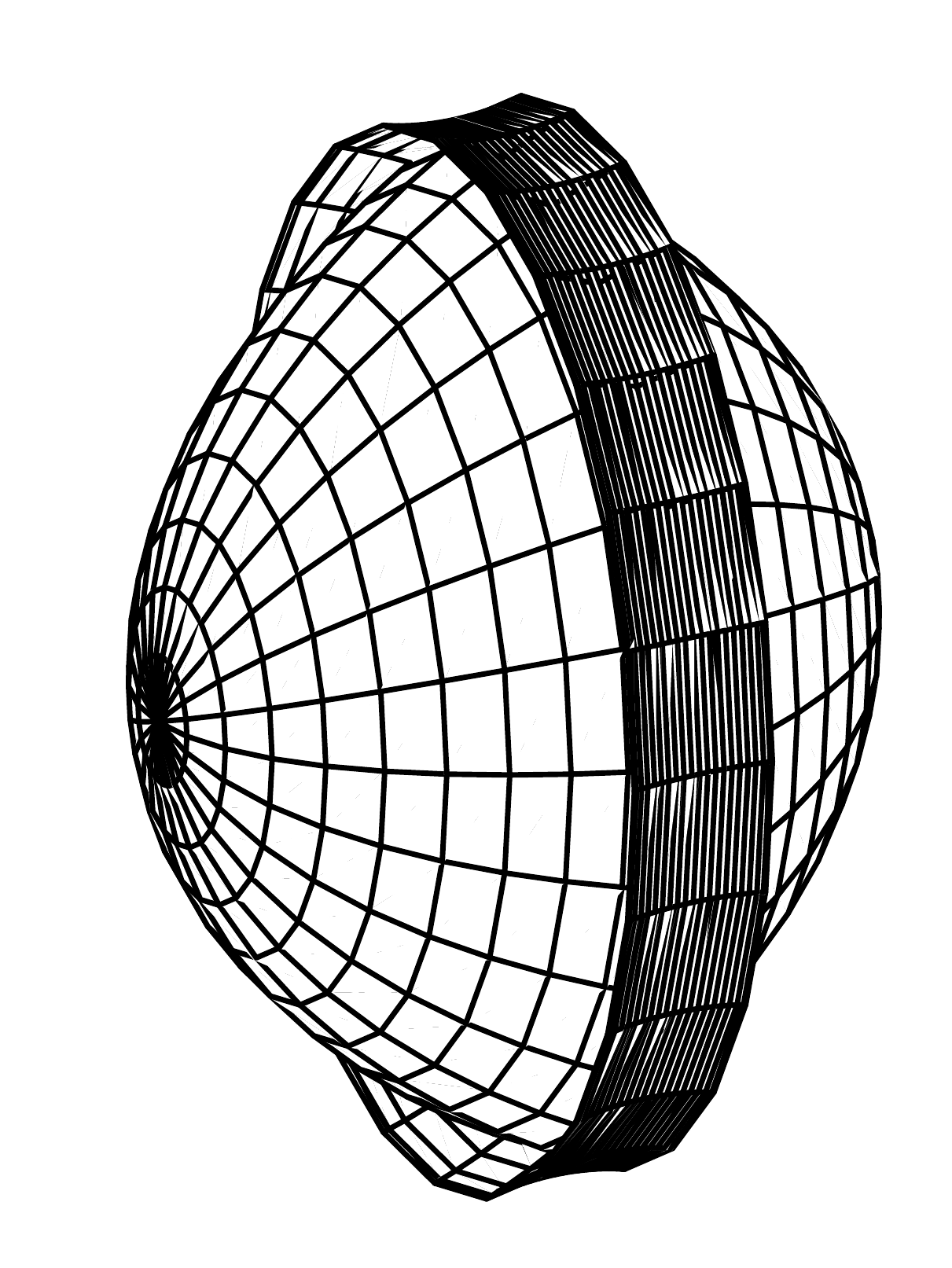}\hskip 1.5 cm
\includegraphics[scale=0.3]{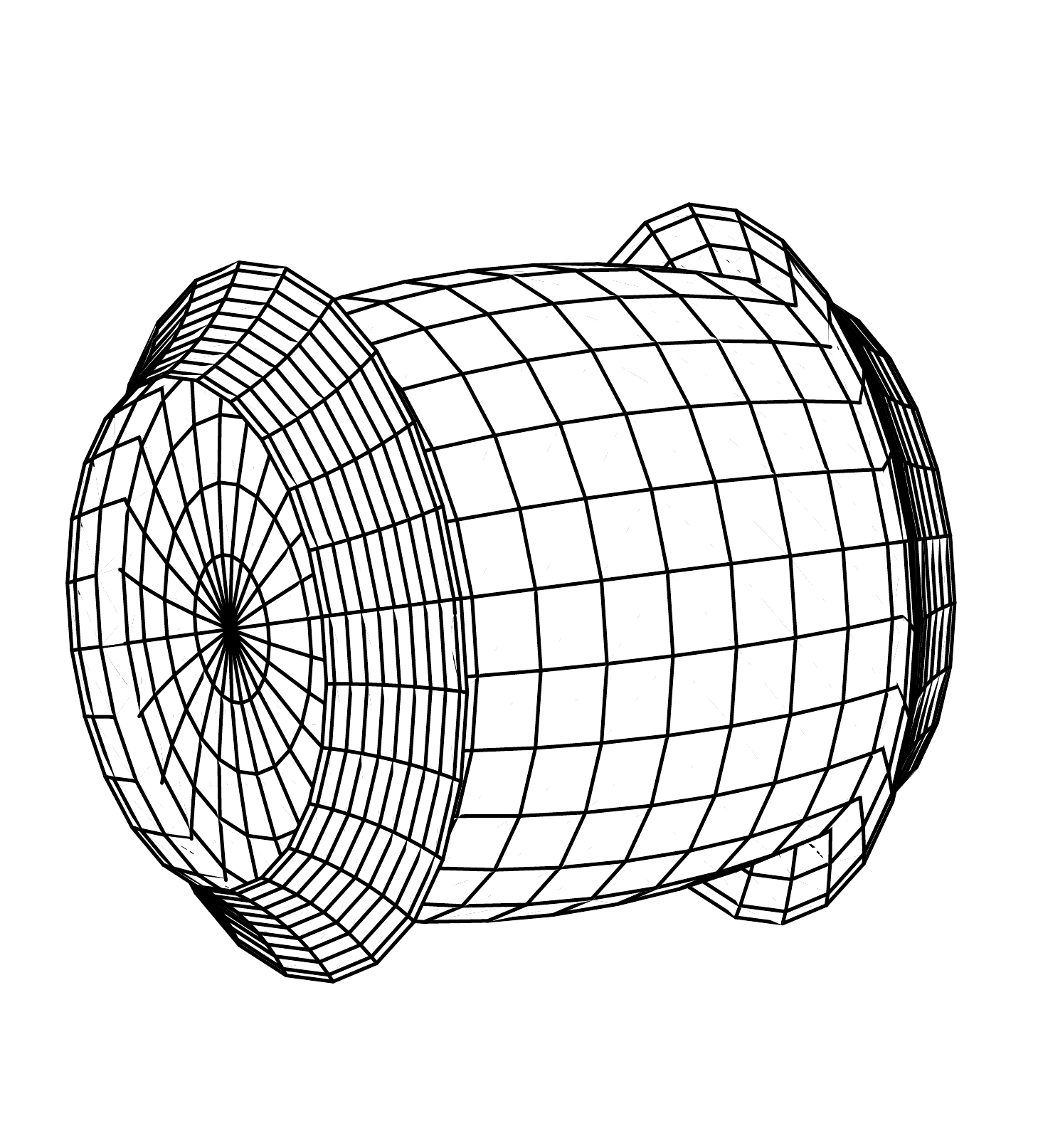}
\caption{The images of the map $X$ in equation \eqref{eq:X} after the support functions of the surfaces in Figure \ref{fig:loss_smooth} are evolved for a short time according to equation \eqref{eq:evolsuppfn}.}

\label{fig:growhorns}
\end{figure}

The corollary now follows:  Since the flow \eqref{eq:evolsuppfn} is strictly parabolic, the solution depends continuously on the initial data, so nearby smooth, uniformly convex data (such as $s_0+\varepsilon$ for $\varepsilon>0$ small) will evolve so that the radii of curvature become negative while the inradius does not shrink to zero.  On the corresponding hypersurface, there is a unique smooth
solution for a short time, for which the curvature approaches infinity while the inradius and diameter
do not approach zero.  
\qed\end{proof}

The case $\alpha=1$ admits fewer examples where
smoothness is lost:

\begin{theorem}\label{thm:alpha1losesmooth}
Suppose that $\alpha=1$, and that condition \ref{cond:extendf*} holds but the boundary inverse-concavity condition \ref{cond:invconvinvbdry} does not hold.  Then there exists a uniformly convex non-smooth hypersurface $M_0$ with smooth support function $s_0$ such that
the solution of Equation \eqref{eq:evolsuppfn} has $\min\{\rr(v,v):\ v\in TS^n,\ \|v\|=1\}<0$ for all small $t>0$.
\end{theorem}

\begin{proof}
The idea of the construction is the same as for Theorem \ref{thm:curv_blowup}:  We produce a non-smooth uniformly convex hypersurface with smooth support function, with a point and direction in which the radius of curvature is zero and has negative time derivative.  The construction is more difficult in this case since counterexamples cannot in general be constructed in the class of axially symmetric hypersurfaces.

We choose the point $e_0=(0,1)\in S^n\subset\RR^{n+1}\simeq\RR^n\times\RR$, and the direction of zero radius of curvature will be the vector $(e_1,0)$.  We construct the support
function locally first, and then show how to patch this smoothly into a suitable sphere.  For 
points $\vec{z}=(z,z_0)$ in the upper hemisphere we write the support function $s$ in the form
\begin{equation}\label{eq:sfromsigma}
s(\vec{z}) = z_0 \sigma\left(\frac{z}{z_0}\right),
\end{equation}
and observe that the non-negativity of the matrix $\rr_{ij}$ is equivalent to convexity of $\sigma$.
In fact we have the following identities:  Write $E_i = \bar\nabla_i z - \frac{\bar\nabla_iz_0}{z_0}z$, and note that since $\bar\nabla_i\bar\nabla_j\vec{z}=-\bar g_{ij}\vec{z}$ we have
$\bar\nabla_i E_j = -\frac{\bar\nabla_jz_0}{z_0}E_i$.  Then differentiating \eqref{eq:sfromsigma} gives
\begin{equation}
\bar\nabla_is\big|_{\vec{z}}=\bar\nabla_i z_0\sigma+D\sigma(E_i).
\end{equation}
Differentiating again gives an expression for $\rr_{ij}$ by Equation \eqref{E:r}:
\begin{equation}\label{eq:rij}
\rr_{ij}\big|_{\vec{z}} = z_0^{-1}D^2\sigma(E_i,E_j).
\end{equation}
Differentiating further, we obtain expression for covariant derivatives of $\rr_{ij}$:
\begin{align}
\bar\nabla_k\rr_{ij}\big|_{\vec{z}}&=z_0^{-2}D^3\sigma(E_i,E_j,E_k)-\frac{\bar\nabla_kz_0}{z_0^2}D^2\sigma(E_i,E_j)\label{eq:nbrij}\\
&\quad\null-\frac{\bar\nabla_iz_0}{z_0^2}D^2\sigma(E_j,E_k)-\frac{\bar\nabla_jz_0}{z_0^2}D^2\sigma(E_k,E_i).
\notag
\end{align}
A final differentiation gives an expression for second covariant derivatives, which we evaluate at the point $(0,1)$ where $z_0=1$, $\bar\nabla z_0=0$, and $E_i=e_i$:
\begin{equation}
\bar\nabla_l\bar\nabla_k\rr_{ij}\big|_{\vec{z}}=D_lD_kD_iD_j\sigma+\bar g_{kl}D_iD_j\sigma+\bar g_{il}D_jD_k\sigma+\bar g_{jl}D_iD_k\sigma.
\end{equation}
At this point the equations \eqref{eq:rij} and \eqref{eq:nbrij} become
\begin{align*}
\rr_{ij}&=D_iD_j\sigma;\\
\bar\nabla_k\rr_{ij} &= D_iD_jD_k\sigma.
\end{align*}
Now we make our choice of $\sigma$:  
By assumption, the restriction of $\psi$ to $\tilde\Gamma_0$ is not inverse-concave.  Therefore by 
Lemma \ref{T:ic} there exists a point $(0,a_2,\dots,a_n)$ with $a_i>0$ for $i=2,\dots,n$, and a unit vector $b=(0,b_2,\dots,b_n)$ such that 
\begin{equation}\label{eq:chooseab}
\sum_{i,j=2}^n\left(\ddot f_*^{ij}+\frac{2}{r_i}\dot f_*^j\delta_{ij}\right)b_ib_j=-\delta<0.
\end{equation}
From the proof of Lemma \ref{T:ic} (or by adding a suitable multiple of $a_i$ to $b_i$) it is possible to choose $b_i$ such that $\sum_{i=2}^n\dot f_*^ib_i=0$.
In a suitably small ball $B_r(0)$ to be chosen below, we define
\begin{equation}\label{eq:defsigma}
\sigma(z) = \frac{c_1}{24}z_1^4+\frac12\sum_{k=2}^n \left(a_k+b_kz_1+\frac12c_kz_1^2\right)z_k^2,
\end{equation}
where $a_k$ and $b_k$ are as above, and we choose $c_1=\frac{\delta}{2n\dot f_*^1}$ and $c_k = \frac{2b_k^2}{a_k}+\frac{\delta}{2n\dot f_*^k}$ for $k\geq 2$.  With this choice we have at the point $(0,1)$ that $\rr_{kl}=\text{\rm diag}(0,a_2,\dots,a_n)$, $\bar\nabla_1\rr_{kl}=\text{\rm diag}(0,b_2,\dots,b_n)$, and $\bar\nabla_k\bar\nabla_l\rr_{11} = \text{\rm diag}(c_1,c_2,\dots,c_n)$.  Since $f_*$ is smooth on $\tilde\Gamma_0$ we can assume that $a_k\neq a_l$ for each $k\neq l$.

With these choices we check that $\frac{\partial}{\partial t}\rr_{11}<0$ under the flow \eqref{eq:evolsuppfn} at the point $(0,1)$:  By Lemma \ref{T:evlnr}, since $\alpha=1$ and $\rr_{11}=0$,
\begin{align*}
\frac{\partial}{\partial t}\rr_{11} &= F_*^{-2}\dot F_*^{kl}\bar\nabla_k\bar\nabla_l\rr_{11} +F_*^{-2}\ddot F_*^{kl,pq}\bar\nabla_1\rr_{kl}\bar\nabla_1\rr_{pq}-2F_*^{-3}(\nabla_1F_*)^2\\
&=f_*^{-2}\sum_{k=1}^n\dot f_*^kc_k + f_*^{-2}\sum_{k,l=2}^n\ddot f_*^{kl}b_kb_l\\
&\quad\null + 2f_*^{-2}\sum_{1\leq k<l\leq n}\frac{\dot f_*^k-\dot f_*^l}{a_k-a_l}\bar\nabla_1\rr_{kl}^2
-2f_*^{-3}\left(\sum_{i=2}^n\dot f_*^ib_i\right)^2,
\end{align*}
where we used the expression from Lemma \ref{lem:d2Fvsd2f} in the last line.
The last two terms vanish by our choice of $b_i$ and since $\bar\nabla_1\rr_{kl}=0$ for $k\neq l$, and the choice of $c_k$ gives
$$
\frac{\partial}{\partial t}\rr_{11} = f_*^{-2}\sum_{k=1}^n\dot f_*^k\frac{\delta}{2n\dot f_*^k}
+f_*^{-2}\sum_{k,l=2}^n\left(\ddot f_*^{kl}+2\frac{\dot f_*^k}{a_k}\delta_{kl}\right)b_kb_l = -\frac{\delta}{2f_*^2}<0,
$$
where we used equation \eqref{eq:chooseab} in the second and third terms.

Now we show that we can choose $r>0$ sufficiently small so that $\sigma$ is convex on $B_r(0)$ (and in fact has positive definite Hessian away from the origin).  To prove this, we first compute the Hessian of $\sigma$:
\begin{align}\label{eq:hess.sigma}
D^2\sigma(v,v) &= \sum_{k=1}^n\frac{c_k}{2}z_k^2v_1^2+2\sum_{k=2}^n(b_kz_k+c_kz_1z_k)v_1v_k
+\sum_{k=2}^n\left(a_k+b_kz_1+\frac{c_k}{2}z_1^2\right)v_k^2\notag\\
&=\frac{c_1}{2}z_1^2v_1^2+\sum_{k=2}^n\frac{c_k}{2}z_1^2v_k^2
+\sum_{k=2}^n\frac{1}{2c_k}\left(c_kz_kv_1+2b_kv_k\right)^2\notag\\
&\quad\null+\sum_{k=2}^n\left(a_k-\frac{2b_k^2}{c_k}\right)v_k^2
+\sum_{k=2}^nb_kz_1v_k^2+2\sum_{k=2}^nc_kz_1z_kv_1v_k\notag\\
&\geq \frac{c_1}{2}z_1^2v_1^2+\sum_{k=2}^n\frac{c_k}{2}z_1^2v_k^2
+\sum_{k=2}^n\frac{\delta a_k}{2nc_k\dot\psi^k}v_k^2\notag\\
&\quad\null+\sum_{k=2}^nb_kz_1v_k^2+2\sum_{k=2}^nc_kz_1z_kv_1v_k,
\end{align}
where we substituted our particular choice of $c_k$ to obtain the last term on the first line.
Now the first term on the last line is controlled by the third term on the first line for $|z|$ sufficiently small, and the last term is controlled by the first and third terms on the first line
for $|z|$ small by the Cauchy-Schwartz inequality.  This shows that for $r>0$ small enough, the function $\sigma$ is convex on $B_r(0)$, and has strictly positive definite Hessian on $B_r(0)\setminus\{0\}$.

This completes the construction locally.  The remaining steps of the proof are to show that
this local convex function can be `pasted in' to the support function of a sphere in a $C^\infty$ way to produce a non-smooth convex hypersurface with smooth support function.  We do this in two stages:  First we produce a non-smooth convex function which equals the given function
$\sigma$ in a small ball and equals the support function of a sphere outside a large ball.  Then we perform a smoothing to produce a $C^\infty$ convex function, taking care not to change the function near the origin and at large distances.

The first stage is straightforward:   The support function of a sphere with centre on the line $\{0\}\times\RR\subset\RR^{n+1}$ is equal to $\sigma'(z)=C_1\sqrt{1+|z|^2}+C_2$ (when restricted to the line $\RR^n\times\{1\}$), for $C_1>0$ and $C_2$ arbitrary.  We will choose $C_1$ and $C_2$ in such a way that the function $\tilde\sigma=\max\{\sigma,\sigma'\}$ is equal to $\sigma$ on $B_{r/3}(0)$ and equal to $\sigma'$ on $B_r(0)\setminus B_{2r/3}(0)$ (and so extends to $\RR^n$ by taking $\tilde\sigma=\sigma'$ outside $B_r(0)$, and is a globally defined convex function).  To do this we first choose $C_2$ to make $\sigma'=0$ on $\partial B_{r/3}(0)$:  That is, we choose $C_2=-C_1\sqrt{1+(r/3)^2}$.  Since $\sigma\geq 0$ on $B_{r/3}(0)$ we then have $\sigma>\sigma'$ and hence $\tilde\sigma=\sigma$ on $B_{r/3}(0)$, for any $C_1>0$.  Then we
choose $C_1$ sufficiently large to ensure that $C_1\left(\sqrt{1+(2r/3)^2}-\sqrt{1+(r/3)^2}\right)>\sup_{B_r(0)}\sigma$.   This guarantees that $\sigma'\geq C_1\left(\sqrt{1+(2r/3)^2}-\sqrt{1+(r/3)^2}\right)>\sigma$ on $B_r(0)\setminus B_{2r/3}(0)$.

Now we proceed to the second stage:  We want to smooth $\tilde\sigma$ in the region $B_{3r/4}(0)\setminus B_{r/4}(0)$ while leaving it unchanged elsewhere and keeping it convex.  
We choose a function $\eta:\ [0,\infty)\to\RR$ which is $C^\infty$, zero on $[0,r/4]\cup[3r/4,\infty)$, and strictly positive on $[r/4,3r/4]$ with $\eta\leq 1$, $|\eta'|\leq 1$ and $|\eta''|\leq 1$.  We also choose a $C^\infty$ function $\rho:\ \RR^n\to\RR$ which is zero outside $B_1(0)$, positive inside $B_1(0)$, and has $\int_{\RR^n}\rho=1$.
Then for small positive $\varepsilon$ we define
$$
\hat\sigma_\ep(z) = \int_{B_1(0)}\tilde\sigma(z+\ep\eta(|z|)w)\rho(w) dw^n=\ep^{-n}\eta(|z|)^{-n}\int_{\RR^n}\tilde\sigma(y)\rho\left(\frac{z-y}{\ep\eta(|z|)}\right)dy^n.
$$
Then $\hat\sigma_\ep$ is equal to $\tilde\sigma$ in $B_{r/4}(0)$ and on $\RR^n\setminus B_{3r/4}(0)$, 
and is $C^\infty$ (by the first expression on $B_{r/3-\ep}(0)$, and on $\RR^n\setminus B_{2r/3+\ep}(0)$, and by the second expression on $B_{3r/4-\ep}(0)\setminus B_{r/4+\ep}(0)$, thus covering every point provided $\ep<r/24$).  

We claim that $\hat\sigma_\ep$ is convex for sufficiently small $\ep>0$.  To show this we first note that there exists $\beta>0$ such that both $\sigma-\frac{\beta}{2}|z|^2$ and $\sigma'-\frac{\beta}{2}|z|^2$ are locally convex on $B_r(0)\setminus B_{r/6}(0)$, and hence $D_iD_j\tilde\sigma\geq \beta \delta_{ij}$ (weakly) also.  Also, there exists $K$ such that $|D\sigma|\leq K$ and $|D\sigma'|\leq K$, and hence $|D\tilde\sigma|\leq K$, on the same set.  Now we compute the Hessian of $\hat\sigma_\ep$:
\begin{align*}
D_iD_j\hat\sigma_\ep (z)&= \int_{B_1(0)}D_lD_k\tilde\sigma\big|_{z+\ep\eta w}\left(\delta_{ki}+\ep \eta'\frac{z_iw_k}{|z|}\right)\left(\delta_{lj}+\ep \eta'\frac{z_jw_l}{|z|}\right)\rho(w)dw^n\\
&\quad\null+\ep \int_{B_1(0)}\left(\eta''\frac{z_iz_j}{|z|^2}+\frac{\eta'}{|z|}\left(\delta_{ij}-\frac{z_iz_j}{|z|^2}\right)\right) D\tilde\sigma(w)\rho(w)\,dw^n\\
&\geq \int_{B_1(0)}\left(\beta(1-\ep)^2-\ep K(1+6/r)\right)\delta_{ij}\,\rho(w)\,dw^n\\
&\geq \frac{\beta}{2}\delta_{ij}
\end{align*}
for $\ep>0$ sufficiently small, provided $r/4\leq |z|\leq 3r/4$.  

This completes our construction:  We now have a $C^\infty$ convex function $\hat\sigma_\ep$ on $\RR^n$ which agrees with $\sigma$ on $B_{r/4}(0)$, has strictly positive Hessian away from the origin, and agrees with $C_1\sqrt{1+|z|^2}+C_2$ outside $B_r(0)$.  The function $s$ on $S^n$ defined by Equation \eqref{eq:sfromsigma} on the upper hemisphere, and equal to $C_1+C_2z_0$ on the lower hemisphere, is $C^\infty$  and is the support function of a non-smooth weakly convex hypersurface for which the solution of equation \eqref{eq:evolsuppfn} has $\rr_{11}<0$ at the point $(0,1)$ for any small $t>0$.  This completes the proof of Theorem \ref{thm:alpha1losesmooth}.
%%%%%%%%% pictures???
\qed
\end{proof}

\begin{remark}\label{rmk:lossreg}
\begin{enumerate}[label={(\arabic*).}, ref={(\arabic*)}]
\item The analogue of Corollary \ref{cor:losesmooth} % and \ref{cor:weaknonsmooth} 
also holds.  
\item Note that the assumption of Theorem \ref{thm:alpha1losesmooth} is related to a lack of concavity of the function $f$:  If we define
$\tilde f(\kappa_2,\dots,\kappa_n) = \lim_{x\to\infty} f(x,\kappa_2,\dots,\kappa_n)$, then under the conditions of the theorem $\tilde f$ exists, and the requirement is that $\tilde f$ be concave.  The condition then amounts to a kind of `concavity at infinity' for $f$.  In particular if $f$ is concave then this is always satisfied.  

Examples of flows covered by Theorem \ref{thm:alpha1losesmooth} include $f=H_r$ for $r<-1$.  Note that in Theorem \ref{thm:alpha1flatsides} we show that a necessary and sufficient condition for flat sides in the initial hypersurface to persist is 
$\lim_{x\to\infty}f(x,1,\dots,1)=\infty$, or equivalently $f_*\big|_{\partial\Gamma_+}=0$.  This condition rules out the situation described in the theorem.
\item Let us be clear on what this result means:  It implies that the classical solution starting from smooth, uniformly convex initial data may cease to exist before the hypersurface shrinks to a point.  This does not
rule out the possibility that a convex barrier solution might continue to exist until the hypersurface shrinks to a point, or even become smooth and uniformly convex again before the final time.  As indicated in figure \ref{fig:growhorns} the map $X$ defined by equation \eqref{eq:X} from the solution of equation \eqref{eq:evolsuppfn} ceases to be an immersion in these examples, and does not describe a barrier solution of equation \eqref{E:theflow}.
\item\label{item:microscopic.convexity}
The argument presented above modifies easily to give a natural condition for solutions of fully nonlinear parabolic equations of the form
$$
\frac{\partial}{\partial t}u = F(D^2u,Du,u,x,t)
$$
to preserve local convexity :  For any  unit vector $v$, let $M_v=\{A\in\text{Sym}(n):\ Av=0,\ w^TAw>0\text{\ for\ }w^Tv=0\}$.  For $A\in M_v$ there is a unique $i_v(A)\in M_v$ such that $i_vA\circ Aw=w$ for all $w$ with $w^Tv=0$.  Then we require that for every $p\in\RR^n$, unit vector $v$, and time $t$, the function $(A,z,x)\in M_v\times\RR\times\RR^n\mapsto F(i_v(A),p,z,x,t)$ is concave.  This condition is weaker than inverse-concavity of $F$.  

In the case where $F=F(D^2u)$ is $O(n)$-invariant, then we can write $F(A)=f(\lambda)$ where $\lambda=(\lambda_1,\dots,\lambda_n)$ are the eigenvalues of $A$ and $f$ is a symmetric smooth function.  Then the condition is equivalent to the requirement that $f(0,\lambda_2^{-1},\dots,\lambda_n^{-1})$ is concave in $\lambda_2,\dots,\lambda_n$, which corresponds to the condition in Theorem \ref{thm:alpha1losesmooth}.

If $F$ is strictly parabolic and smooth on a cone containing the positive cone, then the argument above shows that this condition is necessary to avoid loss of convexity.  Provided convexity is not lost at boundary points, a maximum principle argument shows that this condition is also sufficient to keep smooth locally convex solutions locally convex.  The argument in the case where $F$ is inverse-concave is given in \cite{A3}*{Theorem 3.3} (see also the papers \cite{CGM} and \cite{GuanBian}). 
\end{enumerate}
\end{remark}

\section{Examples:  Hypersurfaces which lose convexity}\label{sec:lossconv}

In this section we give a similar construction to produce flows in which smooth, uniformly convex initial hypersurfaces can evolve to become non-convex.  While these examples are more difficult to construct, they still point to the necessity for certain concavity conditions on the boundary of the positive cone.  

\begin{theorem}\label{thm:loss.convex}
Suppose condition \ref{cond:extendf} holds, but 
the boundary inverse-concavity condition \ref{cond:invconvbdry} fails.  Then there exists a smooth, weakly convex compact hypersurface $M_0=X_0(M)$ for which the solution $X_t$ of Equation \eqref{E:theflow} has negative minimum principal curvature for all small $t>0$.
\end{theorem}

Again by the continuous dependence of solutions on initial data we have the corollary that there are smooth, uniformly convex hypersurfaces for which the classical solution (and hence the generalized solution) loses convexity, under the assumptions of the Theorem.  In this sense the boundary inverse-concavity condition is necessary to guarantee that smooth, uniformly convex hypersurfaces remain convex.

\begin{proof}

The construction is parallel to that in Theorem \ref{thm:curv_blowup}, but we work with the hypersurfaces as graphs rather than using the support function.  The result is an instance of the observation \ref{item:microscopic.convexity} in Remark \ref{rmk:lossreg}.

We construct a region $\Omega_0=\left\{\!\vec{x}=(x,x_0):\ u(x)\leq x_0\leq\!\sqrt{C_1^2-|x|^2}+C_2,\ |x|\leq C_1\!\right\}$, where $u$ is a smooth convex function with $u(x) = C_2-\sqrt{C_1^2-|x|^2}$ outside a ball $B_r(0)$.  The boundary $M_0$ of $\Omega_0$ is a then a smooth convex hypersurface, and we will construct this in such a way that the smallest principal curvature is strictly negative for  every small $t>0$, under the flow \eqref{E:theflow}.

We begin by constructing the function $u$ locally:  We will choose $u(0)=0$ and $Du(0)=0$, in which case a direct computation gives the following identities at $z=0$:
\begin{align}
h_{ij}&=D_iD_ju;\label{eq:hvsD2u}\\
\nabla_kh_{ij}&=D_kD_iD_ju;\label{eq:DhvsD3u}\\
\nabla_l\nabla_kh_{ij}&=D_lD_kD_iD_ju-D_iD_ju\,D_kD_pu\,D_lD_pu\label{eq:D2hvsD4u}\\
&\quad\null-D_kD_iu\,D_jD_pu\,D_lD_pu-
D_kD_ju\,D_iD_pu\,D_lD_pu.\notag
\end{align}
%%% details here on existence of $b$, $a$
By assumption there exists $a_2,\dots,a_n>0$ and $b_2,\dots,b_n$ such that
$$
\sum_{i,j=2}^n\left(\ddot f^{ij}+\frac{2\dot f^j}{a_j}\delta_{ij}\right)b_ib_j=-\delta<0,
$$
where $\dot f$ and $\ddot f$ are evaluated at $(0,a_2,\dots,a_n)$.  Furthermore, by adding a suitable multiple of $a_i$ to $b_i$ if necessary, we can assume that  $\sum_{i=2}^n\dot f^ib_i=0$.   In analogy with the choice
\eqref{eq:defsigma} in the proof of Theorem \ref{thm:alpha1losesmooth}, we choose
$$
u(x) = \frac{c_1}{24}x_1^4+\frac12\sum_{k=2}^n \left(a_k+b_kx_1+\frac12c_kx_1^2\right)x_k^2,
$$
where $c_1=\frac{\delta}{2n\dot f^1}$ and $c_k = \frac{2b_k^2}{a_k}+\frac{\delta}{2n\dot f^k}$ for $k\geq 2$.    We observe that by \eqref{eq:hvsD2u}--\eqref{eq:D2hvsD4u} we have $h_{kl}=\text{diag}(0,a_2,\dots,a_n)$, $\nabla_1h_{kl}=\text{diag}(0,b_2,\dots,b_n)$, and $\nabla_k\nabla_lh_{11}=\text{diag}(c_1,\dots,c_n)$ at $x=0$.

The proof of Theorem \ref{thm:alpha1losesmooth} shows that $u$ is convex on $B_r(0)$ for $r>0$ sufficiently small, and uniformly convex except at the origin.   By the same argument as given there, we can modify $u$ outside $B_{r/4}(0)$, keeping it uniformly convex and smooth, to make it equal to $C_1-\sqrt{C_2^2-|x|^2}$ outside $B_r(0)$, for suitable $C_1>0$ and $C_2\in\RR$.  Finally, we verify that the minimal principal curvature becomes negative for small positive times:  At $x=0$ we have $h_{11}=0$, and from the last evolution equation in Lemma \ref{T:evlneqns} we have (since $\sum_i\dot f^ib_i=0$ and $h_{1i}=0$ for every $i$, and using Lemma \ref{lem:d2Fvsd2f})
\begin{align*}
\frac{\partial}{\partial t}h_{11} &= \alpha f^{\alpha-1}
\sum_{k=1}^n\dot f^kc_k+\alpha f^{\alpha-1}\sum_{k,l=2}^n\ddot f^{kl}b_kb_l\\
&\quad\null+2\alpha f^{\alpha-1}\sum_{1\leq k<l\leq n}\frac{\dot f^k-\dot f^l}{a_k-a_l}\left(\nabla_ih_{kl}\right)^2+\alpha(\alpha-1)f^{\alpha-2}\left(\sum_{k=2}^n\dot f^kb_k\right)^2\\
&=\alpha f^{\alpha-1}\sum_{k=1}^n\dot f^k\frac{\delta}{2n\dot f^k}
+\alpha f^{\alpha-1}\sum_{k,l=2}^n\left(\ddot f^{kl}+\frac{2\dot f^k}{a_k}\delta_{kl}\right)b_kb_l\\
&=-\frac{\alpha f^{\alpha-1}\delta}{2}\\
&<0.
\end{align*}
It follows that the smallest principal curvature immediately becomes negative under the flow \eqref{E:theflow}, and the proof is complete.
\qed\end{proof}

\begin{example}\label{ex:loss.conv}
Examples where Theorem \ref{thm:loss.convex} applies only occur for $n\geq 3$, since the boundary inverse-concavity condition always holds for $n=2$.  One such example for $n=3$ is given by 
$$
f(x,y,z) = \frac{xy}{\sqrt{x^2+y^2}}+\frac{xz}{\sqrt{x^2+z^2}}+\frac{yz}{\sqrt{y^2+z^2}},
$$
which is smooth and increasing on the cone $\Gamma=\{x+y>0,\ x+z>0,\ y+z>0\}$.  In this case we have
$$
\tilde f(y,z) = f^{-1}(0,y^{-1},z^{-1}) = yz\sqrt{y^{-2}+z^{-2}}=\sqrt{y^2+z^2},
$$
which is certainly not concave.   
\end{example}

%%%%%%%% remark on expansion flows???

\section{Convergence to points and spherical limits for uniformly convex initial data}\label{sec:topoints}

In this section we prove a result which characterises which flows will deform smooth, uniformly convex initial hypersurfaces to points.  Informed by the various examples constructed in previous sections, our result is exhaustive --- that is, it gives a necessary and sufficient condition for this result to hold.

We remark that an earlier paper \cite{Han} gives conditions under which smooth, uniformly convex
initial data give rise to solutions which remain smooth and uniformly convex until they contract to points in finite time.  His result is as follows:

\begin{theorem}[\cite{Han}]\label{thm:Han}
Suppose $f_*$ is concave, and either 
\begin{enumerate}[label={(\roman*).}, ref={(\roman*)}]
\item\label{item:i} $f_*$ extends to a continuous function on the closed positive cone $\bar{\Gamma}_+$, and is zero on the boundary $\partial\Gamma_+$; or
\item\label{item:ii} $\alpha=1$, $f$ is concave, and $\sum_{i=1}^n\frac{\partial}{\partial r_i}\left(\gamma\circ f_*\right)$ extends to a continuous function on the closed positive cone, where $\gamma$ is a positive $C^1$ function with positive derivative at positive arguments.
\end{enumerate}
Then under the flow \eqref{E:theflow}, any smooth, uniformly convex initial hypersurface $M_0=X_0(M)$ has a solution $M_t=X_t(M)$ on a finite maximal time interval $[0,T)$, where $X$ is a $C^\infty$ family of embeddings, each hypersurface $M_t$ is uniformly convex, and $M_t\to p$ as $t\to T$, for some
$p\in\RR^{n+1}$.
\end{theorem}

\begin{remark}\label{rem:Han}
\begin{enumerate}[label={(\arabic*).}, ref={(\arabic*)}]
\item\label{rmk1}
The assumption in \ref{item:i} that $f_*$ extends to a continuous function on the closed positive cone is superfluous (see Remark \ref{rmk.conds} \ref{rmk:extend.to.bdry}).
\item\label{rmk2} The assumption in \ref{item:ii} on the trace of the derivatives extending to the boundary is not needed:  This is used to obtain a curvature bound using the evolution equation for the second fundamental form, but one can instead use the evolution equation for $h_{ij}/S$.  Under the assumption that $f$ is concave, the maximum eigenvalue of $h_{ij}/S$ does not increase, and an upper bound on principal curvatures follows.
\item Apart from the points mentioned above, the assumptions in Han's theorem are rather sharp, as the counterexamples we have provided in this paper demonstrate:  At least in the case where $f_*$ extends smoothly to the boundary, Theorem \ref{thm:curv_blowup} and Corollary \ref{cor:losesmooth} show that it is necessary to assume that $f_*=0$ on $\partial\Gamma_+$ if $\alpha\neq 1$;  in the case $\alpha=1$, Theorem \ref{thm:alpha1losesmooth} shows that some kind of concavity of $f$ is required (at least `at infinity').  The assumption of concavity of $f_*$ is necessary (again only `at infinity') by the result of Theorem \ref{thm:loss.convex}.
\end{enumerate}
\end{remark}

We next provide a result which extends Theorem \ref{thm:Han}, weakening the assumptions as much as possible within the class of inverse-concave speeds:

\begin{theorem}\label{thm:contract.point}
Suppose $f_*$ is concave, and either
\begin{enumerate}[label={(\roman*).}, ref={(\roman*)}]
\item\label{item:i2} $\alpha\neq 1$ and $f_*$ approaches zero on the boundary of the positive cone, or
\item\label{item:ii2} $\alpha=1$, and either
\begin{enumerate}[label={(\alph*)}]
\item\label{subitem1} $f_*$ approaches zero on the boundary of the positive cone; or
\item\label{subitem2} $f$ is concave; or
\item\label{subitem3} $n=2$; or
\item\label{subitem4} Conditions \ref{cond:extendf*} and \ref{cond:invconvinvbdry} hold.
\end{enumerate}
\end{enumerate}
Then the conclusion of Theorem \ref{thm:Han} holds, and furthermore in case \ref{item:ii2}\ref{subitem1}--\ref{subitem3} the hypersurfaces become spherical in the sense that $\frac{X(z,t)-p}{\sqrt{T-t}}$ converges in $C^\infty$ to an embedding $X_T$ with image the unit sphere.
\end{theorem}

Note that \ref{subitem2} implies \ref{subitem4} assuming condition \ref{cond:extendf*} holds. 

\begin{proof}
Case \ref{item:i2} can be proved by essentially the same argument as given in \cite{Han}, given the comments in remark \ref{rem:Han}\ref{rmk1}.  For completeness we briefly recall the argument in this case.

There exists a smooth solution on a maximal time interval $[0,T)$.  $T$ is finite, since a sufficiently large sphere enclosing the initial hypersurface is a barrier which contracts to a point in finite time.

The evolution equation for $S$ in Lemma \ref{T:evlneqns} implies by the maximum principle that the minimum value of $S$ is non-decreasing in $t$, so there exists a constant $S_->0$ such that $S(x,t)\geq S_-$ for all $x\in M$ and $t\in [0,T)$.   

The inverse-concavity assumption allows us to deduce that the radii of curvature remain bounded by their initial maximum value, using the last evolution equation for $\rr_{ij}$ in Lemma \ref{T:evlnr}:  Since $F_*$ is concave by the concavity of $f_*$ and Lemma \ref{lem:fconvvsFconv}, the gradient terms are all non-positive, as is the first term on the second line.  If $0<\alpha\leq 1$ then the last term is also non-positive, so the maximum eigenvalue of $\rr_{ij}$ does not increase.  In the case $\alpha>1$ we observe that $F_*\leq r_{\max}$, and $\dot F_*^{ij}\bar g_{ij}=\sum_i\dot f_*^i\geq 1$ by Lemma \ref{T:lowerdf}, so at a maximum eigenvalue the last two terms in the evolution equation can be estimated as follows:
$$
- \alpha F_*^{- \left( \alpha+1 \right)}  \dot F_*^{kl} \overline{g}_{kl} \rr_{ij} - \left( 1-\alpha \right) F_*^{-\alpha} \overline{g}_{ij} \leq -F_*^{-(1+\alpha)}r_{\max}\leq 0,
$$
so the result still holds, and we have $\rr_{ij}\leq \rr_+\bar g_{ij}$ on $M\times [0,T)$. 

An upper bound on the speed can be deduced as long as the inradius stays positive:  Denoting by $R_-$ the inradius of $M_{t_0}$ for some fixed $t_0\in(0,T)$, choose the origin at the centre of an insphere of $M_{t_0}$, and note that $s(x,t)\geq R_-$ for all 
$x\in M$ and $t\in[0,t_0]$.  Then from the evolution equations for $s$ and $S$ in Lemma \ref{T:evlneqns}, we derive an evolution equation for $Q=\frac{S}{2s-R_-}$ as follows:
\begin{equation*}
\frac{\partial Q}{\partial t}
= {\mathcal{L}}Q+4\dot S^{ij}\frac{\nabla_is}{2s-R_-}\nabla_jQ
+Q^2\left(1+\alpha-R_-\frac{\dot S^{ij}h_{ip}h^p{}_j}{S}\right).
\end{equation*}
We can estimate the term in the last bracket using Lemma \ref{lem:lowerdf2}:
$$
\dot S^{ij}h_{ip}h^p{}_j = \alpha f^{\alpha-1}\dot f^i\kappa_i^2\geq \alpha f^{\alpha+1} = \alpha S^{\frac{1+\alpha}{\alpha}}.
$$
Noting that $2s-R_-\geq R_-$, this gives the estimate
$$
\frac{\partial Q}{\partial t}\leq {\mathcal L}Q+4\dot S^{ij}\frac{\nabla_is}{2s-R_-}\nabla_jQ
+Q^2\left(1+\alpha-\alpha R_-^{1+\frac1\alpha}Q^{\frac1\alpha}\right).
$$
It follows that if a new maximum value of $Q$ is attained, then the left-hand side is non-negative, while the first term on the right is non-positive and the second is zero, and hence the bracket is non-negative.  
It follows that
$$
Q(x,t)\leq \max\left\{\sup_{M_0}Q,\ \left(1+\frac1\alpha\right)^\alpha R_-^{-(1+\alpha)}\right\}
$$
for $(x,t)\in M\times[0,t_0]$.
Noting that $s\leq R_+$, this gives the estimate
$$
S(x,t)\leq 2R_+\max\left\{\sup_{M_0}Q,\ \left(1+\frac1\alpha\right)^\alpha R_-^{-(1+\alpha)}\right\}.
$$
That is, we have proved that on any time interval $[0,t_0]$ where the inradius remains positive, there exists a constant $S_+$ such that $S(x,t)\leq S_+$ for all $(x,t)\in M\times[0,t_0]$.

Next we deduce a positive lower bound on the eigenvalues of $\rr_{ij}$ (equivalently, an upper bound on principal curvatures):  The speed bounds give $S_+^{-1/\alpha}\leq f_*\leq S_-^{-1/\alpha}$.  Since $f_*$ is continuous on $\bar\Gamma_+$ (see Lemma \ref{lem:cts.ext}), the set
$K=\{S_+^{-1/\alpha}\leq f_*\leq S_-^{-1/\alpha}\}\cap \bar\Gamma_+\cap\{r_{\max}\leq\rr_+\}$ is compact, and so the continuous function $r_{\min}$ attains a minimum value on $K$.  If this minimum is zero, then the point $\vec{r}$ lies in $\partial\Gamma_+\cap K$.   But this is a contradiction, since by assumption $f_*=0$ on $\partial\Gamma_+$, so $\vec{r}\notin K$.  Therefore $\rr_-=\inf_Kr_{\min}>0$.

We have proved that the principal curvatures remain in the compact region $\{r_{\min}\geq \rr_-\}\cap\{r_{\max}\leq \rr_+\}$ of $\Gamma_+$ on $M\times[0,t_0]$.  Since $f_*$ is inverse-concave, the flow admits H\"older estimates for second derivatives (see Appendix), and by Schauder estimates the solution is bounded in $C^k$ for every $k$.  It follows that the solution can be extended beyond $t_0$, proving that the inradius converges to zero as $t$ approaches the maximal time of existence.  Since the radius of curvature is bounded, the circumradius also approaches zero, and the hypersurfaces approach a point.

Now we proceed to the case \ref{item:ii2}.  The sub-cases \ref{subitem2} and \ref{subitem3} are proved in \cite{A4} and \cite{A3} respectively, so we need only consider sub-cases \ref{subitem1} and \ref{subitem4}.  We begin with sub-case \ref{subitem1}.   The proof of convergence to a point in this case is exactly as above.
To prove that the limiting shape is spherical we prove a new pinching estimate which requires only inverse-concavity:

For this we prove a new pinching estimate which requires only inverse-concavity:

\begin{lemma}\label{lem:inv.conc.pinch}
If $\alpha=1$ and $f_*$ is concave, then under equation \eqref{eq:evolsuppfn} the scaling-invariant quantity $\sup_{v\in T_z{\mathbb S}^n,\ \|v\|=1}\left(\frac{\rr(v,v)\big|_{(z,t)}}{F_*(\rr(z,t))}\right)$ is strictly decreasing in $t$ unless $M_t$ is a totally umbillic sphere.
\end{lemma}

\begin{proof}
We consider the evolution of $T_{ij}=\rr_{ij}-CF_*\bar g_{ij}$, where $C$ is chosen to make $T_{ij}$ negative definite when $t=0$.  We will use the maximum principle introduced in \cite{A3}*{Theorem 3.2}:  If $T$ is a tensor satisfying
$$
\frac{\partial}{\partial t}T_{ij} = a^{kl}\bar\nabla_k\bar\nabla_l T_{ij}+ u^k\bar\nabla_kT_{ij}+N_{ij},
$$
and $N$ satisfies the generalized null eigenvector condition
$$
N_{ij}v^iv^j+2\inf_{\Gamma}a^{kl}\left(2\Gamma_k^p\bar\nabla_lT_{ip}v^i-\Gamma_k^p\Gamma_l^qT_{pq}\right)\leq 0,
$$
whenever $T_{ij}\leq 0$ and $T_{ij}v^j=0$,
and $T_{ij}\leq 0$ initially, then $T_{ij}\leq 0$ remains true for positive times.

The evolution equation for $T_{ij}$ follows from those for $\Psi$ and $\rr_{ij}$ in Lemma \ref{T:evlnr}:
\begin{align*}
\frac{\partial}{\partial t}T_{ij} &=\overline{\mathcal{L}}T_{ij} + F_*^{-2}\ddot F_*^{kl,pq}\bar\nabla_i\rr_{kl}\bar\nabla_j\rr_{pq}-2F_*^{-3}\bar\nabla_iF_*\bar\nabla_jF_*\\
&\quad\null+2CF_*^{-3}\dot F_*^{kl}\bar\nabla_kF_*\bar\nabla_lF_*\bar g_{ij}-F_*^{-2}\dot F_*^{kl}\bar g_{kl}T_{ij}
\end{align*}
The last term satisfies the null eigenvector condition, so can be ignored.
It suffices to prove that at any positive definite matrix $\rr_{ij}$ with $T_{ij}=\rr_{ij}-CF_*(\rr)\delta_{ij}\leq 0$ and with $T_{ij}v^iv^j=0$ for some unit vector $v$, we have for any totally symmetric $3$-tensor $B$ satisfying $(B_{kij}-C\delta_{ij}\dot F_*^{pq}B_{kpq})v^iv^j=0$ the inequality
\begin{align*}
\left(\ddot F_*^{kl,pq}-2F_*^{-1}\dot F_*^{kl}\dot F_*^{pq}\right)B_{ikl}B_{jpq}v^iv^j+2CF_*^{-1}\dot F_*^{ij}\dot F_*^{kl}\dot F_*^{pq}B_{ikl}B_{jpq}\|v\|^2&\\ 
\null+ 2\dot F_*^{kl}
\left(2\Gamma_k^p\left(B_{lip}-C\delta_{ip}\dot F_*^{rs}B_{lrs}\right)v^i-\Gamma_k^p\Gamma_l^qT_{pq}\right)&\leq 0
\end{align*}
for some choice of $\Gamma$.  By continuity it suffices to prove this when $\rr_{ij}$ has all eigenvectors distinct.  In this case we rotate so that $\rr=\text{diag}(\rr_1,\dots,\rr_n)$ with $\rr_j<\rr_i$ for $j>i$, and observe that the the null eigenvector condition implies that $v=e_1$ and $\rr_1=Cf_*$.  The condition on $B$ is then $B_{k11}=C\dot f_*^iB_{kii}$ for each $k$.  Using this and Lemma \ref{lem:d2Fvsd2f} the required inequality becomes (choosing $\Gamma_k^1=0$ for all $k$):
\begin{align*}
0&\geq\sum_{k,l}\ddot f_*^{pq}B_{1pp}B_{1qq}+2\sum_{k>1}\frac{\dot f_*^1-\dot f_*^k}{\rr_1-\rr_k}B_{k11}^2
+2\sum_{1<k<l}\frac{\dot f_*^k-\dot f_*^l}{\rr_k-\rr_l}B_{1kl}^2\\
&\quad\null - \frac{2}{C^2f_*}B_{111}^2+\frac{2}{Cf_*}\sum_{k=1}^n\dot f_*^kB_{k11}^2
-2\sum_{k=1}^n\sum_{p=2}^n\frac{\dot f_*^k}{\rr_1-\rr_p}\left(B_{1kp}\right)^2\\
&\quad\null
+2\sum_{k=1}^n\sum_{p=2}^n(\rr_1-\rr_p)\left(\Gamma_k^p+\frac{B_{1kp}}{\rr_1-\rr_p}\right)^2.
\end{align*}
It follows that the optimal choice of $\Gamma_k^p$ is the one which makes the last terms vanish.  With this choice we show that the remaining terms do indeed satisfy the required inequality:  The terms involving $\ddot f_*$ are indeed non-positive since $f_*$ is concave by assumption.  Of the remaining terms we consider next those terms involving $B_{111}$:  These are (observing $C=\rr_1/f_*$)
$$
\left(-\frac{2}{C^2f_*}+\frac{2}{Cf_*}\dot f^1\right)B_{111}^2=\frac{2}{C^2f_*^2}\left(-f_*+\dot f_*^1\rr_1\right)B_{111}^2=-\frac{2}{C^2f_*^2}\sum_{k=2}^n\dot f_*^k\rr_kB_{111}^2\leq 0,
$$
where we used the Euler identity $\sum_{k=1}^n\dot f_*^k\rr_k=f_*$.  Next we consider the terms
involving $B_{k11}$ for $k>1$, which are (since $Cf_*=\rr_1$):
$$
\left(2\frac{\dot f_*^1-\dot f_*^k}{\rr_1-\rr_k}+\frac{2\dot f^k}{\rr_1}-2\frac{\dot f_*^1}{\rr_1-\rr_k}\right)B_{k11}^2=-2\frac{\dot f_*^k\rr_k}{\rr_1(\rr_1-\rr_k)}B_{k11}^2\leq 0.
$$
The next terms we consider are those involving $B_{1kk}$ for $k>1$.  These are:
$$
-2\frac{\dot f_*^k}{\rr_1-\rr_k}B_{1kk}^2\leq 0.
$$
The only remaining terms are those which involve $B_{1kl}$ with $1<k<l$.  These are:
$$
2\left(\frac{\dot f_*^k-\dot f_*^l}{\rr_k-\rr_l}-\frac{\dot f_*^k}{\rr_1-\rr_l}-\frac{\dot f^l}{\rr_1-\rr_k}\right)B_{1kl}^2.
$$
The first term in the bracket is nonpositive by Lemma \ref{lem:fdotineq}, and the others are negative, so the inequality holds.  

Furthermore we observe that the inequality is strict unless $B_{111}=0$, $B_{k11}=B_{1kk}=0$ for $k>1$, and $B_{1kl}=0$ for $1<k<l$.  By the strong maximum principle, the inequality becomes strict (and hence the strict decrease of the theorem is proved) unless there is a parallel vector field $v=e_1$ such that $\rr_{ij}v^iv^j=Cf_*$.  Differentiating gives $C\bar\nabla_kf_*=\bar\nabla_k\rr_{11}=B_{k11}=0$, so $f_*$ is constant.  But now we can deduce that the mean radius of curvature $S_1=\bar g^{ij}\rr_{ij}$ is constant, using the Simons-type identity
$$
0=\bar\Delta f_* = \dot f_*^{ij}\nabla_i\nabla_jS_1+\bar g^{ij}\ddot f_*^{kl,pq}\bar\nabla_i\rr_{kl}\bar\nabla_j\rr_{pq}+\frac12\sum_{i,j}(\dot f_*^i-\dot f_*^j)(\rr_i-\rr_j)\leq f_*^{ij}\bar\nabla_i\bar\nabla_jS_1.
$$
where we used the concavity of $f_*$ and Lemma \ref{lem:fdotineq}.
The maximum principle applies to show $S_1$ is constant.  Now it follows that the hypersurface is a totally umbillic sphere from the result of \cite{EckerHuisken}.
\qed
\end{proof}

\begin{remark}  The argument of the last paragraph constitutes a new result for compact convex hypersurfaces for which an inverse-concave function of curvatures is constant.  This extends the main result of \cite{EckerHuisken} to allow curvature functions which are either concave or inverse-concave.
\end{remark}

The result on the asymptotic shape of the contracting hypersurfaces in sub-case \ref{subitem1} now follows as in \cite{A1}, using the following observation:

\begin{lemma}
If $f_*$ approaches zero on $\partial\Gamma_+$, then for any $C>0$ there exists $C'>0$ such that
if $\rr\in\Gamma_+$ and $\rr_{\max}\leq Cf_*(\rr)$, then $\rr_{\max}\leq C'\rr_{\min}$.
\end{lemma}

\begin{proof}
Otherwise there exists a sequence $\left(\rr(n)\right)$ in $\Gamma_+$ with $\rr_{\max}(n)\leq Cf_*(\rr(n))$ but $\rr_{\max}(n)\geq n\rr_{\min}(n)$.  But then $\tilde\rr(n)=\frac{\rr(n)}{\rr_{\max}(n)}$ has
$\tilde\rr_{\max}(n)=1$, $f_*(\tilde\rr(n))\geq 1/C$, and $\tilde\rr_{\min}(n)\leq\frac{1}{n}$.  Since the set $\{\rr\in\bar\Gamma_+:\ \rr_{\max}\leq 1\}$ is compact, there is a subsequence $\left(\tilde\rr(n')\right)$ which converges to a limit $\rr$ with $\rr_{\max}\leq 1$, $f_*(\rr)\geq 1/C$, and $\rr_{\min}=0$.  But this contradicts the assumption that $f_*$ is zero on $\partial\Gamma_+$.
\qed\end{proof}

To complete the argument in sub-case \ref{subitem4} we will require some further regularity, which we can deduce from standard results once we have the following observation:

\begin{lemma}\label{lem:unif.par}
In case \ref{subitem4}, for any $C>0$ there exist constants $0<a_-\leq a_+$ such that at any point $\rr\in\bar\Gamma_+$ with $f_*(\rr)\geq C\rr_{\max}$, $a_-\leq \dot f_*^i\leq a_+$ for $i=1,\dots,n$.  
\end{lemma}

\begin{proof}
Since $f_*$ is homogeneous of degree one, $\dot f_*^i$ is homogeneous of degree zero, and it suffices to obtain the required bounds on the set ${\mathcal K}=\{\rr\in\bar\Gamma_+:\ f_*(\rr)\geq C\rr_{\max},\ \rr_1+\dots+\rr_n=1\}$, which is compact.  On this set, $\dot f_*^i$ is continuous and non-negative, since the extension to the boundary is $C^1$ at boundary points with $f_*>0$.  Therefore $f_*^i$ attains its minimum, and the lemma follows if we can rule out the possibility of a point $\rr\in{\mathcal K}$ which is a boundary point of $\Gamma_+$, at which $\dot f_*^i=0$ for some $i$ (
by Conditions \ref{T:Fconds}\ref{condii} this cannot occur at an interior point of $\bar\Gamma_+$).   Order the components of $\rr$ so that $\rr_1\geq\rr_2\geq\dots\geq \rr_n=0$, and let $k=\max\{i:\ \rr_i>0\}>0$.  By concavity of $f_*$ and Lemma \ref{lem:fdotineq}, $\dot f_*^1\leq \dot f_*^2\leq \dots\leq \dot f_*^n$, and so by assumption $\dot f_*^1=0$.  The restriction of $f_*$ to the first $k$ components is inverse-concave on $\Gamma_+^{(k)}$, so by Lemma \ref{lem:fstardotineq} (with $f$ replaced by $f_*$) we have $\dot f_*^1\rr_1^2\geq \dot f_*^i\rr_i^2$ for $i=1,\dots,k$, and so $\dot f_*^1\geq \frac{\rr_k}{k\rr_1^2}\sum_{i=1}^k\dot f_*^i\rr_i
=\frac{f_*(\rr)\rr_k}{k\rr_1^2}>0$.  This contradicts the assumption that $\dot f_*^1=0$, and proves the lemma.
\qed\end{proof}

Given the uniform bounds on $\dot f_*$, we can now deduce some regularity estimates.  The first step is to deduce H\"older continuity of second derivatives of the support function.

On any time interval $[0,t_0]$ on which $r_-(M_t)$ has a positive lower bound, we have $\rr(z,t)\in{\mathcal A}=\left\{\rr:\ 0<f_-\leq f_*(\rr),\ \rr_{\max}\leq C\right\}$ from the estimates above.   Since $f_*$ is concave, the set ${\mathcal A}$ is convex, and we have $\frac{a_-}{C^2}\leq \dot\psi^i\leq \frac{a_+}{f_-^2}$ at each point of ${\mathcal A}$.  We extend $\psi=-1/f_*$ to a concave, uniformly monotone function $\theta$ defined on all of $\RR^n$, by setting $\theta(\rr) = \sup\left\{\psi(\tilde\rr)+D\psi(\tilde\rr)(\rr-\tilde\rr):\ \tilde\rr\in{\mathcal A}\right\}$.   This is a supremum of linear functions, hence concave, and agrees with $\psi$ on $\mathcal A$, and furthermore $\theta$ is Lipschitz and satisfies $\frac{a_-}{C^2}\leq \dot\theta^i\leq \frac{a_+}{f_-^2}$ everywhere on $\RR^n$.

Given any point $\bar z\in\Sn^n$, we choose coordinates for $\RR^{n+1}$ so that $\bar z=(0,1)\in\Sn^n\subset\RR^n\times\RR\simeq\RR^{n+1}$, consider the evolution of the function $\sigma$ on $\RR^n$ defined by equation \eqref{eq:sfromsigma}.   In this parametrization the matrix $\rr_{ij}$ is given by  \eqref{eq:rij} (with $z_0=(1+|x|^2)^{-1/2}$), and the metric $\bar g_{ij}$ and its inverse $\bar g^{ij}$ are given by 
\begin{align}\label{eq:sph.metric}
\bar g_{ij} &= \frac{1}{1+|x|^2}\left(\delta_{ij}-\frac{x^ix^j}{1+|x|^2}\right);\\
\bar g^{ij} &= (1+|x|^2)\left(\delta^{ij}+x^ix^j\right)=(1+|x|^2)Q^{ik}Q^{jk};
\end{align}
where $Q^{ij} = \delta^{ij}+\frac{x^ix^j}{1+\sqrt{1+|x|^2}}$.  Using this we derive the evolution equation
\begin{equation}\label{eq:evol.sigma}
\frac{\partial}{\partial t}\sigma(x) = -\frac{\sqrt{1+|x|^2}}{F_*\left(\rr,\bar g\right)}
=\frac{\Theta\left(Q\circ D^2\sigma\circ Q\right)}{1+|x|^2}.
\end{equation}
Here $\Theta(A)=\theta(\vec{a})$, where $\vec{a}$ is the vector of eigenvalues of $A$.   The right-hand side is concave in the components of $D^2\sigma$, and the matrix of derivatives with respect to the components of $D^2\sigma$ are given by
$$
\frac{Q^{ik}Q^{jl}\dot\Theta^{kl}}{1+|x|^2},
$$
which has eigenvalues between $\frac{a_-}{C^2(1+|x|^2)}$ and $\frac{a_+}{f_-^2}$.  In particular, on the unit ball about the origin equation \eqref{eq:evol.sigma} is concave and uniformly parabolic, and we can apply the Krylov $C^{2,\alpha}$ estimates \cite{Krylov} (or more conveniently \cite{L}*{Corollary 14.9}), yielding 
in particular that for any $\tau\in(0,1)$
\begin{equation}\label{eq:c2alphaest}
|D^2\sigma(x_2,t_2)-D^2\sigma(x_1,t_1)|\leq C(\tau)\left(|x_2-x_1|+\sqrt{|t_2-t_1|}\right)^\beta
\end{equation}
for some constant $C(\tau)>0$ and $\beta\in(0,1)$, for $(x_2,t_2)$ and $(x_1,t_1)$ in $B_{1/2}(0)\times[\tau t_0,t_0]$.  This implies a global H\"older continuity estimate on the matrix $\rr_{ij}$.

Our next step is to deduce bounds on the third derivatives of the support function, and in particular on the covariant derivatives of the matrix $\rr_{ij}$.  To do this 
we return to the evolution equation \eqref{eq:evol.sigma}, in the form
$$
\frac{\partial}{\partial t}\sigma = \frac{\Psi(Q\circ D^2\sigma\circ Q)}{1+|x|^2}.
$$
Differentiating in a spatial direction, we derive the following equation for $v=\frac{\partial\sigma}{\partial x^i}$:
$$
\frac{\partial v}{\partial t}=\frac{1}{1+|x|^2}\dot\Psi^{kl}Q^{kp}Q^{lq}D_pD_qv+W
$$
where $W$ is a smooth function of $x$ and $D^2\sigma$, and hence is bounded in $C^{0,\beta}$.  Further, the coefficients $\dot\Psi^{kl}Q^{kp}Q^{lq}/(1+|x|^2)$ are uniformly elliptic and bounded in $C^{0,\beta}$.
We can apply a Schauder estimate \cite{L}*{Theorem 4.9} to deduce bounds in $C^{2,\beta}$ for $v$ (again these hold globally since we have such bounds for short time from the local existence theorem).  In particular, we have $|\bar\nabla_k\rr_{ij}|\leq C_3$ on $\Sn^n\times[0,t_0]$.

Armed with this estimate, we can finally deduce a lower bound on the radii of curvature.  Recall that for some $C>0$ and $\gamma>0$, we have for all $(z,t)\in\Sn^n\times[0,t_0]$ that $\rr(z,t)\in{\mathcal A}_{C/2,2\gamma} = \left\{\rr\in\bar\Gamma_+:\ \rr_{\max}\leq C/2,\ f_*(\rr)\geq 2\gamma\right\}$.  The estimates above and the assumptions imply that there exist constants $a_-$, $a_+$, $M_2$ and $M_3$ such that $\rr\in{\mathcal A}_{C,\gamma}$ implies $a_-\leq \dot\psi^i(\rr)\leq a_+$, $|\ddot\psi(\rr)|\leq M_2$ and $|\dddot\psi(\rr)|\leq M_3$.  It follows in particular that there exists $\ep_1>0$ such that if $\rr\in{\mathcal A}_{C/2,2\gamma}$ and $\rr'\in\bar\Gamma_+$ with $|\rr'-\rr|\leq \ep_1=\min\left\{\frac{C}{2},\frac{\gamma}{na_+}\right\}$, then $\rr'\in{\mathcal A}_{C,\gamma}$.

We use the evolution equation from Lemma \ref{T:evlnr} to derive the following equation for $T_{ij}=\rr_{ij}-\ep\bar g_{ij}$, where $\ep=\ep_0\E^{-Lt}$, $\varepsilon>0$ is small, and $L$ is to be chosen:
$$
\frac{\partial}{\partial t}T_{ij} = \dot\Psi^{kl}\bar\nabla_k\bar\nabla_lT_{ij}+\ddot\Psi^{pq,rs}\bar\nabla_i\rr_{pq}\bar\nabla_j\rr_{rs}-\dot\Psi^{kl}\bar g_{kl}\rr_{ij}+L\ep\bar g_{ij}.
$$
For $\ep_0>0$ sufficiently small, $T_{ij}$ is positive definite at $t=0$.  We will choose $L$ to ensure that $T_{ij}$ remains positive definite on the entire interval $[0,t_0]$.  If this is not the case, then there is some first point and time $(z,t)\in\Sn^n\times(0,t_0]$ and direction $v\in T_z\Sn^n$ for which $T_{ij}v^iv^j=0$. 
We choose an orthonormal basis $\{e_i\}$ for $T_z\Sn^n$ such that $e_1=v$, and such that $\rr_{ij}=\text{\rm diag}(\rr_1,\dots,\rr_n)$ at $(z,t)$ with $\rr_{i}$ in non-decreasing order.  
Let $m=
\max\{i:\ \rr_{i}=\ep\}$.  By continuity we can assume $\rr_{i+1}>\rr_i$ for $i> m$.  We have $\bar\nabla_k\rr_{ab}=0$ for $1\leq a,b\leq m$, 
$\frac{\partial}{\partial t}\rr_{11}\leq 0$, and (as in \cite{A3}*{Theorem 3.2})
$$
0\leq \dot\Psi^{kl}\bar\nabla_k\bar\nabla_lT_{11}-\sup_{\Gamma}2\dot\Psi^{kl}
\left(2\Gamma_k^p\bar\nabla_l\rr_{p1}-\Gamma_k^p\Gamma_l^q(\rr_{pq}-\ep\bar g_{pq})\right).
$$
We make the (optimal) choice $\Gamma_k^p=\frac{\bar\nabla_k\rr_{p1}}{\rr_p-\ep}$ for $m<k,p\leq n$, and $\Gamma_k^p=0$ otherwise.  This yields
$$
\dot\Psi^{kl}\bar\nabla_k\bar\nabla_lT_{11}\geq 2\sum_{k,l>m}\frac{\dot\psi^k(\bar\nabla_1\rr_{kp})^2}{\rr_p-\ep}.
$$
Using the identity from Lemma \ref{lem:d2Fvsd2f}, we arrive at the inequality
\begin{align}
0&\geq \sum_{p,q>m}\left(\ddot\psi^{pq}+\frac{2\dot\psi^p}{\rr_p-\ep}\delta_{pq}\right)\bar\nabla_1\rr_{pp}
\bar\nabla_1\rr_{qq}\notag\\
&\quad\null+2\sum_{m<p<q\leq n}\left(\frac{\dot\psi^q-\dot\psi^p}{\rr_q-\rr_p}+\frac{\dot\psi^p}{\rr_q-\ep}+\frac{\dot\psi^q}{\rr_p-\ep}\right)(\bar\nabla_1\rr_{pq})^2+(L-\dot\psi^{pq}\bar g_{pq})\ep.\label{eq:connect.ineq}
\end{align}
Now we use the assumption that $f_*$ is inverse-concave on $\partial\Gamma_+$:
Let $\bar\rr_i = 0$ for $i=1,\dots,m$, and $\bar\rr_i=\rr_i$ for $i=m+1,\dots,n$.  Then $|\bar\rr-\rr|_\infty=\ep\leq\ep_0$, so provided $\ep_0\leq\ep_1$ we have $\bar\rr\in{\mathcal A}_{C,\gamma}$.  It follows that
$$
\sum_{p,q>m}\left(\ddot\psi^{pq}+\frac{2\dot\psi^p}{\rr_p}\delta_{pq}\right)\Big|_{\bar\rr}\bar\nabla_1\rr_{pp}
\bar\nabla_1\rr_{qq}\geq 0,
$$
and so 
\begin{align*}
\sum_{p,q>m}&\left(\ddot\psi^{pq}+\frac{2\dot\psi^p}{\rr_p-\ep}\delta_{pq}\right)\Big|_{\rr}\bar\nabla_1\rr_{pp}
\bar\nabla_1\rr_{qq}\\
&=\sum_{p,q>m}\left(\ddot\psi^{pq}+\frac{2\dot\psi^p}{\rr_p}\delta_{pq}+\frac{2\ep\dot\psi^p}{\rr_p(\rr_p-\ep)}\right)\Big|_{\rr}\bar\nabla_1\rr_{pp}
\bar\nabla_1\rr_{qq}\\
&\geq\sum_{p,q>m}\left[\left(\ddot\psi^{pq}+\frac{2\dot\psi^p}{\rr_p}\delta_{pq}\right)\Big|_{\bar\rr}
+\left(-\ep M_3-\frac{\ep M_2}{\rr_p}+\frac{2\ep a_-}{\rr_p^2}\right)\delta_{pq}\right]
\bar\nabla_1\rr_{pp}
\bar\nabla_1\rr_{qq}\\
&\geq -\ep\left(M_3+\frac{M_2^2}{8a_-}\right)\sum_{p>m}(\bar\nabla_1\rr_{pp})^2\\
&\geq -\ep\left(M_3+\frac{M_2^2}{8a_-}\right)C_3^2.
\end{align*}
Similarly we can bound the second term in \eqref{eq:connect.ineq}:  By the assumption of inverse-concavity on the boundary we have at $\bar\rr$ for each $m<p<q\leq n$
$$
\frac{\dot\psi^q-\dot\psi^p}{\rr_q-\rr_p}+\frac{\dot\psi^p}{\rr_q}+\frac{\dot\psi^q}{\rr_p}
=\frac{\dot\psi^q\rr_q^2-\dot\psi^p\rr_p^2}{\rr_p\rr_q(\rr_q-\rr_p)}\geq 0
$$
by Lemma \ref{lem:fstardotineq}.  Using this we find
\begin{align*}
&\left(
\frac{\dot\psi^q-\dot\psi^p}{\rr_q-\rr_p}+\frac{\dot\psi^p}{\rr_q-\ep}+\frac{\dot\psi^q}{\rr_p-\ep}
\right)\Big|_{\rr}\\
&=\left(
\frac{\dot\psi^q-\dot\psi^p}{\rr_q-\rr_p}+\frac{\dot\psi^p}{\rr_q}+\frac{\dot\psi^q}{\rr_p}
+\frac{\ep\dot\psi^p}{\rr_q(\rr_q-\ep)}+\frac{\ep\dot\psi^q}{\rr_p(\rr_p-\ep)}
\right)\Big|_{\rr}\\
&\geq \left(
\frac{\dot\psi^q-\dot\psi^p}{\rr_q-\rr_p}+\frac{\dot\psi^p}{\rr_q}+\frac{\dot\psi^q}{\rr_p}
\right)\Big|_{\bar\rr}\\
&\quad\null + \ep\sum_{i=1}^m\int_0^1\frac{\ddot\psi^{iq}-\ddot\psi^{ip}}{\rr_q-\rr_p}\Big|_{s\bar\rr+(1-s)\rr}\,ds-\frac{\ep M_2}{\rr_p}-\frac{\ep M_2}{\rr_q} + \frac{\ep a_-}{\rr_p^2}+\frac{\ep a_-}{\rr_q^2}\\
&\geq \ep\left(-M_3-\frac{M_2}{\rr_p}-\frac{M_2}{\rr_q}+\frac{a_-}{\rr_p^2}+\frac{a_-}{\rr_q^2}\right)\\
&\geq -\ep\left(M_3+\frac{M_2^2}{2a_-}\right).
\end{align*}
The inequality \eqref{eq:connect.ineq} becomes
\begin{equation}\label{eq:L.ineq}
0\geq \ep\left(-CM_3-C\frac{M_2^2}{a_-}\right)C_3^2+L-Ca_+,
\end{equation}
where $C$ is a constant depending only on $n$.  But the right hand side is strictly positive if we choose $L$ suitably large depending on $M_2$, $M_3$, $a_{\pm}$, $C_3$ and $n$, and we have a contradiction.  This proves that the inequality $\rr_{ij}\geq \ep_0\E^{-Lt}\bar g_{ij}$ is preserved, provided $\ep_0\leq \ep_1$ and $L$ is chosen to make the right-hand side of \eqref{eq:L.ineq} positive.

The proof of subcase \ref{subitem4} can now be completed exactly as before:  We have established that as long as the inradius remains positive the principal radii of curvature remain in a compact region of the open positive cone, and higher regularity follows by standard Schauder estimates, proving that the solution can be continued to a longer time interval.   This implies that the inradius converges to zero at the end of the interval of existence.  The circumradius also converges to zero since $\kappa_{\min}$ has a positive lower bound.
\qed\end{proof}

%%% coounterexamples to sphere limit if only weak inv-conv on bdry???

\begin{remark}\label{rmk:caveat}
\begin{enumerate}[label={\arabic*.}]
\item We do not know whether the limiting shape in case \ref{subitem4} is always spherical.  This is the case if the inverse-concavity of $f_*$ on the boundary of the positive cone is strict in non-radial directions, since then one can deduce that $f_*$ is inverse-concave on the sets $\rr_{1}=\ep(\rr_2+\rr_n)$ for small $\ep>0$, and so the inequality $\rr_1\geq \ep(\rr_2+\dots+\rr_n)$ is preserved, by a slight modification of the argument in \cite{A3}.  It then follows that the ratio of circumradius to inradius remains bounded, and one can extract a smooth hypersurface as a limit of rescalings as the final time is approached.  The result of Lemma \ref{lem:inv.conc.pinch}  then implies that the limit hypersurface is a sphere.  However it seems possible that with only weak inverse-concavity on the boundary, the ratio of inradius to circumradius could become unbounded.
\item
One should be cautious about the sense in which the conditions are necessary and sufficient:  
The conditions are necessary and sufficient in the case where $f$ and $f_*$ are $C^{2,1}$ up to the 
boundary of the positive cone at points where they are non-zero (in particular, if $f$ and $f_*$ extend
smoothly beyond the boundary at such points).
We have not ruled out the possibility that there could be flows with $\alpha\neq 1$ for which $f_*$ does not vanish on the boundary of the positive cone, but does not extend as a strictly monotone smooth function to a larger cone, and for which solutions do not lose smoothness.  In particular, is seems possible that if the coefficients of $\dot f_*$ become unbounded near the boundary of the positive cone, the resulting fast diffusion could have the effect of preventing loss of smoothness.  The same considerations apply for the requirement for concavity at infinity.
However we have not so far been able to construct any examples where this occurs.
\end{enumerate}
\end{remark}

In the case $\alpha=1$ we also have a result of this kind when the speed $f$ is concave:

\begin{theorem}\label{thm:fconc}
 Suppose $\alpha=1$, $f$ is concave, condition \ref{cond:invconvbdry} holds, and $f\in C^{2,1}_{\text{loc}}(\bar\Gamma_+\cap\{f>0\})$.  Then the solution of \eqref{E:theflow} for any smooth, uniformly convex initial hypersurface exists and remains smooth and uniformly convex on a finite maximal time interval $[0,T)$, and the inradius $r_-(M_t)$ converges to zero as $t$ approaches $T$.
\end{theorem}

In this generality we are unable to prove that the circumradius also converges to zero, so it is conceivable that the solution could collapse onto a lower-dimensional convex body of positive diameter.

In the case where $f$ is concave and zero on $\partial\Gamma_+$, it was proved in \cite{A1} that arbitrary smooth, uniformly convex initial hypersurfaces contract to points with spherical limiting shape under Equation \eqref{E:theflow}.  We do not know whether the limiting shape must be spherical under the assumptions of Theorem \ref{thm:fconc}, but as in case \ref{subitem4} of Theorem \ref{thm:contract.point}, such a result does hold if the inverse-concavity on the boundary is strict in non-radial directions.

\begin{proof}
In the case where $f$ is concave there is a well-known curvature pinching estimate, much easier to prove in this case than the result of Lemma \ref{lem:inv.conc.pinch}:

\begin{lemma}\label{lem:pinch.curv.conc}
Let $\{M_t\}_{0\leq t<T}$ be a family of uniformly convex hypersurfaces evolving by Equation \eqref{E:theflow} with $f$ concave and $\alpha=1$, and suppose at $t=0$ the inequality $h_{ij}\leq Cfg_{ij}$ holds for some $C>0$.  Then the same is true for each $t\in(0,T)$.  
\end{lemma}

\begin{proof}
From the evolution of $h_{ij}$, $f$ and $g_{ij}$ we find:
\begin{align*}
\frac{\partial}{\partial t}\left(h_{ij}-Cfg_{ij}\right)
&= {\mathcal L}\left(h_{ij}-Cfg_{ij}\right) + \ddot F^{pq,rs}\nabla_ih_{pq}\nabla_jh_{rs}\\
&\quad\null + \dot f^k\kappa_k^2\left(h_{ij}-Cfg_{ij}\right)-2fh_i{}^p\left(h_{pj}-Cfg_{pj}\right).
\end{align*}
Hamilton's maximum principle for tensors \cite{Ha}*{Theorem 9.1} applies to show that the inequality $h_{ij}-Cfg_{ij}\leq 0$ is preserved.
\qed\end{proof}

Since $f$ is concave and inverse-concave on the boundary, we have exactly as in Lemma \ref{lem:unif.par} that for any $C>0$ there exist constants $0<a_-\leq a_+$ such that at any point $\kappa\in\bar\Gamma_+$ with $f(\kappa)>0$ and $f(\rr)\geq C\kappa_{\max}$, $a_-\leq \dot f^i\leq a_+$.    The evolution equation for $S$ implies that a lower bound $f\geq f_->0$ is preserved, and we also have an upper bound $f\leq f_+$ while the inradius remains positive using the same argument as in Theorem \ref{thm:contract.point}, since that argument depends only on an estimate of the form $\sum_i\dot f^i\kappa_i^2\geq \gamma f^2$, and this holds since $\sum_i\dot f^i\kappa_i^2\geq \frac{1}{a_+}\sum_i(\dot f^i\kappa_i)^2\geq \frac{1}{na_+}\left(\sum_i\dot f^i\kappa_i\right)^2 = \frac{f^2}{na_+}$.

The next step is to deduce bounds on the derivatives of the second fundamental form.  
These can be produced in two steps as in the proof of Theorem \ref{thm:contract.point}:  The evolving hypersurfaces can be written locally as the graph of a function $u$ evolving by a uniformly parabolic scalar equation, and the concavity of $f$ allows the application of the Evans-Krylov estimates to deduce $C^{2,\beta}$ estimates.  Differentiating the scalar equation in a spatial direction gives a uniformly
parabolic equation for $v=\partial_ku$, and Schauder estimates then produce $C^{2,\beta}$ estimates for $v$, which imply bounds for $\|\nabla_ih_{jk}\|$.

Finally, we can prove that $\kappa_{\min}$ remains positive while the inradius is positive, by considering
the evolution equation for $T_{ij}=h_{ij}-\ep_0\E^{-Lt}g_{ij}$.  The details are formally almost identical to the corresponding computation in the proof of Theorem \ref{thm:contract.point}.  As before, bounds on all higher derivatives follows by Schauder theory, and it follows that the solution continues to exist until the inradius converges to zero.
%%% PROBLEM:   DOES CIRCUMRADIUS APPROACH ZERO?
\qed\end{proof}
\begin{remark}
We expect that there should be a result which holds without assuming that either $f$ or $f_*$ is concave.
One might conjecture that if conditions \ref{cond:invconvbdry}, \ref{cond:invconvinvbdry} and \ref{cond:admit.reg} hold (and $f$ and $f_*$ extend as $C^{2,1}$ functions up to the boundary of the positive cone near points where they are non-zero), and $f_*=0$ on the boundary of the positive cone if $\alpha\neq 1$, then a similar result should hold.  Under these assumptions, maximum principle arguments show that if the hypersurface remains smooth, then it remains uniformly convex, and if the support function remains smooth, then the principal curvatures remain bounded.  At present we have a result in this direction (but still requiring stronger hypotheses) only in the case $n=2$:
\end{remark}

\begin{theorem}\label{thm:nis2}
If $n=2$, and either $0<\alpha<1$ with $f\geq CH$, or $\alpha=1$, then the solution of \eqref{E:theflow} for any smooth, uniformly convex initial hypersurface exists and remains smooth and uniformly convex on a finite maximal time interval $[0,T)$, and the inradius $r_-(M_t)$ converges to zero as $t$ approaches $T$.  In the case $\alpha=1$ the surfaces $M_t$ contract to a point as $t\to T$ and the limiting shape is spherical.
\end{theorem}

\begin{proof}
The case $\alpha=1$ was proved in \cite{A4}.  In the case $0<\alpha<1$ we also use the computations of \cite{A4}:   Choosing $G=\frac{(\kappa_2-\kappa_1)^2}{(\kappa_2+\kappa_1)^2}$, we have from the computations there (equation 6 and the following discussion) that at any spatial maximum point of $G$, 
$$
\frac{\partial G}{\partial t}\leq 4(1-\alpha)\frac{KSG}{H}=(1-\alpha)SGH(1-G).
$$
An upper bound on the speed $S$ holds as long as the inradius has a positive lower bound, by the argument given in the proof of Theorem \ref{thm:contract.point}.  By assumption we also have that $H\leq f/C$ is bounded.  Therefore $1-G$ has an exponentially decaying lower bound, so remains
strictly positive.  Since the speed $S$ also has a positive lower bound (from the evolution equation for $S$), the flow is uniformly parabolic and the solution remains smooth and uniformly convex, on any time interval where the inradius has a positive lower bound (here the regularity results of \cite{Andrews2D} are used).  The result follows.
\qed
\end{proof}

\section{Bounds on displacement and flat sides}\label{S:moving}

In the previous section we were concerned with flows which keep smooth, strictly convex hypersurfaces smooth and strictly convex.  In the remaining sections of the paper we are concerned with the behaviour of
the flows starting from initial data which are boundaries of open bounded convex regions, but may be neither smooth nor uniformly convex.

In this section we are concerned with the displacement of the hypersurface in short positive time intervals.  The easy direction is to bound the displacement from above, which we accomplish using spherical barriers.  Positive lower bounds on displacement are more subtle, and require the construction of suitable barriers.  The main result of this section is a complete characterisation of the flows for which such a lower displacement bound exists.  The results for $\alpha\neq 1$ can be seen as higher-dimensional analogues of the results proved for evolving curves in \cite{A9}, but the intermediate case $\alpha=1$ introduces some new features.

\subsection{Upper displacement bound}\label{SS:upperdisp}
%% importance of hom deg 1 as a critical case
%% weakly convex, strictly $k$-convex?

As in \cite{A2}, we obtain an upper bound on the displacement of the evolving hypersurfaces using spheres enclosed within $M_{0}$ as barriers.  We denote by $R_+$ and $R_-$ the circumradius and inradius of $M_0$.  

\begin{theorem} \label{T:du}
For any smooth, uniformly convex solution $\left\{ M_{t}\right\}_{t>0}$ of \eqref{E:theflow},
\begin{equation} \label{E:hl}
  s\left( z, t\right) \geq s\left( z, 0 \right) - C(\alpha)\frac{R_{+}}{R_{-}} t^{\frac1{1+\alpha}}
\end{equation}
for all $t\in \left[0, \frac{1}{1+\alpha} R_{-}^{1+\alpha} \right]$ in the interval of existence of the solution and all $z\in S^n$.
\end{theorem}

%%%% pictures??
\begin{proof} Choose the origin $O$ to be the centre of a ball of radius $R_{-}$ enclosed by $M_{0}$.  Fix $z\in \mathbb{S}^{n}$ and define for each $\delta\in\left( 0, 1\right]$ a sphere
$$
S_{\delta} = (1-\delta)X(z,0)+\delta S_{R_-}(O)=
\left\{ y\in \mathbb{R}^{n+1} : \left| y - \left( 1- \delta \right) X\left( z, 0\right) \right| = \delta \, R_{-} \right\}  \mbox{,}
$$
where $X(.,0)$ is the embedding of $M_0$ given by \eqref{eq:X}.
$S_{\delta}$ is contained within the convex hull of $B_{R_{-}}\left( O\right)$ and $X\left( z, 0\right)$ and so, since $M_{0}$ is convex, $S_{\delta}$ is enclosed by $M_{0}$.  Under equation \eqref{E:theflow}, this sphere evolves to $S_{r\left( t\right)}\left((1-\delta)X(z,0)\right)$, where $r\left( t\right)$ evolves according to
$$
\frac{dr}{dt} = -\frac{1}{r^\alpha} \mbox{.}
$$
With initial condition $r\left( 0\right) = \delta \, R_{-}$, this ODE has solution
$$
r\left( t\right) = \left(\delta^{1+\alpha}R_{-}^{1+\alpha} - (1+\alpha)t\right)^{\frac1{1+\alpha}},
$$
and the sphere shrinks to a point at time $t= \frac{\delta^{1+\alpha}R_{-}^{1+\alpha}}{1+\alpha}$.  By the comparison principle,  each of these shrinking spheres remains enclosed by $M_t$ as long as both exist, so in particular $(1-\delta)X(z,0)$ is enclosed by $M_t$ for $0\leq t\leq \frac{\delta^{1+\alpha}R_{-}^{1+\alpha}}{1+\alpha}$.  This gives the estimate (choosing $\delta$ so that $\delta^{1+\alpha}R_-^{1+\alpha}=(1+\alpha)t$)
$$
s(z,t)=\sup\{\langle y,z\rangle:\ y\in\Omega_t\}
\geq \langle (1-\delta)X(z,0),z\rangle = (1-\delta)s(z,0).
$$
Therefore
$$
s(z,t)-s(z,0)\geq -\delta s(z,0)\geq -\frac{(1+\alpha)^{\frac1{1+\alpha}}t^{\frac1{1+\alpha}}}{R_-} (2R_+) = C(\alpha)\frac{R^+}{R^-}t^{\frac1{1+\alpha}}.
$$
\qed\end{proof}

\subsection{Lower displacement bound}\label{SS:lowerdisp}

Next we consider lower displacement bounds.  This is simplest in the case $0<\alpha<1$, where spherical barriers suffice to prove that the hypersurfaces move:

\begin{theorem}\label{thm:alphasmall}
If $0<\alpha<1$ then 
\begin{equation}\label{eq:alphasmallldb}
s(z,t)\leq s(z,0)-(1-\alpha)\alpha^{\frac{\alpha}{1-\alpha}}R_+^{-\frac{2\alpha}{1-\alpha}}t^{\frac{1}{1-\alpha}}
\end{equation}
for $0\leq t\leq \alpha^{-1}R_+^{1+\alpha}$ while the solution continues to exist.
\end{theorem}

%%% pictures
\begin{proof}
Choose the origin $O$ to be at the centre of a sphere of radius $R_+$ which encloses $M_0$.  Then for any fixed $z\in S^n$, $M_0$ is contained in the region $B_{R_+}(O)\cap\{y:\ \langle y,z\rangle\leq s(z,0)\}$.  For any $r>R_+$ we can choose a unique sphere of radius $r$ which encloses the latter region and intersects the hyperplane
$\{\langle y,z\rangle=s(z,0)\}$ in the same set as $S_{R_+}(0)$.  The centre of this sphere is at the point $O-\left(\sqrt{r^2+s^2-R_+^2}-s\right)z$, where $s=s(z,0)$.  By the comparison principle, $M_t$ is contained in the sphere with the same centre of radius 
$r(t)=\left(r^{1+\alpha}-(1+\alpha)t\right)^{1/(1+\alpha)}\leq r-r^{-\alpha}t$, and so
\begin{align*}
s(z,t)- s(z,0)&\leq-\sqrt{r^2+s^2-R_+^2}+r-r^{-\alpha}t\\
&=\frac{R_+^2-s^2}{r+\sqrt{r^2+s^2-R_+^2}}-r^{-\alpha}t\\
&\leq \frac{R_+^2}{r}-r^{-\alpha}t.
\end{align*}
For $0<t<\frac{R_+^{1+\alpha}}{\alpha}$ we choose $r=\left(\frac{R_+^2}{\alpha t}\right)^{\frac{1}{1-\alpha}}$, yielding
$$
s(z,t)- s(z,0)\leq -(1-\alpha)\alpha^{\frac{\alpha}{1-\alpha}}R_+^{-\frac{2\alpha}{1-\alpha}}t^\frac{1}{1-\alpha}.
$$
\qed\end{proof}

As a consequence of this result, an initial hypersurface containing a flat side (that is, an open subset contained in a supporting hyperplane) must immediately move.  We will see in Section \ref{S:posspeed}
that this implies that the speed of motion immediately becomes positive, so that the flat side becomes curved. 
Next we investigate the situation in the case $\alpha>1$.  Here we have the reverse situation, and
flat sides always persist:

\begin{theorem}\label{thm:alphalarge}
Suppose $\Omega_0$ be a bounded open convex region, and $M_0=\partial\Omega_0$.  Suppose that for some $z\in S^n$, $M_0$ has a flat side in direction $z$.  That is, there exists $y_0\in M_0$ with $\langle y_0,z\rangle=h(z)=\sup\{\langle y,z\rangle:\ y\in\Omega_0\}$ and $r>0$ such that the hemisphere $B_r^{n+1}(y_0)\cap\{\langle y,z\rangle<h(z)\}$ is contained in $\Omega_0$.   Then if $\{M_t=\partial\Omega_t\}_{0\leq t\leq T}$ is any generalized solution of the flow \eqref{E:theflow} satisfying Conditions \ref{T:Fconds} with $\alpha>1$, then
$y_0\in M_t$ for $0\leq t\leq C(\alpha)r^{1+\alpha}$.
\end{theorem}

%%%%% pictures
\begin{proof}
We construct a self-similar subsolution of \eqref{E:theflow} with a flat side, similar to the example in \cite{A9}*{Theorem II1.13}.   Define $A=\{(x,y)\in\RR\times\RR^n:\ x>0,\ G(x,y)<0\}$ where $G(x,y) = |y|^2+2x-1-2x^{\frac{\alpha-1}{2\alpha-1}}$.  $G$ is a convex function, so $A$ is a convex set which is rotationally symmetric about the $x$ axis.  Furthermore $A$ has $C^\infty$ boundary in $\{x>0\}$ with strictly positive principal curvatures, and convexity guarantees that $s=\langle X,\nu\rangle$ is also strictly positive on $\partial A\cap\{x>0\}$.  

The boundary of $A$ contains the disk $\{(0,y):\ |y|\leq 1\}$.  We look in more detail at the behaviour of $\partial A$ near the boundary curve $\{(0,y):\ |y|=1\}$.  For $|y|>1$ near this curve the curvature in the radial direction is given by
$$
\kappa_r = \frac{x_{yy}}{(1+x_y^2)^{3/2}}=\frac{\alpha(2\alpha-1)}{(\alpha-1)^2}(|y|-1)^{\frac{1}{\alpha-1}}(1+o(1))\quad\text{\rm as\ }|y|\to 1;
$$
so in particular the boundary of $A$ is $C^{2+\frac{1}{\alpha-1}}$.  The curvatures in the orthogonal directions satisfy
$$
\kappa_\perp = \frac{x_y}{|y|\sqrt{1+x_y^2}}=\frac{2\alpha-1}{\alpha-1}\left(|y|-1\right)^{1+\frac{1}{\alpha-1}}\left(1+o(1)\right)\quad\text{\rm as\ }|y|\to 1.
$$
Also we can compute
\begin{align*}
s &= \frac{2\alpha-1}{\alpha-1}(|y|-1)^{\frac{\alpha}{\alpha-1}}\left(1+o(1)\right)\\
&=\frac{2\alpha-1}{\alpha-1}\left(\frac{(\alpha-1)^2}{\alpha(2\alpha-1)}\right)^\alpha\kappa_r^\alpha\left(1+o(1)\right).
\end{align*}
Since $s$ is strictly positive and $\kappa_r$ and $\kappa_\perp$ are bounded away from $\{x=0\}$, there exists $\beta>0$ such that 
$$
s\geq \beta\max\{\kappa_r, \kappa_\perp\}^\alpha\geq \beta S
$$
everywhere on $\partial A$ (since $f(1,\dots,1)=1$, so that $S=f^\alpha\leq \kappa_{\max}^\alpha$ by the monotonicity and homogeneity of $f$).  It follows that $\left(r_0^{1+\alpha}-\frac{1+\alpha}{\beta}t\right)^{\frac{1}{1+\alpha}}A$ is a subsolution of \eqref{E:theflow}. 

\begin{figure}\hskip 1.5 cm
\includegraphics[scale=0.4]{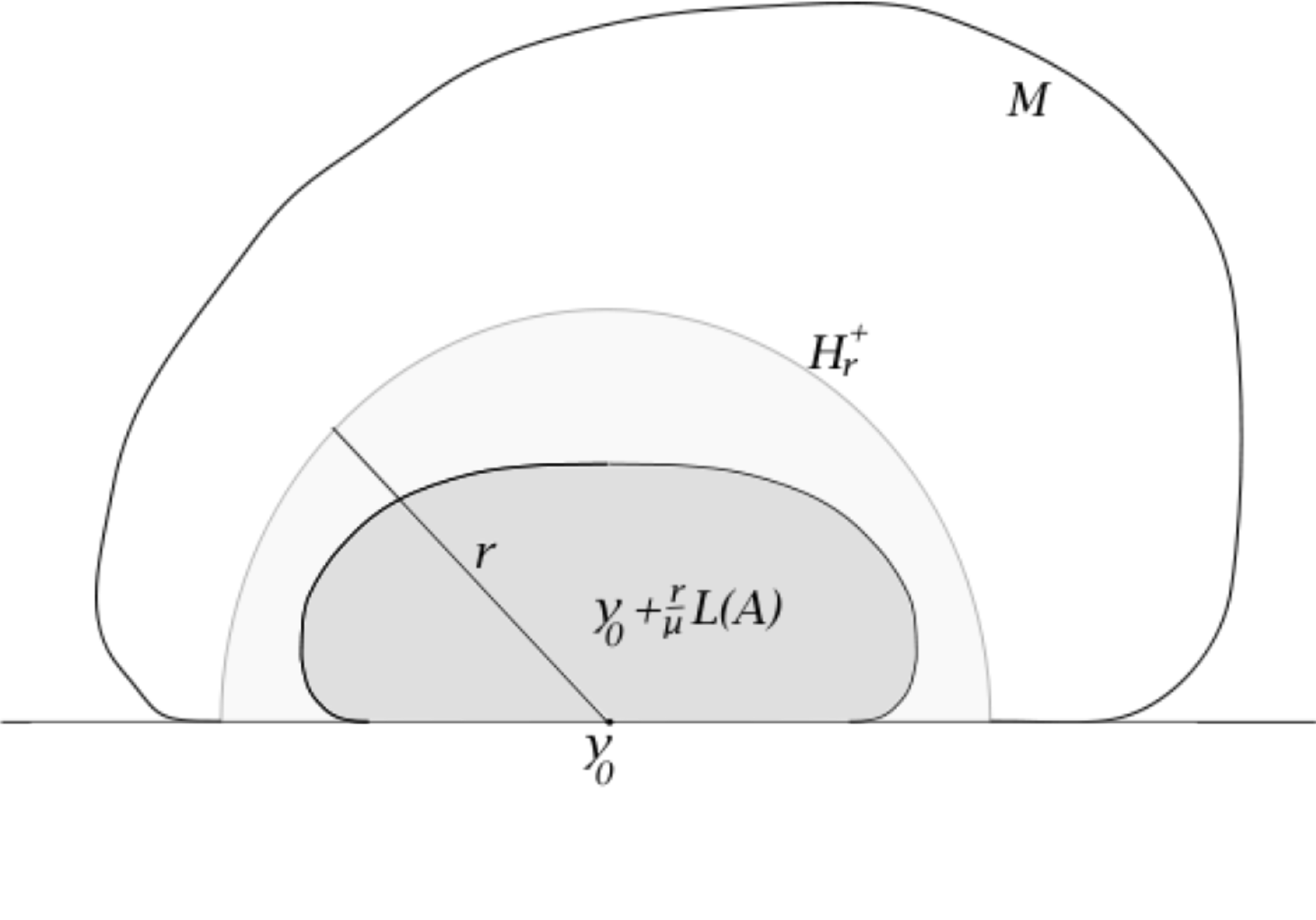}
\caption{A weakly convex hypersurface, containing a hemisphere $H_r=B_r^{n+1}(y_0)\cap\{\langle y,z\rangle<h(z)\}$ and contained in the corresponding halfspace $\{\langle y,z \rangle \leq h(z)\}$.  Inside the hemisphere we place a homothetically shrinking barrier $y_0+\frac{r(t)}{\mu}L(A)$, showing that the flat side on the hypersurface cannot move.}\label{fig:flat_side} 
\end{figure}

Let $\mu>0$ such that $A\subset B_\mu(0)$ (for example, $\mu=4$ suffices for any $\alpha>1$), and let $L$ be a rotation which takes $(1,0)$ to $-z$ and $\RR^n$ to $\{\langle y,z\rangle = 0\}$.  Then 
$y_0+\frac{r}{\mu}L\left(A\right)\subset \Omega_0$.  The comparison principle then gives 
that $y_0+\left(\left(\frac{r}{\mu}\right)^{1+\alpha}-\frac{1+\alpha}{\beta}t\right)^{\frac{1}{1+\alpha}}L(A)\subset \Omega_t$, and hence $y_0\in M_t$, for $0\leq t\leq \frac{\beta}{1+\alpha}\left(\frac{r}{\mu}\right)^{1+\alpha}$ while the solution exists.
\smartqed\end{proof}

The remaining case is $\alpha=1$, which will be of most interest in the present paper.  In this case there is a simple dichotomy between those flows for which flat sides persist and those where there is a lower bound on the displacement for positive times, and we provide a very simple criterion which distinguishes the two:

\begin{theorem}\label{thm:alpha1flatsides}
If $\alpha=1$, then flat sides move under \eqref{E:theflow} if and only if $f_*=0$ on the boundary of the positive cone.  That is, 
\begin{enumerate}
\item If part \ref{cond:invvanishbdry} of Conditions \ref{cond:extra} holds, then under \eqref{E:theflow} we have for all $t>0$ in the interval of existence and all $z\in\Sn^n$
\begin{equation}\label{eq:HD1dlb}
s(z,t)\leq s(z,0)-R_+3^{-\chi\left(\frac{3R_+^2}{t}\right)},
\end{equation}
where 
$\chi(\xi)$ is the inverse function of the function $\hat f$ defined by $\hat f(x) = f(x,1,\dots,1)$.  
\item If part \ref{cond:invvanishbdry} of Conditions \ref{cond:extra} does not hold, there exists $C>0$ such that if $\{M_t\}_{0\leq t\leq T}$ is a generalised solution of \eqref{E:theflow}, and $y_0\in M_0$ and $z\in S^n$ with $\langle y_0,z\rangle=s_0(z)=\sup\{\langle y,z\rangle:\ y\in\Omega_0\}$, and the hemisphere $B_r(y_0)\cap\{\langle y,z\rangle<s_0(z)\}$ is contained in $\Omega_0$ for some $r>0$, then $y_0\in M_t$ for $0\leq t\leq \min\{Cr^2,T\}$.
\end{enumerate}
\end{theorem}

\begin{proof}
We begin with the first case.  We will construct barriers which force the hypersurface to move, as follows:
Fix $z\in S^n$, and rotate so that $e_{n+1}=z$.   We construct graphical barriers, of the form 
$\{(x,y)\in\RR^n\times\RR:\ y=u(|x|,t)\}$.  We require these to be supersolutions of the flow, so that
$$
\frac{\partial u}{\partial t}\geq -\sqrt{1+(u')^2}f\left(\frac{-u''}{(1+(u')^2)^{3/2}},\frac{-u'}{|x|\sqrt{1+(u')^2}},\dots,\frac{-u'}{|x|\sqrt{1+(u')^2}}\right),
$$
and also that they be of the form $u(|x|,t)=w(|x|+vt)$.  This gives $\frac{\partial u}{\partial t} = vw'$, and
the required inequality becomes
$$
vw'\geq \frac{w'}{|x|}f\left(\frac{|x|w''}{w'(1+(w')^2)},1,\dots,1\right) = \frac{w'}{|x|}\hat f\left(\frac{|x|w''}{w'(1+(w')^2)}\right),
$$
where $\hat f(x):=f(x,1,\dots,1)$.  The assumption that $f_*=0$ on $\partial\Gamma_+$ is equivalent to 
the assumption that $\lim_{x\to\infty}\hat f(x)=\infty$, since $\hat f(x,1,\dots,1)=f_*(x^{-1},1,\dots,1)^{-1}$ (the implication in the other direction follows since if $f_*\neq 0$ anywhere on the boundary of the positive cone then $f_*(0,1,\dots,1)\neq 0$ by monotonicity).
We seek solutions which are concave, so that we have $w''\leq 0$ and $w'\leq 0$.  Thus the inequality
is equivalent to
$$
v|x|\leq \hat f\left(\frac{|x|w''}{w'(1+(w')^2)}\right).
$$
Choosing the origin to be at the centre of a ball of radius $R_+$ containing $M_0$, we have that for all $t>0$, $M_t$ is contained in the half-cylinder $\{(x,y):\ |x|\leq R_+,\ y\leq s(z,0)\}$.  We construct our graphical solution on the domain $\{|x|\leq 3R_+\}$, so that $v|x|\leq 3vR_+$, and so it suffices to find a solution of
$$
3vR_+\leq \hat f\left(\frac{|x|w''}{w'(1+(w')^2)}\right).
$$
Since $\hat f$ is increasing, there is a well-defined inverse function and we can write equivalently
$$
\frac{w''}{w'(1+(w')^2)} \geq \frac{1}{|x|}\hat f^{-1}\left(3vR_+\right).
$$
An explicit solution in the equality case with $w(R_+)=s(z,0)$ is given by
\begin{align*}
w(r) &= s(z,0)-3R_+\int_{1/3}^{r/(3R_+)} \frac{\xi^{\sigma}}{\sqrt{1-\xi^{2\sigma}}}\,d\xi\\
&\leq s(z,0)+\frac{3R_+}{1+\sigma}\left(\left(\frac13\right)^{1+\sigma}-\left(\frac{r}{3R_+}\right)^{1+\sigma}\right),
\end{align*}
where $\sigma=\hat f^{-1}(3vR_+)$.  Note that in the cylinder $\{|x|\leq R_+\}$ these solutions lie above the hyperplane $y = s(z,0)$, and so above the initial hypersurface $M_0$.  By the comparison principle, we have that $M_t$ lies below the hypersurface
$$
\{(x,y):\ y=w(|x|+vt),\ 0\leq |x|\leq 3R_+-vt\}
$$
(note that the boundary of this hypersurface lies outside the cylinder and so cannot meet $M_t$)
For fixed $t>0$ we choose $v=\frac{2R_+}{t}$, yielding that $M_t$ lies below
$$
\{(x,y):\ y=w(|x|+2R_+)\leq w(2R_+),\ 0\leq|x|\leq R_+\}.
$$
In particular this gives
$$
s(z,t)-s(z,0)\leq
\frac{3R_+}{1+\sigma}\left(\left(\frac{1}{3}\right)^{1+\sigma}-\left(\frac23\right)^{1+\sigma}\right)
\leq -R_+3^{-\hat f^{-1}\left(\frac{3R_+^2}{t}\right)}.
$$
This establishes \eqref{eq:HD1dlb} since $\chi=\hat f^{-1}$.

Next we consider the second case.  As in Theorem \ref{thm:alphalarge}, we will prove the existence of a self-similar subsolution with a flat side.  We use the set $A$ constructed in the proof of Theorem \ref{thm:alphalarge} (take $\alpha=2$ for concreteness).   Then we have 
$$
\kappa_\perp = 3(|y|-1)^2(1+o(1))\quad\text{\rm as\ }|y|\to 1,
$$
and 
$$
s = 3(|y|-1)^2(1+o(1))\quad\text{\rm as\ }|y|\to 1.
$$
It follows that there exists $\beta>0$ such that
$$
s\geq \beta \kappa_{\perp}
$$
everywhere on $\partial A$.  Let $f_\infty=\lim_{x\to\infty}\hat f(x)<\infty$.  Now we observe that for $x>0$,
$$
S=f(\kappa_r,\kappa_{\perp},\dots\kappa_\perp)
=\kappa_\perp f\left(\frac{\kappa_r}{\kappa_\perp},1,\dots,1\right) = \kappa_\perp \hat f\left(\frac{\kappa_r}{\kappa_\perp}\right)\leq f_\infty\kappa_\perp,
$$
and therefore $s\geq \frac{\beta}{f_\infty} S$ everywhere on $\partial A$.

It follows that $\left(r_0^2-\frac{2f_\infty}{\beta}t\right)^{1/2}A$ is a subsolution of the flow \eqref{E:theflow}, and the result follows as in Theorem \ref{thm:alphalarge}.
\qed\end{proof}

\section{Upper speed bounds}\label{S:boundspeed}

To obtain an upper speed bound, we use methods similar to those in \cite{A1}, \cite{A2}, \cite{A5}, \cite{M1} and \cite{M2}, which originate from the work of Chou on Gauss curvature flow \cite{T}.

\begin{theorem} \label{T:usb}
Let $f$ be such that $\sum_i\dot f^i\kappa_i^2\geq C_1f^2$ on $\Gamma_+$.  Then there exists $C_2$ depending on $\alpha$ and $C_1$ such that
for any smooth, uniformly convex solution $\left\{ M_{t} \right\}_{\left[ 0, T\right)}$ of \eqref{E:theflow} with $R_{-}\leq r_{-} \left( M_{t} \right) \leq r_{+}\left( M_{t}\right) \leq R_{+}$,
$$
S\left( x, t\right) \leq   C_2\frac{R_{+}}{R_{-}}  \left( R_-^{-\alpha} + t ^{-\frac{\alpha}{1+\alpha}}\right).$$
\end{theorem}

%% check C1 in proof
\begin{proof}
We use the method of Chou \cite{T} discussed previously in the proof of Theorem \ref{thm:contract.point}:
Consider the function $Q= \frac{S}{s -\delta}$, for a small constant $\delta$ to be chosen.  Its evolution equation is
\begin{align*}
\frac{\partial}{\partial t} Q &= \mathcal{L}Q + \frac{2}{\left( s-\delta \right)} \dot{S}^{kl} \nabla_{k} s \nabla_{l} Q -\frac{\delta } {\left( s-\delta \right)} \dot{S}^{kl} h_{km}{h^{m}}_{l}Q
+ \left(1+\alpha \right) Q^{2}\\
&\leq  \mathcal{L}Q+ \frac{2}{\left( s-\delta \right)} \dot{S}^{kl} \nabla_{k} s \nabla_{l} Q
-C_1\alpha\delta \left( s-\delta \right)^{1/\alpha}Q^{2+1/\alpha}+ \left( 1+\alpha \right) Q^{2}
\end{align*}
since $\dot S^{kl}h_k{}^ph_{pl} = \alpha f^{\alpha-1}\dot f^i\kappa_i^2\geq C_1\alpha S^{\frac{1+\alpha}{\alpha}}$.  Taking $\delta = \frac{R_{-}}{2}$, we find 
\begin{equation*} 
    \frac{\partial}{\partial t} Q \leq \mathcal{L} Q + \frac{2}{\left( s-\delta \right)} \dot{S}^{kl} \nabla_{k} s \nabla_{l} Q  + \left(1+\alpha - \alpha\delta^{1+1/\alpha} Q^{1/\alpha}\right) Q^{2} \mbox{,}
\end{equation*}
from which it follows from the maximum principle that
$$
\max_{\mathbb{S}^{n}} Q\left( \cdot, t\right) \leq \max \left\{ \left(\frac{2(1+\alpha)}{\alpha}\right)^\alpha\delta^{-(1+\alpha)}, \left(\frac{2}{1+\alpha}\right)^\frac{\alpha}{1+\alpha}\delta^{-1}t^{-\frac{\alpha}{1+\alpha}} \right\}.
$$
The result now follows from the definition of $Q$ and the fact that $s\leq 2 R_{+}$.  \smartqed\end{proof}

\begin{remark}\label{rmk:invconc}
The upper speed bound can also be deduced from the upper displacement bound using the Harnack estimate (Theorem \ref{T:Harnack}) as in \cite{A9}*{Lemma II2.1}.  The proof above is more robust since the Harnack inequality requires inverse-concavity.  The assumption on $\dot f^i\kappa_i^2$ is always satisfied (with $C_1=1$) in the inverse-concave case by Lemma \ref{lem:lowerdf2}.
\end{remark}

\section{Lower speed bounds}\label{S:posspeed}

In this section we translate the bounds on displacement into bounds above and below on the speed for positive times.  In \cite{A9} and \cite{A2} estimates of this kind were achieved by combining the lower displacement bounds with the Harnack estimate \cite{A6}*{Theorem 5.6}, which requires that the speed be inverse-concave:

\begin{theorem} \label{T:Harnack}
  For any smooth, uniformly convex solution of \eqref{eq:evolsuppfn} on $\mathbb{S}^{n}\times \left( 0, T\right)$,
  $$\frac{\partial}{\partial t} \left( \Psi t^\frac{\alpha}{1+\alpha} \right) \leq 0 \mbox{.}$$
\end{theorem}

Following \cite{CDS}*{Theorem 2.2} and \cite{AM}*{Lemma 12.2} we here deduce the lower speed bounds using an estimate first proved by Smoczyk \cite{Smoc}*{Proposition 4} for the mean curvature flow.  The result is that lower speed bounds hold precisely when lower displacement bounds do:

\begin{theorem}\label{T:lsb}
For any flow of the form \eqref{E:theflow} satisfying Conditions \ref{T:Fconds} with $\alpha>0$, and any $z\in\Sn^n$ and $t_2>t_1\geq 0$,
$$
-\Psi(z,t_2)\geq \frac{s(z,t_1)-s(z,t_2)}{(1+\alpha)(t_2-t_1)}.
$$
\end{theorem}

Recall that $\Psi=-S$, so Theorem \ref{T:lsb} amounts to a lower bound on the speed $S$ in any direction in which the support function has strictly decreased.

\begin{proof}
The evolution equation \eqref{eq:evolsuppfn}  for the support function $s$ 
can be rewritten in the form
$$
\frac{\partial s}{\partial t} = \dot\Psi^{ij}\bar\nabla_i\bar\nabla_js+s\dot\Psi^{ij}\bar g_{ij} +(1+\alpha)\Psi.
$$
Combining this with the evolution equation for $\Psi$ in Lemma \ref{T:evlnr}, we find:
$$
\frac{\partial}{\partial t}\left(s-(1+\alpha)(t-t_1)\Psi\right)
= \overline{\mathcal L}\left(s-(1+\alpha)(t-t_1)\Psi\right)+\left(s-(1+\alpha)(t-t_1)\Psi\right)\dot\Psi^{ij}\bar g_{ij}.
$$
It follows from the maximum principle that if $s-(1+\alpha)(t-t_1)\Psi\geq 0$ for all $z\in\Sn^n$ at $t=t_1$, then 
$s-(1+\alpha)(t-t_1)\Psi\geq 0$ for all $z\in\Sn^n$ at $t=t_2>t_1$.  In particular, if we fix $z\in \Sn^n$ and choose the origin to be at $X(z,t_1)$, then $s(z,t_1)=0$ while $s(z',t_1)\geq 0$ for all $z'$, so the inequality holds at $t=t_1$, and therefore we have
$$
0\leq s(z,t_2)-(1+\alpha)(t_2-t_1)\Psi(z,t_2) = s(z,t_2)-s(z,t_1)-(1+\alpha)\Psi(z,t_2),
$$
which gives the result of the Theorem.\qed
\end{proof}

\begin{corollary}\label{cor:lower.speed.bound}
If $0<\alpha<1$ then
$$
S\geq \frac{1-\alpha}{1+\alpha}\alpha^{\frac{\alpha}{1-\alpha}}R_+^{-\frac{2\alpha}{1-\alpha}}\min\left\{t,\alpha^{-1}R_+^{1+\alpha}\right\}^{\frac{\alpha}{1-\alpha}}.
$$
If $\alpha=1$ and $f_*=0$ on $\partial\Gamma_+$ then 
$$
S\geq \frac{R_+}{2t}3^{-\chi\left(\frac{3R_+^2}{t}\right)}.
$$
\end{corollary}

\begin{proof}
This follows immediately from Theorem \ref{T:lsb} using the displacement estimates proved in Theorem \ref{thm:alphasmall} and Theorem  \ref{thm:alpha1flatsides}.
\qed\end{proof}

\section{Existence and uniqueness of barrier solutions}\label{sec:barrier}

This section explores the existence and uniqueness of barrier solutions, and their dependence on the initial data.  The results are largely consequences of the existence theorems and displacement bounds proved in the previous sections. 

When working with convex regions it is natural to work with the Hausdorff distance, but in a more general setting it is useful to use a stronger metric which is suited to the fact that our notion of barrier solution involves inclusion of sets but is intended as an evolution of the \emph{boundaries} of these sets:  For closed sets $U$ and $V$ we define
$$
\dHaus(U,V) = \max\left\{\Haus(U,V),\Haus(\partial U,\partial V)\right\}.
$$

\begin{lemma}\label{lem:convhaus}
Suppose $U$ and $V$ are compact convex sets.  Then $\dHaus(U,V)=\Haus(U,V)$.  Furthermore if $\Omega_1$ and $\Omega_2$ are closed sets such that $U\subset\Omega_i\subset V$ for $i=1,2$, then $\dHaus(\Omega_1,\Omega_2)\leq \Haus(U,V)$.
\end{lemma}

\begin{proof}
$\Haus(\partial U,\partial V)$ is equal to the larger of $\sup_{x\in\partial U}d(x,\partial V)$ and $\sup_{y\in\partial V}d(y,\partial U)$, where $d(x,C)=\inf\{d(x,y):\ y\in C\}$.   Let $x\in\partial U$.  By compactness $d(x,\partial V)=d(x,y)$ for some $y\in\partial V$.  We consider two cases:  If $x\notin V$ then $y$ also attains the minimum distance from $x$ to $V$, so $d(x,\partial V)\leq d(x,V)\leq \Haus(U,V)$.  If $x\in V$, then let $H=\{z:\ z\cdot v\leq k\}$ be any supporting hyperplane for $U$ at $x$, so that $x\cdot v=k$ and $U\subset H$.  Let $w$ be a point in $V$ which maximizes $w\cdot v$.  Then $B_{d(x,y)}(x)\subset V$, so $w\cdot v\geq k+d(x,y)$.  It follows that $d(w,U)\geq d(x,y)$, and hence $d(x,y)\leq \Haus(U,V)$.  Thus $d(x,\partial V)\leq \Haus(U,V)$ for every $x\in\partial U$, so $\sup_{x\in\partial U}d(x,\partial V)\leq\Haus(U,V)$.  The same argument with $U$ and $V$ interchanged gives $\sup_{y\in\partial V}d(y,\partial U)\leq\Haus(U,V)$.  Therefore $\Haus(\partial U,\partial V)\leq\Haus(U,V)$ and $\dHaus(U,V)=\Haus(U,V)$ as required.

Now we prove the second statement:  We have $\sup_{x\in\Omega_1}d(x,\Omega_2)\leq\sup_{x\in V}d(x,U)\leq \Haus(U,V)$, and similarly $\sup_{x\in\Omega_2}d(x,\Omega_1)\leq\Haus(U,V)$, so $\Haus(\Omega_1,\Omega_2)\leq\Haus(U,V)$.  It remains to show $\Haus(\partial\Omega_1,\partial\Omega_2)\leq\Haus(U,V)$.  To see this, let $x\in\partial\Omega_1$, and let $z$ be the nearest point in $U$ to $x$.  If $x\neq z$ let $v=\frac{x-z}{|x-z|}$, and otherwise let $v$ be any unit normal vector to $U$ at $z$.  In either case $U\subset \{w:\ w\cdot v\leq z\cdot v\}$.  We note that $z+s\notin V$ for $s\geq\Haus(U,V)$, since $d(z+sv,U)=s$ for $s>0$.  Let $s_*=\sup\{s:\ z+sv\in\Omega_2\}\leq \Haus(U,V)$, and let $y=z+s_*v$.  Then $y\in\partial\Omega_2$, and $d(y,x)\leq \Haus(U,V)$ since $x$ and $y$ are both in the interval $\{z+sv:\ 0\leq s\leq\Haus(U,V)\}$.  Therefore $d(x,\partial\Omega_2)\leq\Haus(U,V)$, and since $x\in\partial\Omega_1$ is arbitrary, $\sup_{x\in\partial\Omega_1}d(x,\partial\Omega_2)\leq\Haus(U,V)$.  The same argument with $\Omega_1$ and $\Omega_2$ interchanged gives $\sup_{x\in\partial\Omega_2}d(x,\partial\Omega_1)\leq\Haus(U,V)$, and hence $\Haus(\partial\Omega_1,\partial\Omega_2)\leq\Haus(U,V)$ as required.
\qed\end{proof}

In comparing with non-convex sets the following lemma is useful:

\begin{lemma}\label{lem:hauscontain}
Let $\Omega$ be a compact convex set with nonempty interior, and suppose (by translating if necessary) that $r_-B\subset\Omega\subset r_+B$, where $B$ is the closed unit ball in $\RR^{n+1}$.  Then there exist constants $K$ and $L$ depending only on $r_-$ and $r_+$ such that the following hold:
\begin{enumerate}[label={(\roman*).}]
\item If $\Omega'$ is any closed set with $\dHaus(\Omega',\Omega)\leq d$ and $d< 1/K$, then 
$$
(1-Kd)\Omega\subset\Omega'\subset(1+Kd)\Omega.
$$
\item If $\Omega'$ is closed and $(1-a)\Omega\subset\Omega'\subset(1+a)\Omega$ for $a<1$, then $\dHaus(\Omega',\Omega)\leq La$.
\end{enumerate}
\end{lemma}

Note that $\Omega'$ need not be convex, though $\Omega$ must be.

\begin{proof}
\begin{enumerate}[label={(\roman*).}]
\item By assumption we have $\Haus(\Omega',\Omega)\leq d$, so $\Omega'\subset\Omega+dB\subset \Omega+\frac{d}{r_-}\Omega=\left(1+\frac{d}{r_-}\right)\Omega$, where the addition is in the sense of Minkowski, so that $U+V=\{x+y:\ x\in U,\ y\in V\}$.  The last equality used convexity of $\Omega$. 

To get the reverse inclusion we must work harder:  Assume that $d< \frac{r_-}{4}$, and 
let $x\in\left(1-\frac{4d}{r_-}\right)\Omega$.   Then $x+4dB\subset x+\frac{4d}{r_-}\Omega\subset \left(1-\frac{4d}{r_-}\right)\Omega+\frac{4d}{r_-}\Omega=\Omega$.  Therefore $d(x,\partial\Omega)\geq 4d$.  But we also have that $x\in\Omega\subset \Omega'+dB$, so there exists $y\in\Omega'$ and $z\in B$ such that $x=y+dz$.  We claim that $y+dB\subset\Omega'$:  Otherwise there exists $w\in y+dB$ with $w\in\partial\Omega'$.  But $\Haus(\partial\Omega',\partial\Omega)\leq d$, so there exists $q\in\partial\Omega$ with $d(q,w)\leq d$.  But then $d(x,\partial\Omega)\leq d(x,q)\leq d(x,y)+d(y,w)+d(w,q)\leq 3d$, which is impossible since we proved $d(x,\partial\Omega)\geq 4d$.  In particular we have $x\in y+dB\subset\Omega'$.  Since $x\in\left(1-\frac{4d}{r_-}\right)\Omega$ is arbitrary, we have $\left(1-\frac{4d}{r_-}\right)\Omega\subset\Omega'$,
completing the proof of (i) with $K=\frac{4}{r_-}$.
\item By assumption $\Omega'\subset(1+a)\Omega=\Omega+a\Omega\subset\Omega+ar_+B$.  Also, $\Omega=(1-a)\Omega+a\Omega\subset\Omega'+a\Omega\subset\Omega'+ar_+B$.  Therefore $\Haus(\Omega,\Omega')\leq ar_+$.

Now let $x\in \partial\Omega'$.   Since $(1-a)\Omega\subset\Omega'$ we have $(1-a)\text{\rm int}(\Omega)\subset\text{\rm int}(\Omega')$, so $x\in\Omega'\setminus(1-a)\text{\rm int}(\Omega)\subset (1+a)\Omega\setminus(1-a)\text{\rm int}(\Omega)$.  It follows that $x=(1+s)z$ for some $s\in[-a,a]$ and $z\in\partial\Omega$.  In particular $d(x,\partial\Omega)\leq |s||z|\leq ar_+$.  That is, we have $\partial\Omega'\subset\partial\Omega+ar_+B$.

Finally, let $x\in\partial\Omega$.  Then $(1-a)x\in(1-a)\Omega\subset\Omega'$, but $(1+a)x\notin\text{\rm int}(\Omega')$.  Therefore there exists $s\in[-a,a]$ such that $(1+s)x\in\partial\Omega'$.  Thus $d(x,\partial\Omega')\leq ar_+$, and we have proved $\partial\Omega\subset\partial\Omega'+ar_+B$.  

This proves (ii) with $L=r_+$.
\end{enumerate}
\qed\end{proof}

We note that such a result does not hold if $\dHaus(\Omega,\Omega')$ is replaced by the Hausdorff distance $\Haus(\Omega,\Omega')$:  For example we could take $\Omega'$ to be a finite set of points approximating $\Omega$ as closely as desired in Hausdorff distance, in which case $\Omega$ cannot be included in any scaled copy of $\Omega'$.

The proof of Theorem \ref{T:du} was given for classical solutions.  However, the argument involves nothing more than comparison with classical solutions (shrinking spheres), so we can apply the same argument to barrier solutions to give the following:

\begin{proposition}\label{prop:cty}
Let $\{\Omega_t\}_{a\leq t\leq b}$ be a barrier solution of \eqref{E:theflow}, such that $\Omega_a$ is a compact convex set with nonempty interior.  Then $\dHaus(\Omega_t,\Omega_a)\leq C(t-a)^{1/(1+\alpha)}$ for $a<t\leq b$, for some $C$ depending on the inradius and circumradius of $\Omega_a$.
\end{proposition}

Next we prove the following result concerning continuous dependence of barrier solutions on initial data:

\begin{theorem}\label{thm:ctsdep}
Let $\{\Omega_t\}_{0\leq t\leq T}$ be a barrier solution of equation \eqref{E:theflow}, such that $\Omega_t$ is convex and bounded with nonempty interior.  By a suitable translation we assume that there exists $0<r_-\leq r_+$ such that $r_-B\subset\Omega_t$ and $\Omega_t\subset r_+B$ for each $t\in[0,T]$.  Then there exist $C_1$ and $C_2$ depending only on $r_-$, $r_+$, $\alpha$ and $n$ such that for any other barrier solution $\{\Omega'_t\}_{0\leq t\leq T}$, the following holds:
$$
\sup_{0\leq t\leq T}\dHaus(\Omega'_t,\Omega_t)\leq C_1\dHaus(\Omega'_0,\Omega_0)^{\frac{1}{1+\alpha}}+C_2\dHaus(\Omega'_0,\Omega_0).
$$
\end{theorem}

\begin{proof} Let $d=\dHaus(\Omega_0,\Omega'_0)$.
By the previous Lemma, we have $(1-Kd)\Omega_0\subset\Omega_0'\subset(1+Kd)\Omega_0$ for some $K$ depending on $r_-$ and $r_+$.  The scaling invariance \eqref{E:scaling} of the flow enables us to write down the solutions of \eqref{E:theflow} with initial data given by the right and left sides of the above inclusion, in terms of the sets $\Omega_t$.  The comparison principle therefore gives
$$
(1-Kd)\Omega_{\left(1-Kd\right)^{-(1+\alpha)}t}\subset \Omega'_t
\subset(1+Kd)\Omega_{\left(1+Kd\right)^{-(1+\alpha)}t}.
$$
Since $(1-Kd)\Omega_0\subset\Omega_0\subset(1+Kd)\Omega_0$ we also have
$$
(1-Kd)\Omega_{\left(1-Kd\right)^{-(1+\alpha)}t}\subset \Omega_t
\subset(1+Kd)\Omega_{\left(1+Kd\right)^{-(1+\alpha)}t}.
$$
Therefore by Lemma \ref{lem:convhaus} we have 
\begin{align*}
\dHaus(\Omega_t,\Omega'_t)&\leq \Haus\left((1-Kd)\Omega_{\left(1-Kd\right)^{-(1+\alpha)}t},(1+Kd)\Omega_{\left(1+Kd\right)^{-(1+\alpha)}t}\right)\\
&\leq\Haus\left((1-Kd)\Omega_{\left(1-Kd\right)^{-(1+\alpha)}t},\Omega_{\left(1-Kd\right)^{-(1+\alpha)}t}\right)\\
&\quad\null+\Haus\left(\Omega_{\left(1-Kd\right)^{-(1+\alpha)}t},\Omega_{\left(1+Kd\right)^{-(1+\alpha)}t}\right)\\
&\quad\null+\Haus\left(\Omega_{\left(1+Kd\right)^{-(1+\alpha)}t},(1+Kd)\Omega_{\left(1+Kd\right)^{-(1+\alpha)}t}\right)\\
&\leq Kdr_++Kdr_++\Haus\left(\Omega_{\left(1-Kd\right)^{-(1+\alpha)}t},\Omega_{\left(1+Kd\right)^{-(1+\alpha)}t}\right).
\end{align*}
The last term is controlled by the displacement bound of Theorem \ref{T:du}, which gives 
\begin{align*}
\Haus\left(\Omega_{\left(1-Kd\right)^{-(1+\alpha)}t},\Omega_{\left(1+Kd\right)^{-(1+\alpha)}t}\right)
&\leq C\left((1+Kd)^{-(1+\alpha)}t-(1-Kd)^{-(1+\alpha)}t\right)^{\frac{1}{1+\alpha}}\\
&\leq \tilde Cd^{\frac{1}{1+\alpha}}t^{\frac{1}{1+\alpha}}.
\end{align*}
This completes the proof.
\qed
\end{proof}
 
A rather general uniqueness result follows immediately from Theorem \ref{thm:ctsdep}:  Convex barrier solutions are unique.

\begin{theorem}\label{thm:uniqbarrier}
Let $\{\Omega_t\}_{0\leq t<T}$ and $\{\Omega'_t\}_{0\leq t<T}$ be barrier solutions of \eqref{E:theflow}, and suppose that $\Omega_t$ is a compact convex set with nonempty interior, for each $t$, and $\Omega_0=\Omega'_0$.  Then $\Omega_t=\Omega'_t$ for each $t\in(0,T)$.
\end{theorem}

Note that we do not assume $\Omega'_t$ is convex.
Our next result is a compactness theorem for convex solutions, which we will employ in proving existence:

\begin{theorem}\label{thm:solncompact}
For each $k\in\NN$, suppose $\Omega^{(k)}\subset\RR^{n+1}\times[0,T]$ is a barrier solution of \eqref{E:theflow}, such that for each $k$ and each $t$, $\Omega^{(k)}_t=\{x:\ (x,t)\in\Omega^{(k)}\}$ is a compact convex set.  Assume that the inradii of $\Omega^{(k)}_t$ are bounded below by $r_->0$, and there exists $r_+$ such that $\Omega^{(k)}_t\subset B_{r_+}(0)$ for all $k$ and $t$.  Then there exists a subsequence $\Omega^{(k')}$ and a barrier solution $\Omega\subset\RR^{n+1}\times[0,T]$ of \eqref{E:theflow} with $\Omega_t=\{x:\ (x,t)\in\Omega\}$ convex for each $t$, such that $\lim_{k'\to\infty}\sup_{0\leq t\leq T}\Haus(\Omega^{(k')}_t,\Omega_t)=0$.
\end{theorem}

\begin{proof}
Passing to a subsequence, we can assume by the Blaschke selection theorem that $\Omega_0^{(k)}$ converges to a compact convex set $\Omega_0$.  Passing to a further subsequence we can assume $\dHaus(\Omega_0^{(k)},\Omega_0)\leq 2^{-k}$.  For a suitable choice of origin we have $B_{r_-}(0)\subset\Omega_0\subset B_{r_+}(0)$.  Since $\dHaus(\Omega_0^{(k)},\Omega_0^{(l)})\leq 2^{-k}+2^{-l}$ by the triangle inequality, Theorem \ref{thm:ctsdep} gives $\sup_{0\leq t\leq T}\dHaus(\Omega_t^{(k)},\Omega_t^{(l)})\leq C_1(2^{-k}+2^{-l})^{\frac{1}{1+\alpha}}+C_2(2^{-k}+2^{-l})$.  Therefore $\Omega_t^{(k)}$ is Cauchy in Hausdorff distance for each $t$, and hence converges to some compact convex set $\Omega_t$ with $\sup_{0\leq t\leq T}\dHaus(\Omega_t^{(k)},\Omega_t)\leq C_12^{-\frac{k}{1+\alpha}}+C_22^{-k}$.  It remains to prove that the limit is a barrier solution.  Suppose $\{A_t\}_{a\leq t\leq b}$ is a classical supersolution with $0\leq a<b\leq T$ and $\Omega_a\subset A_a$.  It follows from Lemma \ref{lem:hauscontain} that $(1-2^{-k}K)\Omega_a^{(k)}\subset\Omega_a\subset A_a$, so the barrier solution $(1-2^{-k}K)\Omega_{a+(1-2^{-k}K)^{-(1+\alpha)}(t-a)}^{(k)}$ remains inside $A_t$ for each $t$ where both are defined.  Letting $k\to\infty$ and using the continuity of $\Omega_t$ from Proposition \ref{prop:cty} we deduce that $\Omega_t\subset A_t$ for all $t\in[a,b]$.  The argument for comparison with classical subsolutions is similar.
\qed\end{proof}

Next we prove the main existence result:

\begin{theorem}\label{existbarrier}
Let $f$ satisfy the hypotheses of Theorem \ref{thm:contract.point}, \ref{thm:fconc} or \ref{thm:nis2}.  Then for any compact convex region $\Omega_0$ with nonempty interior in $\RR^{n+1}$ there exists a unique barrier solution  $\{\Omega_t\}_{0\leq t<T}$ of equation \eqref{E:theflow}.  The inradius of $\Omega_t$ approaches zero as $t\to T$.
\end{theorem}

\begin{proof}
Let  $\Omega_0^{(k)}$ be a sequence of smoothly, uniformly convex regions converging to $\Omega_0$ in Hausdorff distance. For each $k$ there is a smooth solution $\{M^{(k)}_t=\partial\Omega^{(k)}_t\}_{0\leq t<T_k}$ of equation \eqref{E:theflow} with initial data $\partial\Omega_0^{(k)}$, with the inradius of $M^{(k)}_t$ converging to zero as $t\to T_k$ (by Theorem  \ref{thm:contract.point}, \ref{thm:fconc} or \ref{thm:nis2}).  This gives a lower bound on $T_k$ independent of $k$:  If $B_{r_-}(p)$ is contained in the interior of $\Omega_0$, then $B_{r_-}(p)\subset\Omega_0^{(k)}$ for all sufficiently large $k$.  Therefore by comparison $B_{r(t)}(p)$ is contained in $\Omega^{(k)}_t$ for all large $k$ and all $t\in[0,t_-)$, where $r(t)^{1+\alpha}=r_-^{1+\alpha}-(1+\alpha)t$ and $t_-=\frac{r_-^{1+\alpha}}{1+\alpha}$.  Therefore $T_k\geq t_-$ for all large $k$.   By Theorem \ref{thm:solncompact} the solutions $\{\Omega^{(k)}_t\}_{0\leq t\leq t_-}$ have a subsequence which converges to a barrier solution $\{\Omega_t\}_{0\leq t\leq t_-}$ with the desired initial condition.  Since this solution can be extended as long as the inradius remains positive (and since the inradius is non-increasing in $t$) the inradius must converge to zero at the end of the maximal interval of existence.
\qed\end{proof}

\section{Examples:  Convex hypersurfaces which collapse to disks}\label{S:cylinders}

In this section we will investigate some interesting examples of curvature flows in which some weakly convex hypersurfaces do not become uniformly convex, even though the speed of motion becomes strictly positive.  Surprisingly, in these cases there are weakly convex regular $n$-dimensional hypersurfaces which do not even contract to points, instead shrinking to line segments, or more generally disks of dimension less than $n$.  These are the first such examples known, and for our purposes they serve to indicate the necessity of further conditions on the speed beyond those required to obtain positive speed in the previous sections.  We note that this behaviour can occur for arbitrary values of $\alpha$, in contrast to the situation for flat sides.

We are motivated by the following highly degenerate example:  If $f(x_1,\dots,x_n) = \max\{x_1,\dots,x_n\}$ and $\alpha>0$, then there is a $C^{1,1}$ generalized solution given by taking a cylinder of arbitrary length $L$ and radius $1$, and capping it with unit hemispheres.  This evolves to give the cylinder of the same length $L$ and radius $\left(1-(1+\alpha)t\right)^{1/(1+\alpha)}$, capped with semispheres of the same radius.   In particular the length of the cylinder does not decrease while the radius tends to zero.  Notice that this is just the Minkowski sum of the interval of length $L$ with a sphere of shrinking radius.   
Similarly, there are generalized solutions given by the Minkowski sum of a disk of dimension $k$ with a shrinking sphere, and these collapse onto the disk.
In the constructions below we modify these to find similar examples involving classical solutions for a class of flows which satisfy Conditions \ref{T:Fconds}, and indeed for which arbitrary smooth, uniformly convex initial hypersurfaces have smooth solutions which remain uniformly convex and contract to points.

Our main result is as follows:

\begin{theorem}\label{thm:cylind}
Let $\hat f(x)=f(x,1,\dots,1)$, and $g(x)=\hat f(x)^\alpha-\hat f(0)^\alpha$.  If $\int_0^1\frac{1}{g^{-1}(z)}\,dz<\infty$, then there exists $\tau>0$ such that if $L>2$ and $\{M_t=\partial\Omega_t\}_{0\leq t<T}$ is any barrier solution of Equation \eqref{E:theflow} such that $B^n_1(0)\times(-L,L)\subseteq \Omega_0\subseteq B^n_1(0)\times\RR$, then
$$
B^n_{r(t)}(0)\times (-(L-2),L-2)\subseteq \Omega_t\subseteq B^n_{r(t)}(0)\times\RR
$$
for $t\in [0,\tau]\cap[0,T)$, where $r(t)=\left(1-(1+\alpha)\hat f(0)^\alpha t\right)^{\frac{1}{1+\alpha}}$.
\end{theorem}

The condition in the theorem amounts to a certain degeneracy of the speed in directions of small curvature.  In particular, if $\hat f(x)^\alpha\leq \hat f(0)^\alpha + Cx^{1+\beta}$ with $\beta>0$, then the Theorem applies.  We give some examples below.

The Theorem says that if the initial hypersurface $M_0$ contains a cylindrical region, then for small positive times the hypersurface $M_t$ contains a slightly shorter cylindrical region (with smaller radius if $\hat f(0)>0$).  This has a surprising corollary:

\begin{corollary}\label{cor:collapsetoline}
If $\hat f(0)>0$ and $\int_0^1\frac{1}{g^{-1}(z)}\,dz<\infty$, then there exists $C>0$ such that if $L>C$ and
$$
B^n_1(0)\times(-L,L)\subseteq \Omega_0\subseteq B^n_1(0)\times\RR
$$
then $\{0\}\times[-(L-C),(L-C)]\subset \Omega_t$ for $0\leq t<T$.  In particular any barrier solution $\{M_t\}$ of \eqref{E:theflow} does not converge to a point as $t\to T$.
\end{corollary}

Before proving this result, we note some examples  where it applies:

\begin{example}\label{ex:cyl1}
Let $S=|A|^\alpha$, i.e. $f(x_1,\dots x_n)=\left(\sum_{i=1}^n x_i^2\right)^{1/2}$.  Then $\hat f(x) = (x^2+n-1)^{1/2}$, and $g(x) = (x^2+n-1)^{\alpha/2}-(n-1)^{\alpha/2}\simeq \frac{\alpha}{2}(n-1)^{\alpha/2-1}x^2$, so $\int_0^1\frac{dz}{g^{-1}(z)}<\infty$, and both Theorem \ref{thm:cylind} and Corollary \ref{cor:collapsetoline} apply, for any $\alpha>0$.  The same holds for $f(x_1,\dots,x_n) = \left(\sum_{i=1}^n x_i^p\right)^{1/p}$, if $\alpha>0$ and $p>1$ (i.e. powers of the power means $H_r$ for $r>1$).\end{example}

\begin{example}\label{ex:cyl2}
If $S = K^{\alpha/n}$, then $\hat f(x) = x^{1/n}$, so $g^{-1}(z) = z^{n/\alpha}$, and
by Theorem \ref{thm:cylind} cylindrical regions persist (in this case without shrinking) for $1-n/\alpha>0$, i.e. $\alpha>n$.  Note that Hamilton verified that cylindrical regions immediately vanish under Gauss curvature flow for $n=2$ (corresponding to the case $\alpha=n$) in \cite{Ha3}.  This was extended in \cite{A2}*{Theorem 23}, where it was proved that cylindrical regions persist under Gauss curvature flows if and only if $\alpha>n$.  It is known in these cases that the ratio of circumradius to inradius is bounded \cite{A2}*{Theorem 4}, so the behaviour described in Corollary \ref{cor:collapsetoline} is impossible here.
\end{example}

\begin{example}\label{ex:cyl3}
If $S=(K/E_k)^{\alpha/(n-k)}$, then we have $\hat f(x) = \left(\frac{x}{C_1x+C_2}\right)^{1/(n-k)}$ for some positive $C_1$ and $C_2$ depending on $n$ and $k$.  In particular $\hat f(x)$ is approximately $(x/C_2)^{1/(n-k)}$ for $x$ small, and cylindrical regions persist if $\alpha>n-k$.
\end{example}

As these examples illustrate, the dependence of the behaviour on $\alpha$ is evident only in the cases where $\hat f(0)=0$, and not in cases where $\hat f(0)>0$.

Now we proceed to the proofs.

\begin{proof}[Proof of Theorem \ref{thm:cylind}]
First assume $\hat f(0)>0$.
We will prove the existence of a homothetically shrinking subsolution of the flow with a cylindrical region.
The construction produces a rotationally symmetric hypersurface of the form $\phi(r,z) = (r,u(r)z)$ for $r\in I\subset\RR$ and $z\in S^{n-1}$.  For a homothetic subsolution with the same shrinking rate as the cylinder we require that 
\begin{equation}\label{eq:homsub}
\hat f(0)^\alpha\langle X,\nu\rangle \geq S.
\end{equation}
This amounts to the following equation for the function $u$:
\begin{align}
\hat f(0)^\alpha\left(u-ru'\right)&\geq \sqrt{1+(u')^2}f\left(\frac{-u''}{(1+(u')^2)^{3/2}},\frac{1}{u\sqrt{1+(u')^2}},\dots,
\frac{1}{u\sqrt{1+(u')^2}}\right)^\alpha\notag\\
&=u^{-\alpha}\left(1+(u')^2\right)^{\frac{1-\alpha}{2}} \hat f\left(\frac{-uu''}{1+(u')^2}\right)^\alpha\label{eq:homsubgraph}
\end{align}
where $u(r)=1$ for $|r|\leq r_0$ (this is the cylindrical region).  Note that the equation is certainly satisfied in this region.   In order to construct $u$ we first make some definitions:  We define
$$
G(z) = \int_0^z\frac{ds}{g^{-1}(s)}.
$$
Then by assumption $G$ is a continuous increasing concave function with $G(0)=0$, and has infinite gradient at $z=0$.  It follows that $G^{-1}$ is increasing, convex and continuously differentiable with $G^{-1}(0)=0$ and $(G^{-1})'(0)=0$.  We set
$$
v(x) = \begin{cases}
	1-\int_{0}^x G^{-1}(s)\,ds,&\quad x>0\\
	1,\quad x\leq 0.
	\end{cases}
$$
Then $v$ is a $C^2$ concave function satisfying
$$
v'(x) = -G^{-1}(r-r_0),
$$
and 
$$
v''(r) = -\frac{1}{G'(-v'(r))} = -g^{-1}(-v')
$$
for $r>0$, so that 
$$
g(-v'') = -v'
$$
and consequently
$$
\hat f(-v'')^\alpha = \hat f(0)^\alpha-v'.
$$
We restrict if necessary to a short interval $0<r<a$, so that $0<v(x)<\frac12$ and $|v'|\leq 1$.  Then we have since $\hat f$ is increasing that 
$$
\hat f\left(\frac{-vv''}{1+(v')^2}\right)^\alpha\leq \hat f(0)^\alpha-v'.
$$
Also, we have $v^{-1}(x)\leq 1-\int_{0}^x \frac{v'(s)}{v(s)^2}\,ds \leq 1-4\int_{r_0}^rv'(s)\,ds\leq 1-4v'(r)$,
and $(1+(v'(r))^2)\leq 1-v'(r)$, so on this interval we have
\begin{equation}\label{eq:vestimate}
v^{-\alpha}\left(1+(v')^2\right)^{\frac{\alpha-1}{2}}\hat f\left(\frac{-vv''}{1+(v')^2}\right)^\alpha
\leq \hat f(0)^\alpha\left(v-C_1v'\right)
\end{equation}
for some $C_1>0$.  
Now we extend $v$ beyond $a$ arbitrarily to a smooth, strictly concave function on $(0,b)$ with $v(x)\to 0$ as $x\to b$, and such that the hypersurface $\{X=(x,v(x)z):\ 0<x\leq b,\ z\in S^{n-1}\}$ is smooth and uniformly convex at the point $(b,0)$.  On the compact region of this hypersurface corresponding to $a\leq x\leq b$,  $\langle (1,0),\nu\rangle$ is strictly positive, and continuous, and hence on this region there exists a constant $C_2$ such that 
\begin{equation}\label{eq:compactest}
\frac{S-\hat f(0)^\alpha\langle X,\nu\rangle}{\langle (1,0),\nu\rangle }\leq C_2.
\end{equation}
Now we define our homothetic subsolution by setting $u(r)=v(r-r_0)$ where $r_0=\max\{C_1,C_2\}$.  Then on the interval $(r_0,r_0+a)$ we have by \eqref{eq:vestimate}
$$
u^{-\alpha}\left(1+(u')^2\right)^{\frac{\alpha-1}{2}}\hat f\left(\frac{-vv''}{1+(v')^2}\right)^\alpha
\leq \hat f(0)^\alpha\left(v-C_1v'\right)\leq \hat f(0)^\alpha\left(v-rv'\right),
$$
as required.  On the remaining interval we have by \eqref{eq:compactest}
\begin{align*}
S-\hat f(0)^\alpha\langle X,\nu\rangle &= S-\hat f(0)^\alpha\langle r_0(1,0)+(r-r_0,v(r-r_0),\nu\rangle\\
& \leq (C_2-r_0)\langle (1,0),\nu\rangle\\
&\leq 0,
\end{align*}
as required. Therefore with this choice of $r_0$ the hypersurface satisfies Equation \eqref{eq:homsub} everywhere.  Let $A_0=\{(x,y):\ |y|<u(|x|),\ |x|\leq r_0+b)\}$ be the region enclosed by this hypersurface.  Then the computation above shows that $A_t=(1-(1+\alpha)\hat f(0)^\alpha t)^{1/(1+\alpha)}A_0$ is a subsolution of Equation \eqref{E:theflow}, and hence by the comparison principle $A_t\subseteq M_t$ if $A_0\subseteq M_0$.  But we also have that the shrinking closed cylinder of radius $(1-(1+\alpha)\hat f(0)^\alpha t)^{1/(1+\alpha)}$ encloses $M_t$, and the result follows.

Next we consider the case where $\hat f(0)=0$. In this case it does not make sense to seek a homothetically shrinking solution since the cylindrical region should not move, and instead we attempt to construct a translating subsolution.  This means that we are seeking a solution of
\begin{equation}\label{eq:transsub}
S\leq V\langle (1,0),\nu\rangle
\end{equation}
where $V>0$ is the speed of motion.  In the graphical parametrization this becomes
\begin{equation}\label{eq:transsubgraph}
u^{-\alpha}\left(1+(u')^2\right)^{\frac{1-\alpha}{2}}\hat f\left(\frac{-uu''}{1+(u')^2}\right)^\alpha \leq -Vu'.
\end{equation}
The construction is similar to that used above:  We define $v(x)$ exactly as before, giving a solution of
$$
\hat f\left(\frac{-vv''}{1+(v')^2}\right)^\alpha\leq -v'.
$$
On an interval $(0,a)$ we have $u$ close to $1$ and $|u'|$ small, and hence we have
$$
u^{-\alpha}\left(1+(u')^2\right)^{\frac{1-\alpha}{2}}\hat f\left(\frac{-uu''}{1+(u')^2}\right)^\alpha \leq -V_1u'
$$
for some $V_1>0$.  Now as before we extend arbitrarily to a smooth, uniformly convex hypersurface, and observe that since $\langle (1,0),\nu\rangle$ is strictly positive on this region, the ratio $S/\langle (1,0),\nu\rangle$ is bounded by some constant $V_2$.  Taking the hypersurface to translate with speed $\max\{V_1,V_2\}$ then gives a translating subsolution of the flow, and the result again follows.
\smartqed\end{proof}

\begin{proof}[Proof of Corollary \ref{cor:collapsetoline}]
In this case the construction in the proof of Theorem \ref{thm:cylind} produces a homothetic subsolution $A_t$ containing a cylindrical part, which shrinks at the same rate as the cylinder of the same radius.  Precisely, $A_0$ is defined as $\{(x,y):\ |y|<u(x),\ |x|\leq r_0+b\}$, where $u$ is a concave function with $u(x)=1$ for $|x|\leq r_0$.  We take $C=r_0+b$.  If $L>C$ then $A_0+(L-C)(1,0)$ and $A_0-(L-C,0)$ are both contained in $M_0$.  The time-dependent set $B_t$ given for each $t$ by the convex hull of $r(t)A_0+(L-C)(1,0)$ and $r(t)A_0-(L-C)(1,0)$ is a subsolution of the flow if $r(t)=\left(1-(1+\alpha)\hat f(0)^\alpha\right)^{1/(1+\alpha)}$, so we have $B_t\subset M_t$ for $0\leq t\leq 1/((1+\alpha)\hat f(0)^\alpha)$.  On the other hand $M_t$ is contained in the cylinder $\{(x,y):\ |y|\leq r(t)\}$, so the maximal time of existence is no greater than $1/((1+\alpha)\hat f(0)^\alpha)$, and $(L-C,L+C)\times\{0\}\subset M_t$ for all $t$ in the interval of existence.
\smartqed\end{proof}

A similar construction gives the following result for hypersurfaces which have a region isometric to $B^k(0)\times S^{n-k}$:  In cases where such cylinders move we prove the existence of homothetic subsolutions, and in cases where the cylinders do not move we prove existence of subsolutions with translating radius, and deduce the following.

\begin{theorem}\label{thm:k-cylind}
Fix $k\in\{2,\dots,n-1\}$ and define for $x,p>0$
$$
f_k(x,p)=f(x,\underbrace{p,\dots,p}_{k-1\ \text{times}},1,\dots,1).
$$
Define 
$$
g_{k,p}(x) = f_k(x,p)^\alpha-f(\underbrace{0,\dots,0}_{k\ \text{times}},1,\dots,1)^\alpha.
$$
Note that $g_{k,p}$ is increasing for each $k$ and $p$, and so has a well-defined inverse $g_{k,p}^{-1}$.
If $\limsup_{p\to 0}\frac{f_k(0,p)^\alpha}{p}<\infty$ and $\int_0^1 \frac{dp}{g_{k,p}^{-1}(Vp)}<\infty$ for sufficiently large $V$, then $(n-k)$-cylindrical regions persist, in the following sense:  If $f_k(0,0)>0$ then there exists $\rho>0$ and $b>0$ such that if $\{M_t\}_{0\leq t\leq T}$ is a viscosity solution of \eqref{E:theflow}, such that $B^k_{(\rho+a)r_0+c}\times S^{n-k}_{r_0}\subset M_0\subset \RR^k\times \bar B^{n+1-k}_{r_0}$, then 
$$
B^k_{\rho r(t)+c}\times S_{r(t)}^{n-k}\subset M_t\subset \RR^k\times \bar B^{n+1-k}_{r(t)},
$$
where $r(t) = \left(r_0^{1+\alpha}-(1+\alpha)f_k(0,0)^\alpha t\right)^{1/(1+\alpha)}$, for $0\leq t<T$.  If $f_k(0,0)=0$ then there exists $a,b,V>0$ such that if
$B^k_R\times S^{n-k}_{r_0}\subset M_0\subset \RR^k\times \bar B^{n+1-k}_{r_0}$, then
$$
B^k_{R(t)}\times S^{n-k}_{r_0}\subset M_t\subset \RR^k\times\bar B^{n+1-k}_{r_0}
$$
where $R(t) = R-ar_0-Vtr_0^{-\alpha}$, as long as $R(t)\geq br_0$.
\end{theorem}

Again, the condition is related to degeneracy of $\dot f^i$ when $\kappa_i=0$  at points of the form $(0,\dots,0,1,\dots,1)$.  

\begin{example}\label{ex:cyl4}
Under the flow with $S = |A|^\alpha$ for any $\alpha>0$, hypersurfaces containing regions isometric to $B_R^k\times S^{n-k}_1$ have an enclosed $k$-disc of positive radius throughout their evolution if $R$ is sufficiently large, for any $k=1,\dots,n-1$.  The same holds for $S = H_r^\alpha$ for any $r>1$ and $\alpha>0$.
\end{example}

\begin{example}\label{ex:cyl5}
The result in the case $S=K^{\alpha/n}$ agrees with the conclusions in \cite{A2}:  Cylindrical regions of the form $B^k\times S^{n-k}$ persist if $\alpha>\frac{n}{k}$.  This is sharp, as proved in \cite{A2}*{Theorem 23}.
\end{example}

\begin{example}\label{ex:cyl6}
For flow by powers of other elementary symmetric functions of the principal curvatures ($S=E_m^{m/\alpha}$) the theorem implies that $k$-cylindrical regions persist provided $m+k>n$ and $\alpha>\frac{m}{m+k-n}$.  This agrees with the previous example in the case $m=n$; the case $m=1$ corresponds to powers of the mean curvature, and in this case $k$-cylindrical regions do not persist for any $k<n$.
\end{example}

\begin{example}\label{ex:cyl7}
Other interesting examples are $S=\left(\frac{E_m}{E_\ell}\right)^{\frac{\alpha}{m-l}}$ where $0\leq l<m\leq n$ and $E_m$ is the $m$th elementary symmetric function of the principal curvatures.  In this case the
Theorem implies that cylindrical regions of the form $B^k\times S^{n-k}$ persist when 
$k+m>n$ and $\alpha>\max\left\{1,\frac{m-\ell}{k+m-n}\right\}$.
\end{example}

\begin{example}\label{ex:cyl8}
Under the flow $S= H-\frac{E_{j+1}}{E_j}$, cylindrical regions of the form $B^k\times S^{n-k}$ persist if $k\geq n-j$.  Thus in the case $j=1$ hypersurfaces can collapse to $(n-1)$-dimensional disks, but not lower-dimensional ones; in the case $j=n-1$ then hypersurfaces can collapse to line segments, or disks of any dimension up to $n-1$.
\end{example}

\section{Examples:  Non-smooth convex hypersurfaces which fail to become smooth}\label{sec:nonsmooth}

In this section we give a construction similar to that of the previous section, which produces generalized solutions for a class of flows from weakly convex initial data, which fail to immediately become smooth.  This result is distinct from that in Section \ref{sec:lossreg}, since the examples constructed here are for flows which keep smooth, uniformly convex hypersurfaces smooth and uniformly convex until they contract to points (and indeed in some cases always have spherical limiting shape).  These examples arise when $f_*$ has degenerate derivative in directions orthogonal to the boundary of the positive cone.  We note that the
condition required to rule out the persistence of flat sides (i.e. $f_*=0$ on the boundary of $\Gamma_+$) excludes the examples we construct here.

\begin{theorem}  Let $\alpha=1$ and suppose $f_*(0,1,\dots,1)>0$, and $f$ satisfies the assumptions of Theorem  \ref{thm:contract.point}, \ref{thm:fconc} or \ref{thm:nis2}.
Let $\hat f_*(x)=f_*(x,1,\dots,1)$, and $g(x) = \hat f_*(0)^{-1}-\hat f_*(x)^{-1}$.  If $\int_0^1\frac{1}{g^{-1}(z)}\,dz<\infty$, then for any $L>0$ there exists $\tau>0$ such that if $\Omega_0$ is an open convex region satisfying
$$
B^n_1(0)\times\{0\}\subset \Omega_0\subset \text{\rm conv}(B_1^n(0)\times\{0\}, (0,L),(0,-L))
$$
then the generalized solution $\{\Omega_t\}$ of Equation \eqref{E:theflow} satisfies
$$
B^n_{r(t)}(0)\times\{0\}\subset \Omega_t\subset \text{\rm conv}(B^n_{r(t)}(0)\times\{0\},\ (0,2L),\ (0,-2L))
$$
for $0<t<\tau$, where $r(t)^2 = 1-2\hat f_*(0)^{-1}t$.  In particular $\Omega_t$ is Lipschitz but not $C^1$ for $t$ in this range.
\end{theorem}

Before embarking on the proof, we mention some examples to which the result applies:  If $f=H_r$ with $r<-1$, then the conditions of the theorem hold.  We note that in these cases $f$ is concave and zero on the boundary of the positive cone, and indeed arbitrary uniformly convex, smooth initial hypersurfaces remain smooth and uniformly convex until the shrink to points, becoming spherical in the process.  Thus the result
complements the second part of Theorem \ref{thm:alpha1flatsides}:  That result shows that flat parts of the hypersurface can persist; the present result shows that non-smoothness in the form of ridges in the surface can also persist for a more restricted class of flows.

\begin{proof}
We work with the evolution equation \eqref{eq:evolsuppfn} for the support function, in an axially symmetric situation.  Instead of working with the support function $s$ as a function on the unit sphere, it is convenient to work with a related function $\sigma$ defined on the unit cylinder $\RR\times S^{n-1}_1(0)\subset \RR\times\RR^n$, defined as follows:  We define $P:\ \RR\times S^{n-1}_1(0)\to S^n_1(0)$ by
$$
P(u,z) = \frac{(u,z)}{\sqrt{1+u^2}},
$$
and then set $\sigma(u,z) = s(P(u,z))\sqrt{1+u^2}$.  We note that the support function may be naturally considered as a homogeneous degree one function on $\RR^{n+1}$, and then $s$ is the restriction of this to $S^{n}$, while $\sigma$ is the restriction of the same function to the cylinder $\RR\times S^{n-1}$.  The axial symmetry assumption corresponds to the requirement that $\sigma$ be independent of $z$.

We compute the principal radii of curvature in terms of $\sigma$ and its derivatives, using the formula \eqref{E:r}.  For this purpose we derive formulae for the metric $\bar g$ and connection $\bar\nabla$ induced by $P$ onto $\RR\times S^{n-1}$.  We relate these to the metric $\tilde g$ and connection $\tilde\nabla$ on $S^{n-1}$, working in local coordinates in which $\partial_0=\frac{\partial}{\partial u}$ and $\partial_i$ is tangent to $S^{n-1}$ for $i=1,\dots,n-1$:  The metric is given by
\begin{align*}
\bar g_{00} &= \frac{1}{(1+u^2)^2};\\
\bar g_{0i}&=0;\\
\bar g_{ij} &= \frac{\tilde g_{ij}}{1+u^2};
\end{align*}
while the connection is given by
\begin{align*}
\bar\nabla_u\partial_u &= -\frac{2u}{1+u^2}\partial_u;\\
\bar\nabla_u\partial_i &= -\frac{u}{1+u^2}\partial_i;\\
\bar\nabla_i\partial_j &= \tilde\nabla_i\partial_j+u\tilde g_{ij}\partial_u.
\end{align*}
From these we find the following expressions for the components of the matrix $\rr$:
\begin{align*}
\rr_{00}&= \frac{\partial^2 s}{\partial u^2}-(\bar\nabla_u\partial_u)s+\bar g_{00}s = \frac{\sigma''}{\sqrt{1+u^2}};\\
\rr_{0i}&=0;\\
\rr_{ij}&=\tilde\nabla_i\tilde\nabla_js-(\bar\nabla_i\partial_j-\tilde\nabla_i\partial_j)s+\bar g_{ij}s=\frac{\sigma-u\sigma'}{\sqrt{1+u^2}}\tilde g_{ij}.
\end{align*}
Since the principal radii of curvature are the eigenvalues of $\rr$ with respect to $\bar g$, these are given by
$$
\rr_0=(1+u^2)^{3/2}\sigma'',\quad \rr_i=\sqrt{1+u^2}(\sigma-u\sigma'),\ i=1,\dots,n-1.
$$
From this and equation \eqref{eq:evolsuppfn} we derive the following evolution equation for $\sigma$:
\begin{align}\label{eq:evol.sigma.cylinder}
\frac{\partial\sigma}{\partial t} &= -\frac{\sqrt{1+u^2}}{f_*\left((1+u^2)^{3/2}\sigma'', \sqrt{1+u^2}(\sigma-u\sigma'),\dots,\sqrt{1+u^2}(\sigma-u\sigma')\right)}\notag\\
&=-\frac{1}{(\sigma-u\sigma')\hat f_*\left(\frac{(1+u^2)\sigma''}{\sigma-u\sigma'}\right)}.
\end{align}
The unit disc $B_1^{n}(0)\times\{0\}$ corresponds to $\sigma(u)=1$ for all $u$.  It follows that the set
$$
{\mathcal A}_-=\left\{(0,y,t)\in\RR\times\RR^n\times\left[0,\frac{\hat f_*(0)}{2}\right]:\ |y|\leq \sqrt{1-\frac{2t}{\hat f_*(0)}}\right\}
$$
is an inner barrier.

To prove the theorem we will combine this inner barrier with an outer barrier:  Precisely, we construct a homothetic supersolution $\sigma_+$  with the same rate of shrinking as the disc, and with $\sigma_+(u)=1$ for $|u|\leq u_0$.   Such a homothetic supersolution satisfies
$$
\frac{1}{(\sigma_+-u\sigma_+')\hat f_*\left(\frac{(1+u^2)\sigma_+''}{\sigma_+-u\sigma_+'}\right)}\geq\frac{\sigma_+}{\hat f_*(0)},
$$
or equivalently
\begin{equation}\label{eq:homothety.GMP}
g\left(\frac{(1+u^2)\sigma_+''}{\sigma_+-u\sigma_+'}\right)\leq\frac{1-\sigma_+(\sigma_+-u\sigma_+')}{\hat f_*(0)}.
\end{equation}
If $\sigma_+$ is convex with $\sigma_+(u)=1$ for $u\leq u_0$, then $\sigma_+\geq 1$ and $(\sigma_+-(u-u_0)\sigma_+')' = -(u-u_0)\sigma_+''\leq 0$ for $u\geq u_0$, so $\sigma_+-(u-u_0)\sigma_+'\leq 1$, or equivalently $(u-u_0)\sigma_+'\geq \sigma_+-1$.  It follows that $1-\sigma_+(\sigma_+-u\sigma_+')\geq u\sigma_+'+1-\sigma_+^2\geq u\sigma_+'-2(\sigma_+-1)\geq \frac12(u\sigma_+'-(\sigma_+-1))+(\sigma_+-1)\left(\frac{u}{2(u-u_0)}-\frac32\right)\geq \frac12(u\sigma_+'-(\sigma_+-1))$ provided $u_0< u\leq \frac32u_0$.  If we arrange that $\sigma_+-u\sigma_+'>0$, then for $u_0\leq u\leq \frac32u_0$ it suffices to find a solution of the inequality
\begin{equation}\label{eq:homothety2}
g\left(\frac{(1+u^2)\sigma_+''}{\sigma_+-u\sigma_+'}\right)\leq \frac1{2\hat f_*(0)}\left(u\sigma_+'-(\sigma_+-1)\right).
\end{equation}
A solution with equality can be found explicitly as follows: We introduce a new variable $v=u\sigma_+'-(\sigma_+-1)$.  Then
$v' = u\sigma_+''$, so the solution is given by
$$
\int_0^v\frac{dz}{(1-z)g^{-1}\left(\frac{z}{2\hat f_*(0)}\right)} = \frac1{2}\log\left(\frac{1+u^2}{1+u_0^2}\right).
$$
The factor $1/(1-z)$ is comparable to $1$ for $z$ small, so the integral on the left exists and defines $v$ as a $C^1$ function of $u$.  The identity $\left(\sigma_+/u\right)' = (1-v)/u$ may then be integrated to produce a $C^2$ function $\sigma_+$ on some interval $[0,u_1]$ with $\sigma_+(u)=1$ for $u\leq u_0$, where $0<u_0<u_1$.  We extend $\sigma_+$ to be even in $u$.

Now we use the barrier $\sigma_+$ to prove the theorem:    From the function $\sigma_+$ we construct a singular hypersurface from the image of a map $X$ from $S^{n-1}\times[-u_1,u_1]$ to $\RR^{n+1}$, defined by $X(z,s) = (\sigma_+',(\sigma_+(u)-u\sigma_+')z)$.  We have $X(u,z)=(0,z)$ for $|u|\leq u_0$, and $X$ is a smooth strictly locally convex embedding for $u_0<|u|\leq u_1$.  In particular $X(u,z)\cdot(1,0)=\sigma_+'(u)$ increases monotonically from $u=u_0$ to $u=u_1$, and the set 
$
A=\left\{(x,y)\in\RR\times\RR^n:\ x=\sigma_+'(u),\ |y|\leq \sigma_+(u)-u\sigma_+'(u),\ |u|\leq u_1\right\}
$
is convex.  From our construction we have that the set
$$
{\mathcal A}_+ = \left\{(x,y,t)\in\RR\times\RR^n\times\left[0,\frac{3\hat f_*(0)}{8}\right]:\ (x,y)\in\sqrt{1-\frac{2t}{\hat f_*(0)}}A,\ |x|\leq \frac{\sigma_+'(u_1)}{2}\right\}
$$
is an outer barrier to the flow.

Now under the assumption that $f$ satisfies the conditions of Theorems \ref{thm:contract.point}, \ref{thm:fconc} or \ref{thm:nis2}, Theorem \ref{existbarrier} guarantees the existence of a unique barrier solution for any initial data given by the boundary of an open convex region, which is therefore equal to the intersection of all outer barriers and to the union of all inner barriers.  In particular if the initial data contains $\mathcal{A}_-$ and is contained in the set $\mathcal{A}_+$ at $t=0$, then this remains true on the interval $\left[0,\frac{3\hat f_*(0)}{8} \right]$, and it follows that the solution has a persisting ridge of infinite curvature on this time interval.
\qed\end{proof}

There are some interesting questions which arise from the above construction:  We proved that the ridges of infinite curvature persist for some time, but we have so far not been able to construct examples where the ridges persist until the hypersurface shrinks to a point.  Can such examples occur?  And further, could there exist an example where the initial hypersurface encloses an open convex region, but the  limiting shape as the hypersurfaces contract to a point is the degenerate shrinking disc?

\section{Strict convexity and smoothness}\label{S:strictconv}

In this section we investigate which flows admit upper and lower bounds on principal curvatures for positive times.  As the examples of the previous section indicate, this is a much more subtle question than the upper and lower speed bounds.  We do not attempt to give a complete
characterization of when such bounds hold, but instead provide such estimates for several natural
classes of flows.

We remark that the examples constructed previously provide a number of necessary conditions:  We must have $\alpha\leq 1$ to avoid flat sides persisting, and if $\alpha=1$ then we must have $f_*=0$ on the boundary of the positive cone;  if $\alpha<1$ then we must also have $f_*=0$ on the boundary of the positive cone to avoid loss of smoothness, while if $\alpha=1$ then this requires $f_*$ to be inverse-concave on the boundary.  Finally, to avoid loss of convexity we must have $f$ inverse-concave on the boundary, and to avoid examples such as shrinking cylinders we must impose a certain non-degeracy near points on the boundary of the positive cone where $f$ is non-zero.

We begin by discussing the case $n=2$, which is somewhat simpler than the higher-dimensional cases.  

\begin{theorem}\label{thm:n=2becomesmooth}
Let $n=2$ and suppose $0<\alpha\leq 1$ and $f_*$ is concave.  Assume that either
\begin{enumerate}[label={(\roman*)}, ref={(\roman*)}]
\item\label{case:fpos} Condition \ref{cond:extendf} holds and $f>0$ on $\bar\Gamma_+\setminus\{0\}$; or
\item\label{case:f0} $f$ and $f_*$ both vanish on $\partial\Gamma_+$, and 
$\int_0^1\frac{1}{f_*(1,x)}dx<\infty$.
\end{enumerate}
Then the barrier solution with initial data given by the boundary of an bounded open convex region is smooth and uniformly convex for small positive times.
\end{theorem}

\begin{proof}
Let $\Omega_0$ be a bounded open convex region in $\RR^3$.  We begin by approximating $\partial\Omega_0$ by smooth, uniformly convex regions $\Omega_0^{(k)}$.  Theorem \ref{thm:contract.point} provides the existence of a solution $\{\Omega^{(k)}_t\}_{0\leq t\leq T_k}$ for each $k$, with $\Omega^{(k)}_t$ converging to a point as $t\to T_k$.  Comparison with an enclosed shrinking sphere gives a positive limit inferior $T$ of $T_k$ as $k\to\infty$.  The results of Section \ref{sec:barrier} imply that the regions $\Omega^{(k)}_t$ converge in the $\dHaus$ distance to the unique barrier solution $\{\Omega_t\}_{0\leq t<T}$ with initial data $\Omega_0$.  We proceed by obtaining estimates on the solutions $\Omega^{(k)}_t$ independent of $k$ for positive $t$, which are then inherited by the limit $\Omega_t$.

Since $f_*$ is zero on the boundary of $\Gamma_+$, Corollary \ref{cor:lower.speed.bound} provides a positive lower bound on the speed $S$ depending only on the inradius and circumradius of $\Omega^{(k)}_0$ and on $t$, for any small $t>0$.  Theorem \ref{T:usb} (with remark \ref{rmk:invconc}) provides an upper bound on $S$.   

We consider case \ref{case:fpos} first.  In this case $\dot f^i$ are uniformly bounded and have a positive lower bound on $\Gamma_+$, by the assumption \ref{cond:extendf}.

Now fix $\varepsilon>0$ small and work on the interval $[\varepsilon,T-\varepsilon]$, so that we have uniform bounds above and below on the speed $S$ and on the inradius and circumradius for the solutions $\Omega^{(k)}_t$.  
It follows that the flow (for example, in the form of an evolving local graph) is uniformly parabolic, and uniform bounds in $C^{2,\beta}$ on $[2\varepsilon,T-\varepsilon]$ follow from the results of \cite{Andrews2D}. 
Schauder estimates then imply uniform bounds in $C^k$ for every $k$ on $[3\varepsilon,T-\varepsilon]$.  It follows that the solution $M_t=\partial\Omega_t$ is a $C^\infty$ hypersurface.  The fact that $M_t$ is uniformly convex now follows from the strong maximum principle of Bony \cite{BonyMP} and Hill \cite{HillMP} as in \cite{A5}*{Proposition 10.1}.

Now consider case \ref{case:f0}:   As before we can work on a time interval where we have uniform upper and lower bounds on the speed, so that $S_-\leq S\leq S_+$ for some positive constants $S_-$ and $s_+$.
The idea is to use the integral assumption to deduce that some components of $\dot f_*$ become large near the boundary of the positive cone, and this will allow us to prove an upper bound on $S_1$. 
The evolution equation for $S_1$ is obtained from a trace of the evolution equation in Lemma \ref{T:evlnr}:
$$
\frac{\partial S_1}{\partial t} = \overline{\mathcal L}S_1+\bar g^{ij}\ddot\Psi^{kl,pq}\bar\nabla_i\rr_{kl}\bar\nabla_j\rr_{pq}-\dot\Psi^{kl}\bar g_{kl}S_1-n(1-\alpha)\Psi.
$$
Since $\Psi$ is negative and concave and $\alpha\leq 1$, all of the terms on the right-hand side are non-positive at a maximum point of $S_1$.  We discard all but the second-last, and use the expression for $\Psi$ in terms of $f_*$ to obtain
$$
\frac{\partial S_1}{\partial t}\leq \overline{\mathcal L}S_1-\alpha S^{-(1+1/\alpha)}\dot F_*^{kl}\bar g_{kl}S_1.
$$
The key term to understand is $\dot F_*^{kl}\bar g_{kl}=\sum_{i=1}^2\dot f_*^i$.  In the two-dimensional case we observe that since $f_*$ is concave, we can write $\dot F_*^{kl}\bar g_{kl}=Q(S_1/f_*)$, where $Q$ is a non-decreasing function.  Thus we have
$$
\dot F_*^{kl}\bar g_{kl} = Q(S_1/f_*)=Q(S^{1/\alpha}S_1)\geq Q(S_-^{1/\alpha}S_1).
$$
Therefore by the maximum principle we have 
$$
S_1\leq S_-^{-1/\alpha}G^{-1}\left(\alpha S_+^{-(1+1/\alpha)}(t-t_0)\right)
$$
where $G(x) = \int_x^\infty\frac{dz}{zQ(z)}$, provided the function $\frac{1}{zQ(z)}$ is integrable.  We relate the latter condition to the integral condition in the theorem, by making the change of variables $z=\frac{1+x}{h(x)}$, where $h(x)=f_*(x,1)$.  We observe that $\lim_{x\to 0}\frac{h(x)}{h'(x)}=0$.  A direct computation at the point $(x,1)$ gives $Q(z) = (1-x)h'(x)+h(x)=h'(x)\left(1-x+\frac{h(x)}{h'(x)}\right)$.   We also compute $dz = -\frac{h'}{h^2}\left(1+x-\frac{h}{h'}\right)dx$, and hence
\begin{align*}
G\left(\frac{1+x}{h(x)}\right) &= \int_z^\infty\frac{dz}{zQ(z)}\\
&=\int_0^x\frac{h'}{h^2}\frac{1+x-\frac{h}{h'}}{\frac{1+x}{h}h'\left(1-x+\frac{h}{h'}\right)}dx\\
&\simeq\int_0^x\frac{dx}{h(x)}.
\end{align*}
The assumption in the theorem is that $\int_0^x\frac{dx}{h(x)}<\infty$, and it follows that we have an upper bound on $S_1$ for each positive time.

Since $f_*=0$ on the boundary of $\Gamma_+$, an upper bound on $S_1$ and a lower bound on $f_*$ imply a lower bound on principal radii of curvature.  Therefore the principal radii of curvature remain in a compact subset of $\Gamma_+$ on which the equation is uniformly parabolic, and uniform higher derivative estimates follow from \cite{Andrews2D} and Schauder estimates as in the previous case.
\qed\end{proof}

Note that in the case $n=2$ the inverse-concavity of $f$ and $f_*$ on the boundary of the positive cone are automatic, since they are linear on the boundary line.  Note that the condition in the second case of the theorem is always satisfied if $f(x,1)$ is bounded by a positive power of $x$, as it is in most classical examples.  An example which does not satisfy the condition but still has $f=f_*=0$ on the boundary of $\Gamma_+$ is given by $f(x,y)=\frac{x+y}{\log(x+y)-\log x-\log y}$.  In examples such as this the equation becomes singularly parabolic near points with $\kappa_{\min}=0$, so one might expect that a Harnack principle or strong maximum principle might imply strict convexity.  Because of this we conjecture that the result should still hold if the integral condition in case \ref{case:f0} is removed.  It also seems likely that the assumption that $f_*$ is concave could be removed.   In this case the result would be essentially sharp in view of the counterexamples constructed previously.

Now we proceed to the higher-dimensional case.  We begin with a class of flows in which the speed is comparable to a power of the mean curvature, analogous to the first case in Theorem \ref{thm:n=2becomesmooth}.

\begin{theorem}\label{thm:n>2likeMCF}
Suppose $0<\alpha\leq 1$, $f$ is concave or convex, and condition \ref{cond:extendf} holds.  Also assume that $f$ is inverse-concave on the boundary of the positive cone (i.e. condition \ref{cond:invconvbdry} holds). Then for any open bounded convex region $\Omega$, the unique barrier solution $\{\Omega_t\}_{0\leq t<T}$ with initial data $\Omega$ is smooth and strictly convex for $0<t<T$.
\end{theorem}

\begin{proof}
The existence of a unique barrier solution is provided by Theorem \ref{existbarrier}.
The remainder of the proof is essentially the same as the first case of the previous Theorem, except that the estimates of Krylov \cite{Krylov} instead of those of \cite{Andrews2D} are used to provide H\"older continuity of second derivatives.
\qed
\end{proof}

Another class includes examples such as powers of the Gauss curvature.  The argument we use here is similar to that
used for Gauss curvature flows in \cite{A2}, and depends on the diffusion becoming fast in directions of small curvature.

\begin{theorem}\label{thm:n>2likeGCF}
Suppose that $0<\alpha\leq 1$, $f$ and $f_*$ both vanish on the boundary of the positive cone, and $f_*$ is concave.  Define $\sigma(r) = \inf\{\dot F_*\big|_A(I):\ \rr_{\max}(A)\geq (r),\ F_*(A)=1\}$, and suppose that $\int_1^\infty\frac{dr}{r\sigma(r)}<\infty$.  Then for any bounded open convex region $\Omega_0$, the unique barrier solution $\{\Omega_t\}_{0\leq t<T}$ with initial data $\Omega$ is smooth and uniformly convex for $0<t<T$.
\end{theorem}

\begin{proof}
As in the previous arguments it suffices to work on a time interval where we have uniform upper and lower bounds on the speed:  $S_-\leq S(x,t)\leq S_+$ for some positive constants $S_-$ and $S_+$.  
We key estimate comes from the evolution of $\rr_{ij}$ from Lemma \ref{T:evlnr}, which implies
$$
\frac{\partial}{\partial t}\rr_{ij}\leq \overline{\mathcal L}\rr_{ij}-\alpha S_+^{-(1+1/\alpha)}\dot F_*^{kl}\bar g_{kl}\rr_{ij}.
$$
It follows from the definition of $\sigma$ that
$$
\dot F_*\big|_{\rr}(I)=\dot F_*\big|_{\rr/F_*(\rr)}(I)\geq \sigma\left(\rr_{\max}S_-^{1/\alpha}\right),
$$
since the largest eigenvalue of $\rr/F_*(\rr)$ is $\rr_{\max}S^{1/\alpha}\geq S_-^{1/\alpha}\rr_{\max}$.  By the maximum principle we therefore have
$$
\rr_{\max}(p,t)\leq S_-^{-1/\alpha}G^{-1}\left(\alpha S_+^{-(1+1/\alpha)}(t-t_0)\right),
$$
where $G(z) = \int_z^\infty\frac{du}{u\sigma(u)}$.  This gives an upper bound on $\rr_{\max}$ for positive times under the assumption of the Theorem.  A lower bound then follows since $f_*$ vanishes on the boundary of the positive cone.  The equation is therefore uniformly parabolic, and uniform higher derivative estimates follow.  Since the barrier solution is the limit of the smooth approximating solutions, the same estimates hold for the barrier solution.
\qed
\end{proof}

We provide some examples where Theorem \ref{thm:n>2likeGCF} applies:  Take
$f=\prod_{k=1}^{n}\left(\frac{E_k}{E_{k-1}}\right)^{\alpha_k}$, where $\sum_k\alpha_k=1$, $\alpha_k\geq 0$ for every $k$, and $\alpha_n>0$ and $\alpha_1>0$.  This is a geometric mean of the functions $E_{k+1}/E_k$, each of which is concave and inverse-concave.  Furthermore, $f$ vanishes on the boundary since $\alpha_n>0$, and $f_*=\prod_{k=1}^n\left(\frac{S_{k}}{S_{k-1}}\right)^{\alpha_{n+1-k}}$ vanishes on the boundary of the positive cone since $\alpha_1>0$.  We note that if we order the principal radii of curvature, so that $\rr_1\geq \rr_2\geq \dots\geq \rr_n$, then $S_{k}/S_{k-1}$ is comparable to $\rr_k$, so we have that $S_1/f_*$ is comparable to $\prod_{k=2}^n\left(\frac{\rr_1}{\rr_k}\right)^{\alpha_{n+1-k}}$, which is bounded by $\left(\frac{\rr_1}{\rr_n}\right)^{1-\alpha_n}$.  But we also have (from differentiating the $S_n/S_{n-1}$ terms) that $\dot F_*(I)\geq n\alpha_1\left(\frac{S_n}{S_{n-1}}\right)^{\alpha_1-1}\prod_{k=1}^{n-1}\left(\frac{S_k}{S_{k-1}}\right)^{\alpha_{n+1-k}}$, which is comparable to
$\prod_{k=1}^{n-1}\left(\frac{\rr_k}{\rr_n}\right)^{\alpha_{n+1-k}}$, hence bounded below by a multiple of $(\rr_1/\rr_n)^{\alpha_n}$.  That is, $\dot F_*(I)\geq c\alpha_1(S_1/f_*)^{\frac{\alpha_n}{1-\alpha_n}}$, and the Theorem applies.

Our final class of flows only allows $\alpha=1$, but we conjecture that a similar result should hold with $0<\alpha<1$ also.  This class involves a mixture of the methods of the previous two cases:  First using fast diffusion to bound the principal curvatures above, and then using uniform parabolicity to obtain lower bounds on principal curvatures.

\begin{theorem}\label{thm:n>2intermediate}
Assume $\alpha=1$, $f$ is concave, $f_*=0$ on $\partial\Gamma_+$, and conditions \ref{cond:invconvbdry} and  \ref{cond:extendf} hold.   Define $\sigma(r) = \inf\{\sum_i\dot f^i\kappa_i^2:\ f(\kappa)=1$, $\kappa_{\max}\geq r\}$, and suppose that $\int_1^\infty\frac{dr}{r\sigma(r)}<\infty$.  Then for any bounded open convex region $\Omega_0$, the unique barrier solution $\{\Omega_t\}_{0\leq t<T}$ with initial data $\Omega$ is smooth and uniformly convex for $0<t<T$.
\end{theorem}

\begin{proof}
We use the evolution of the second fundamental form from equation \eqref{eq:dth} in Lemma \ref{T:evlneqns}.  In the case $\alpha=1$ the concavity of $f$ gives
$$
\frac{\partial}{\partial t}h_i^j\leq {\mathcal L}h_i^j+h_i^j\dot F^{kl}h_k^ph_{pl}.
$$
Working as in the previous theorems on a time interval $[t_0,t_1]$ where the speed $S$ has uniform positive upper and lower bounds, so that $S_-\leq S(x,t)\leq S_+$, we compute
$$
\frac{\partial}{\partial t}\left(\frac{h_i^j}{2S-S_-}\right)\leq 
{\mathcal L}\left(\frac{h_i^j}{2S-S_-}\right)+\frac{4\nabla_kS}{2S-S_-}\nabla_k\left(\frac{h_i^j}{2S-S_-}\right)-\frac{S_-h_i^j}{(2S-S_-)^2}\dot F^{kl}h_k^ph_{pl}.
$$
The definition of $\sigma$ gives $\dot F^{kl}h_k^ph_{pl}\geq \sigma\left(\frac{\kappa_{\max}}{S_+}\right)\geq \sigma\left(\frac{S_-\kappa_{\max}}{S_+(2S-S_-)}\right)$, and it then follows form the maximum principle that
$$
\frac{\kappa_{\max}}{2S-S_-}\leq\frac{S_+}{S_-}G^{-1}\left(\frac{S_-(t-t_0)}{2S_+}\right),
$$
where $G(z)=\int_z^\infty \frac{dr}{r\sigma(r)}$.  This gives an upper bound on $\kappa_{\max}$ for each $t>t_0$.  Condition \ref{cond:extendf} then implies uniform bounds above and below on $\dot f^i$.  Therefore the flow is uniformly parabolic, and we proceed as in Theorem \ref{thm:n>2likeMCF}.
\qed
\end{proof}

Perhaps the most important examples where Theorem \ref{thm:n>2intermediate} applies are the following (analogous to the examples above for Theorem \ref{thm:n>2likeGCF}):  Let 
$F = \prod_{k=1}^n\left(\frac{E_k}{E_{k-1}}\right)^{\alpha_k}$, where $0<\alpha_1<1$, $\alpha_k\geq 0$ for all $k$, and $\sum_{i=1}^n\alpha_1=1$.  Important special cases include $\alpha_1=\dots=\alpha_k=\frac{1}{k}$, which yields $F=E_k^{1/k}$.   In these cases we have $\dot F(h^2)\geq \alpha_1\kappa_{\max}F$, so $\sigma(r)\geq \alpha_1r$, and the condition of the Theorem is satisfied.

We conclude this section and the paper by stating a conjecture, which would amount to a set of conditions which are close to being necessary and sufficient:

\begin{conjecture}
Let $0<\alpha\leq 1$.
Suppose that $f_*$ vanishes on the boundary of the positive cone, and that $f$ extends smoothly to $\partial\Gamma_+$ near points where it is non-zero in such a way that the derivative of $f$ normal to the boundary is positive, and $f$ is inverse-concave on the boundary of $\Gamma_+$ (that is, condition \ref{cond:invconvbdry} holds).  Also assume that condition \ref{cond:admit.reg.extend} holds).  Then for any initial data $M_0$ given by the boundary of an open bounded convex region, there exists a unique barrier solution $\{M_t\}_{0\leq t<T}$ of \eqref{E:theflow}.  The hypersurfaces $M_t$ are smooth and strictly convex for $0<t<T$, and converge to a point as $t\to T$.
\end{conjecture}

The assumption on H\"older estimates is automatically satisfied if $f$ is concave, for example, and is automatic in the two-dimensional case by the results of \cite{Andrews2D}.  The normal derivative condition rules out the shrinking cylinder examples of section \ref{S:cylinders}, and the persisting ridges of section \ref{sec:nonsmooth} are ruled out since $f_*$ vanishes on the boundary of the positive cone.  Progress towards this conjecture would seem something similar to a version of the strong maximum principle of Bony \cite{BonyMP} with weak regularity assumptions.

\end{document}